\documentclass[11]{article} 

\usepackage[totalwidth=480pt, totalheight=680pt, hmargin=2.55cm, vmargin=2.15cm]{geometry}   

\usepackage[american]{babel}

\usepackage{url}                                

\usepackage{abstract}

\usepackage{appendix}

\usepackage[nottoc]{tocbibind}

\input{localstyle.sty}

\theoremstyle{break}
\theoremheaderfont{\normalfont\bfseries}
\theorembodyfont{\itshape}
\theoremseparator{:}
\theoremnumbering{arabic}
\theoremsymbol{\$}

\newtheorem{theorem}{Theorem}[section]
\newtheorem{proposition}{Proposition}[section]

\newtheorem{lemma}{Lemma}[section]

\theorembodyfont{\upshape}
\theoremstyle{nonumberplain}
\newtheorem{proof}{Proof}

\theoremsymbol{}
\newtheorem{remark}{Remark}

\DeclareMathOperator{\cardinal}{card}
\DeclareMathOperator{\GaloisField}{GF}
\DeclareMathOperator{\GaloisTrace}{tr}
\DeclareMathOperator{\PGCD}{PGCD}

\DeclareMathOperator{\NeperBase}{e}
\DeclareMathOperator{\Module}{Mod}
\DeclareMathOperator{\modulus}{mod}
\DeclareMathOperator{\CayleyAlgebra}{CAlg}
\DeclareMathOperator{\UniferousAlgebra}{UAlg}

\newcommand{\CayleyTrace}[1]{\rsfsT_{#1}}
\newcommand{\CayleyNorm}[1]{\rsfsN_{#1}}

\title{Cayley Integers\\
(long version) - math.RA/0506349}
\author{Hubert Holin \\
\url{Hubert.Holin.1982@Polytechnique.org}}
\date{05/05/2005}

\begin{document}


\maketitle

\begin{abstract}
We present here some results of applying the Cayley-Dickson process to certain alternative
algebras (notably built upon Galois fields and congruence rings), in a manner which might
yield new building blocks for cryptographic systems. The results presented here are technical
but not inherently difficult.
\end{abstract}

\begin{center}
\begin{minipage}[t]{12cm}
\footnotesize
\tableofcontents
\end{minipage}
\end{center}
\normalsize

\pagebreak
\section*{Introduction}
\addcontentsline{toc}{section}{\protect\numberline{}{Introduction}}
\indent

The idea of applying the Cayley-Dickson process
to Galois fields is not new (\cite{Schafer(1945), Jacobson(1958), SpringerVeldkamp(2000)}).
Finding the number of units and unimodular elements of the resulting structure, and of
iterations of the doubling procedure on the base Galois field, is a simple application
of a classical exercise (\cite{IrelandRosen(1982)}), as we show here, but the relevance
of the diophantine equations which are classically presented, to the algebraic
structures we are interested in, is not readily apparent in the published literature.
Similar results, but using congruence rings as starting points rather than fields
(for possible uses in cryptographic systems), are
presented here, for what we believe is the first time. One key ingredient of this extension
is \thref{thm:decdbl}, which we believe is new. We focus on enumeration properties
rather than the classification and comparison questions which are extensively studied
elsewhere, at least for algebras over fields (\cite{Schafer(1966),
ZhevlakovSlin'koShestakovShirshov(1982), GoodaireJespersMiles(1996),
SpringerVeldkamp(2000)}).

We would like to emphasize the fact that what we present here is at the confluence of
several rather distinct mathematical traditions, whose terminologies and
notations are sometimes incompatible. We have therefore been led to
make choices, and to be sometimes rather verbose so as to disambiguate the notion used.
It would be worth one's while to take a look at the appendices, if only to ascertain
the meaning we have chosen here for some terms.

\section{Abstract Cayley objects}
\subsection{Definitions and elementary properties}
\subsubsection{Order and index}\label{ssc:ordind}
\indent

Let \((S,\times)\) be a di-associative I.P. loop, whose neutral element we will
denote by \(e\). For \(\omega \in S\) we can define the set
\(P_{\omega} = \{\omega^{n}\thinspace|\thinspace n \in \bbZ\}\).

\(P_{\omega}\) is an abelian group and \([\bbZ \rightarrow P_{\omega};
n \mapsto \omega^{n}]\) is a group homomorphism. We can therefore define \(\omega\)
to be of \emph{infinite order} if that homomorphism is injective, and to be of \emph{finite
order} \(\cardinal(P_{\omega})\) otherwise. These definitions are consistent with their
counterparts for groups. Note that if \(\omega\) is of finite order \(\alpha\), then for any
\(\beta \in \bbZ\), \(\omega^{\beta}\) is also of finite order, that order being
\(\alpha/\PGCD(\alpha,\beta)\), as \((P_{\omega},\times)\) is an abelian group
(recal that \(\PGCD(\alpha,\beta)\) designates the least common multiple of \(\alpha\)
and \(\beta\)).

\(P_{\omega}\) also is a subloop of \((S,\times)\), and we define the two following
binary relations on \(S\):
\begin{gather*}
(\forall\thinspace (x,y)\in S^{2})\thickspace
[\thinspace x \thinspace{\cal{P}}_{\omega}^{R}\thinspace y \Leftrightarrow
x \times y^{-1} \in P_{\omega} \Leftrightarrow
(\exists\thinspace n \in \bbZ)\thickspace x = \omega^{n} \times y \thinspace] \\
(\forall\thinspace (x,y)\in S^{2})\thickspace
[\thinspace x \thinspace{\cal{P}}_{\omega}^{L}\thinspace y \Leftrightarrow
x^{-1} \times y \in P_{\omega} \Leftrightarrow
(\exists\thinspace n \in \bbZ)\thickspace y = \omega^{n} \times x \thinspace]
\end{gather*}

\begin{proposition}\label{prp:powreleq}
\({\cal{P}}_{\omega}^{R}\) and \({\cal{P}}_{\omega}^{L}\) are
equivalency relations on \(S\), the inversion is a bijection from
\(S/{\cal{P}}_{\omega}^{R}\) to \(S/{\cal{P}}_{\omega}^{L}\) and
any element of \(S/{\cal{P}}_{\omega}^{R}\) is equipotent to \(P_{\omega}\).
\end{proposition}

\begin{proof}
We only indicate the proof that \({\cal{P}}_{\omega}^{R}\) is an equivalency relation, as
the one for \({\cal{P}}_{\omega}^{L}\) is perfectly similar.

\((\forall\thinspace x \in S)\thickspace x \times x^{-1} = x^{-1} \times x =
e = \omega^{0} \in P_{\omega}\) and hence \({\cal{P}}_{\omega}^{R}\) is
reflexive.

Let \(x \in S\) and \(y \in S\) such that
\(x \thinspace{\cal{P}}_{\omega}^{R}\thinspace y\). Then
\((\exists\thinspace n \in \bbZ)\thickspace x = \omega^{n} \times y\), so
\(y = \omega^{-n} \times x\), since we are in an I.P. loop (and \(\omega^{-n} =
(\omega^{n})^{-1}\) because of the loop di-associativity). Thus
\({\cal{P}}_{\omega}^{R}\) is symmetrical.

Let \(x \in S\), \(y \in S\) and \(z \in S\) such that
\(x \thinspace{\cal{P}}_{\omega}^{R}\thinspace y\) and
\(y \thinspace{\cal{P}}_{\omega}^{R}\thinspace z\). Then
\((\exists\thinspace n \in \bbZ)\thickspace x = \omega^{n} \times y\) and
\((\exists\thinspace m \in \bbZ)\thickspace y = \omega^{m} \times z\), but then
\(x = \omega^{n} \times (\omega^{m} \times z) = \omega^{n+m} \times z\) thanks
to the di-associativity. Thus
\({\cal{P}}_{\omega}^{R}\) is transitive.

This concludes the proof that  \({\cal{P}}_{\omega}^{R}\) is an equivalency relation;

Let \(x \in S\) and \(y \in S\) such that
\(x \thinspace{\cal{P}}_{\omega}^{R}\thinspace y\). Then
\((\exists\thinspace n \in \bbZ)\thickspace x \times y^{-1} = \omega^{n}\), so
\((y^{-1})^{-1} \times (x^{-1}) = y \times x^{-1} = \omega^{-n}\), which means that
\(y^{-1} \thinspace{\cal{P}}_{\omega}^{L}\thinspace x^{-1}\), and the inversion
is a bijection from \(S/{\cal{P}}_{\omega}^{R}\) to \(S/{\cal{P}}_{\omega}^{L}\).

Let now \(Y \in S/{\cal{P}}_{\omega}^{R}\) and \(y \in Y\), and consider
\(\phi_{y}: P_{\omega}\rightarrow S; z \mapsto z \times y\). Since we are in a loop,
\(\phi_{y}\) is injective, and since the loop is di-associative, \(\phi_{y}\) takes its values
in \(Y\). Furthermore, if \(x \in Y\), then \(x = \omega^{n_{x}} \times y\) for some
\(n_{x} \in \bbZ\), so \(x \times y^{-1} = \omega^{n_{x}}\), thanks to the I.P. of the
loop. But then \(\phi_{y}(\omega^{n_{x}}) = \phi_{y}(x \times y^{-1}) =
(x \times y^{-1}) \times y = y\), again because we are I.P., and thus \(\phi_{y}\) is
surjective. Finally, \(\phi_{y}\) is bijective.
\end{proof}

\begin{remark}
There is no reason for the left and right cosets to agree, in particular there is no reason for
\(P_{\omega}\) to be a normal subloop of \(S\).
\end{remark}

If either \(S/{\cal{P}}_{\omega}^{R}\) or \(S/{\cal{P}}_{\omega}^{L}\) is finite,
then the other one is too, and they have the same cardinal; in that case we will call that
cardinal the \emph{index} of \(P_{\omega}\) in \(S\), and denote it by
\([S:P_{\omega}]\). This definition is also consistent with its counterpart for groups.

\subsubsection{Conjugation}
\indent

Recall (\cite{Bourbaki(A)}) that, given a commutative and associative ring with unit
\((A,+,.)\) and \((E,+,\times,\cdot)\) an \(A\)-algebra with unit, with neutral element
\(e\), a \emph{conjugation} over \(E\) is any function \(\sigma : E \rightarrow E\)
which is bijective, \(A\)-linear and such that:
\begin{enumerate}
\item \(\sigma(e) = e\).
\item \((\forall\thinspace (x,y) \in E^{2})\thickspace
\sigma(x \times y) = \sigma(y) \times \sigma(x)\) (\emph{beware} the interversion of
\(x\) and \(y\)!).
\item \((\forall\thinspace x \in E)\thickspace (x+\sigma(x)) \in A \cdot e\).
\item \((\forall\thinspace x \in E)\thickspace (x \times \sigma(x)) \in A \cdot e\).
\end{enumerate}

These properties imply\footnote{\((x+\sigma(x)) \in A \cdot e \Rightarrow
x \times \sigma(x) = x \times (x + \sigma(x)) - x \times x =
(x + \sigma(x)) \times x - x \times x = \sigma(x) \times x\).}
\((\forall\thinspace x \in E)\thickspace x \times \sigma(x) = \sigma(x) \times x\)
and\footnote{Given \(x \in E\), there exists \(\alpha \in A\) such that
\(x + \sigma(x) = \alpha \cdot e\); the \(A\)-linearity of \(\sigma\) then implies
\(\sigma(x)+(\sigma \circ \sigma)(x) = \sigma(x + \sigma(x)) =
\alpha \cdot \sigma(e)\), and finally, \(\sigma(e) = e\).}
\((\forall\thinspace x \in E)\thickspace (\sigma \circ \sigma)(x) = x\).

Note that if \(E = A\), the identity function is always a conjugation, but that otherwise there
is no guaranty a conjugation can exist. Note also that while, for a given
\(x \in E\), there exists \(T \in A\) and \(N \in A\) such that
\(\CayleyTrace{E}(x) = T \cdot e\) and \(\CayleyNorm{E}(x) = N \cdot e\),
\(T\) and \(N\) may happen not to be unique (see \ref{sbs:manylaws}). 
We will also write
\(\bar{x}\) for \(\sigma(x)\) when no confusion is to be feared.

\subsubsection{Cayley Algebra}
\indent

Recall (\cite{Bourbaki(A)}) that, given a commutative and associative ring with unit
\((A,+,.)\) a \emph{Cayley algebra} over \(A\) is a structure
\((E,+,\times,\cdot,\sigma)\), where \((E,+,\times,\cdot)\) is an \(A\)-algebra with
unit, and \(\sigma\) is a conjugation over \(E\).

Given \((E,+,\times,\cdot,\sigma)\) an \(A\)-Cayley algebra and \(T \subset E\),
we will say \(T\) is a \emph{sub-Cayley algebra} if and only if it is a sub-algebra of an
algebra with unit, of \(E\), and it is stable by conjugation (\emph{i.e.}
\((\forall\thinspace x \in T)\thickspace \sigma(x) \in T\)).

On a Cayley algebra, it is convenient to consider the \emph{Cayley trace} and
\emph{Cayley norm}\footnote{Yet another unfortunate collision of terms, as this ``norm''
is actually quadratic...} defined respectively by \(\CayleyTrace{E}(x) = x + \sigma(x)\)
and \(\CayleyNorm{E}(x) = x \times \sigma(x)\).

We have the following important relations valid for all elements:
\begin{gather}
\CayleyTrace{E}(\sigma(x)) = \CayleyTrace{E}(x) \\
\CayleyNorm{E}(\sigma(x)) = \CayleyNorm{E}(x) \\
\CayleyTrace{E}(y \times x) = \CayleyTrace{E}(x \times y) \label{eqn:comtrace}\\
\CayleyTrace{E}(x \times \sigma(y)) = \CayleyTrace{E}(y \times \sigma(x)) =
\CayleyTrace{E}(x) \times \CayleyTrace{E}(y)-\CayleyTrace{E}(x \times y) =
\CayleyNorm{E}(x+y)-\CayleyNorm{E}(x)-\CayleyNorm{E}(y) \\
x^2 = \CayleyTrace{E}(x) \times x - \CayleyNorm{E}(x) \label{eqn:square}
\end{gather}

Relation (\ref{eqn:comtrace}) is especially noteworthy in light of the fact that no
commutativity or associativity has been assumed\footnote{It is established by noticing that
\(\CayleyTrace{E}(x \times y) = x \times y +
(\CayleyTrace{E}(y)-y) \times (\CayleyTrace{E}(x)-x)\) and that therefore\linebreak
\(\CayleyTrace{E}(x)\times y + \CayleyTrace{E}(y) \times x +
(\CayleyTrace{E}(x \times y)-\CayleyTrace{E}(x) \times \CayleyTrace{E}(y)) =
x \times y + y \times x =
\CayleyTrace{E}(y)\times x + \CayleyTrace{E}(x) \times y +
(\CayleyTrace{E}(y \times x)-\CayleyTrace{E}(y) \times \CayleyTrace{E}(x))\).}.
By contrast, in absence of any form of associativity (or commutativity), little more can be
proved about the Cayley norm\footnote{See \ref{sec:fourtout} for a discussion of what
one might want out of a ``norm''. It should be noted that
\((\forall\thinspace (x,y)\in E^{2})\thickspace
\CayleyNorm{E}(x\times y) = \CayleyNorm{E}(y\times x)\)
may happen even when \(E\) is not alternative (where this is always true; see
\thref{prp:altcay}), for instance in the case of the hexadecimalions.}, apart from the
following proposition (where \(\Module_{A}(x,y)\) denotes the
\(A\)-module generated by \(x \in E\) and \(y \in E\), and \(e\) denotes the neutral
element of \(E\)):

\begin{proposition}\label{lem:cmplx}
Given \(x \in E\), \(\Module_{A}(e,x)\) is stable for ``\(\times\)''; it is a sub-Cayley
algebra of \(E\) which is both associative and commutative.
Furthermore \((\forall\thinspace m \in \bbN)\thickspace
\CayleyNorm{E}(x^{n}) = (\CayleyNorm{E}(x))^{n}\).
\end{proposition}

\begin{proof}
This is a direct application of (\ref{eqn:square}), which we rewrite as
\(x^2 = T \cdot x - N \cdot e\), with \(T \in A\) such that
\(\CayleyTrace{E}(x) = T \cdot e\) and \(N \in A\) such that
\(\CayleyNorm{E}(x) = N \cdot e\).

The last statement of the proposition is proved by induction.
\end{proof}

\begin{remark}
Given \(x \in E\), \((\forall\thinspace n \in \bbN)\thickspace
x^{n} \in \Module_{A}(e,x)\), and we have:
\begin{equation*}
(\forall\thinspace n \in \bbN)
(\forall\thinspace m \in \bbN)\thickspace
x^{n+m} = x^n \times x^m = x^m \times x^n\label{eqn:powN}
\end{equation*}
which is just property (\ref{eqn:AdditiveN}) of \ref{ssc:pitpow} under different notations.
\end{remark}

Invertible elements have a few more properties:

\begin{proposition}\label{prp:xandinverse}
If \(x \in E\) has an inverse for ``\(\times\)'', then \(\Module_{A}(e,x,x^{-1})\) is
stable for ``\(\times\)''; it is a sub-Cayley algebra of \(E\) which is both associative and
commutative.
\end{proposition}

\begin{proof}
Yet another easy consequence of (\ref{eqn:square}), applied to \(x\) as well as to
\(x^{-1}\), and the definition of the inverse.
\end{proof}

\begin{remark}
With no assumption of associativity for \(E\), not much can be said about the relationship
between the Cayley norms of \(x\) and of \(x^{-1}\) in \thref{prp:xandinverse}.
However, if we know that \(\CayleyNorm{E}(x)\) is invertible in \(A \cdot e\), then
of course \(x\) is invertible in \(E\) and
\(x^{-1} = {\CayleyNorm{E}(x)}^{-1} \cdot \sigma(x)\).
\end{remark}

With some rather weak hypotheses, however, we can improve the situation, as is well known
(\cite{Bourbaki(A)}):

\begin{proposition}\label{prp:altcay}
Let \((A,+,.)\) be a commutative and associative ring with unit, and
\((E,+,\times,\cdot,\sigma)\) an \emph{alternative} Cayley algebra over \(A\), with
neutral element \(e\). An element \(x \in E\) is invertible if and only if
\(\CayleyNorm{E}(x)\) is invertible in \(A \cdot e\), in which case the inverse of \(x\)
is \({\CayleyNorm{E}(x)}^{-1} \cdot \sigma(x)\), and \(\Module_{A}(e,x)\) is an
abelian group. Furthermore, \((E,\times)\) is di-associative and we have
\((\forall\thinspace x \in E)(\forall\thinspace y \in E)\thickspace
\CayleyNorm{E}(x \times y) = \CayleyNorm{E}(x) \times \CayleyNorm{E}(y)\).
\end{proposition}

\begin{remark}
The di-associativity of \((E,\times)\) as a magma (see \ref{ssc:altrng}) ensures that
\begin{equation*}
(\forall\thinspace x \in E)(\forall\thinspace y \in E)
(\forall\thinspace n \in \bbN)
(\forall\thinspace m \in \bbN)(\forall\thinspace p \in \bbN)
(\forall\thinspace q \in \bbN) \thickspace
x^{n+m} \times y^{p+q} = x^m \times (x^n \times y^{p+q}) =
(x^{n+m} \times y^p) \times y^q)\label{eqn:powDN}
\end{equation*}
which is just (\ref{eqn:AdditiveDN}) of \ref{ssc:pitpow} under different
notations.
\end{remark}

Alternative Cayley algebras also enjoy a few other nice properties\footnote{This is a 
special case of the theorem which states that the invertible elements of an alternative ring with
unit are a Moufang loop (\cite{GoodaireJespersMiles(1996)}); we only give a proof here
because it is slightly simpler than the general case.}.

\begin{proposition}\label{prp:invloop}
Let \((A,+,.)\) be a commutative and associative ring with unit, and
\((E,+,\times,\cdot,\sigma)\) an \emph{alternative} Cayley algebra over \(A\), with
neutral element \(e\). Let \(E^{*}\) the elements of \(E\) which have an inverse for
``\(\times\)'', then \((E^{*},\times)\) is a Moufang loop.
\end{proposition}

\begin{proof}
An alternative algebra being di-associative, \((E^{*},\times)\) is a di-associative
magma with unit. By definition, every \(x \in E^{*}\) has an inverse \(x^{-1}\), and
since \(E\) is an alternative Cayley algebra,
\(x^{-1} = {\CayleyNorm{E}(x)}^{-1} \cdot \sigma(x)\) and
\((\exists\thinspace T \in A)\thickspace \sigma(x) = T \cdot e - x\), so
\((\exists\thinspace (\alpha,\beta) \in A^{2})\thickspace x^{-1} =
\alpha\cdot e + \beta\cdot x\). Therefore, \(x^{-1}\times(x \times y) =
(\alpha\cdot e + \beta\cdot x)\times(x \times y) =
(\alpha\cdot x)\times y + \beta\cdot(x\times(x\times y))\), but since the algebra
is alternative, \(x\times(x\times y) = (x \times x)\times y\), and we finaly see that
\(x^{-1}\times(x \times y) = (x^{-1}\times x) \times y = y\). We likewise prove
that \((y \times x)\times x^{-1} = y\). This proves that \(E^{*}\) is a quasigroup, hence
a loop since it has a neutral element, and has the I.P. property.

Using \thref{prp:altrngmouf} we find that the loop is actually a Moufang loop.
\end{proof}

\begin{remark}
If \(E\) is actually associative, then of course \(E^{*}\) is a group.
\end{remark}

\begin{remark}
As a consequence of the di-associativity of loops, we have
\begin{equation*}(\forall\thinspace x \in E^{*})(\forall\thinspace y \in E^{*})
(\forall\thinspace n \in \bbZ)
(\forall\thinspace m \in \bbZ)(\forall\thinspace p \in \bbZ)
(\forall\thinspace q \in \bbZ) \thickspace
x^{n+m} \times y^{p+q} = x^m \times (x^n \times y^{p+q}) =
(x^{n+m} \times y^p) \times y^q)\label{eqn:powDZ}
\end{equation*}
which is just (\ref{eqn:AdditiveDZ}) of \ref{ssc:pitpow} under different notations.
\end{remark}

\begin{remark}
The material developed in \ref{ssc:ordind} applies to \(E^{*}\). In particular, if
\(\cardinal(E^{*})\) is finite we have
\allowbreak\((\forall\thinspace \omega \in E^{*})\thickspace
\cardinal(E^{*}) = [E^{*}:P_{\omega}] \cardinal(P_{\omega})\).
\end{remark}

\subsubsection{Cayley Morphism}\label{ssc:caymor}
\indent

Given two Cayley algebras \((E,+,\times,\cdot,\sigma)\) and
\((F,\dagger,*,\bullet,s)\) on a same commutative and associative ring with unit
\((A,+,.)\), a function \(\phi : E \rightarrow F\) is a \emph{(homo)morphism of
Cayley algebras} (or simply a \emph{Cayley (homo)morphism}) if and only if it is a
(homo)morphism of algebras with unit and that furthermore\linebreak
\((\forall\thinspace x \in E)\thickspace \phi(\sigma(x)) = s(\phi(x))\).
Monomorphism, epimorphisms, isomorphisms, endomorphisms and automorphisms of
Cayley algebras are built from homomorphism of Cayley algebras in the usual manner.

As a matter of practicality, it is useful to note that a Cayley morphism which is an
isomorphism of algebras with unit, actually is a Cayley isomorphism, as is trivially seen.

One verifies immediately that one has a category (as \emph{per} \cite{Lang(1971)})
(which we will call \(\CayleyAlgebra\) here) by considering Cayley algebras as objects,
Cayley morphisms as arrows, and the usual composition of functions as composition of
arrows. Of course, if we denote by \(\UniferousAlgebra\) the category of algebras with
unit and corresponding classical morphisms of algebras with units, then associating every
Cayley algebra with its underlying algebra with unit, and every Cayley algebra morphism
with its underlying morphism of algebra with unit, yields a covariant functor from 
\(\CayleyAlgebra\) to \(\UniferousAlgebra\) (a stripping functor).

\subsection{The Cayley-Dickson process}
\subsubsection{A fundamental shortcoming}\label{ssc:funshrt}
\indent

A major problem with \(\CayleyAlgebra\), as categories go, is that products may fail to
exist, if we want the product to be compatible with the stripping functor of \ref{ssc:caymor}.

Consider the following Cayley algebra over \(\bbR\): 
\((E,+,\times,\cdot,\sigma) = (\bbR,\times,+,\times, id_{\bbR})\), with
\(id_{\bbR}\) the identity function on \(\bbR\). Assume we have built the product of
\(E\) with itself, yielding \((G,\boxplus,\boxtimes,\boxdot,\circledS)\). Since we
want the product creation to be compatible with the stripping functor, we see at once that
\((G,\boxplus,\boxtimes,\boxdot)\) must be the product algebra of
\((E,+,\times,\cdot)\) with itself, and that the Cayley algebra
morphisms from from \(G\) to \(E\) (as first component) and \(E\) (as second
component) must be the projections of the algebra \(G\) onto \(E\). Hence
\((x,y) \boxplus (x',y') = (x+x',y+y')\), \((x,y) \boxtimes (x',y') = (xx',yy')\)
(with \((1,1)\) being the neutral element for ``\(\boxtimes\)''),
\(\lambda \boxdot (x,y) = (\lambda x,\lambda y)\) and \(\circledS((x,y)) = (x, y)\).
But then \((1,2) \boxtimes \circledS ((1,2)) = (1,4) \not\in \bbR \boxdot (1,1)\).

There are interesting cases, however, where products \emph{do} exist, which are compatible
with the stripping functor:

\begin{theorem}[Dominated Product]\label{thm:domprod}
Let \(I\) be a set, \(((E_{i},+_{i},\times_{i},\cdot_{i},\sigma_{i}))_{i \in I}\) a
family of Cayley algebras on a same commutative and associative ring with unit
\((A,+,.)\), and \((E,+,\times,\cdot)\) the product algebra of that family, seen
as mere algebras with unit. If there exists a Cayley algebra \((F,\dagger,*,\bullet,s)\)
and a family of Cayley morphisms \((\phi_{i})_{i \in I}\), with
\(\phi_{i} : F \rightarrow E_{i}\) for all \(i \in I\), such that the morphism of
algebras with unit \(\Phi : F \rightarrow E; y \mapsto (\phi_{i}(y))_{i \in I}\)
is \emph{surjective}, then the function
\(\sigma : E \rightarrow E; (x_{i})_{i \in I} \mapsto (\sigma_{i}(x_{i}))_{i \in I}\)
is a conjugation on \((E,+,\times,\cdot)\),  \((E,+,\times,\cdot,\sigma)\) is a
Cayley algebra, for all \(i \in I\) the projections of \(E\) onto \(E_{i}\) is a Cayley
epimorphism and \(\Phi\) is a Cayley epimorphism.
\end{theorem}

\begin{proof}
It is obvious that \(\sigma\) is an \(A\)-linear involution (and thus bijective). Given
\( (x_{i})_{i \in I} \in E\) and \( (x'_{i})_{i \in I} \in E\),
\(\sigma((x_{i})_{i \in I} \times (x'_{i})_{i \in I}) =
\sigma( (x_{i} \times_{i} x'_{i})_{i \in I}) =
(\sigma_{i}(x_{i} \times_{i} x'_{i}))_{i \in I} =
(\sigma_{i}(x'_{i}) \times_{i} \sigma_{i}(x_{i}))_{i \in I} =
\sigma((x'_{i})_{i\in I}) \times \sigma((x_{i})_{i \in I})\).

Let \(y \in F\). we see that \(\sigma(\Phi(y)) = \sigma((\phi_{i}(y))_{i \in I}) =
(\sigma_{i}(\phi_{i}(y)))_{i \in I}\), but since all the \(\phi_{i}\) are Cayley
morphisms, \((\sigma_{i}(\phi_{i}(y)))_{i \in I} = (\phi_{i}(s(y)))_{i \in I} =
\Phi(s(y))\), and thus \(\sigma(\Phi(y)) = \Phi(s(y))\).

Let \((x_{i})_{i \in I} \in E\); since \(\Phi\) is surjective, there exists \(y \in F\)
such that \(\Phi(y) = (x_{i})_{i \in I}\). We see then that
\((x_{i})_{i \in I} + \sigma((x_{i})_{i \in I}) = \Phi(y) + \sigma(\Phi(y))\), but
we have just proved that  \(\sigma(\Phi(y)) = \Phi(s(y))\), so
\( \Phi(y) + \sigma(\Phi(y)) = \Phi(y) + \Phi(s(y))\), and since \(\Phi\) is
\(A\)-linear, \(\Phi(y) + \Phi(s(y)) = \Phi(y \dagger s(y))\). As 
\((F,\dagger,*,\bullet,s)\) is a Cayley algebra, there exists \(\alpha \in A\) such that,
writing \(\varepsilon\) for the neutral element of \(F\), \(y \dagger s(y) =
\alpha \bullet \varepsilon\). Since \(\Phi\) is \(A\)-linear, we find at once that
\(\Phi(y \dagger s(y)) = \alpha \cdot \Phi(\varepsilon) =
\alpha \cdot (\phi_{i}(\varepsilon))_{i \in I}\)
and, denoting by \(e_{i}\) the neutral element of \(E_{i}\) (for each \(i\)),
\(\alpha \cdot (\phi_{i}(\varepsilon))_{i \in I} = \alpha \cdot (e_{i})_{i \in I}\)
since the \(\phi_{i}\) are all Cayley morphisms. Since \( e = (e_{i})_{i \in I}\) is the
neutral element for \(E\), we have proved that
\((x_{i})_{i \in I} + \sigma((x_{i})_{i \in I})\) is the (external) product the neutral
element of \(E\) by of an element of \(A\). We show in exactly the same way that
\((x_{i})_{i \in I} \times \sigma((x_{i})_{i \in I})\) is the (external) product the neutral
element of \(E\) by of an element of \(A\).

This proves that \(\sigma\) is a conjugation, and thus that \((E,+,\times,\cdot,\sigma)\)
is a Cayley algebra. We have also shown that \(\Phi\), which is a morphism of algebra
with unit, verifies \(\sigma(\Phi(y)) = \Phi(s(y))\), so is a Cayley morphism. Since we
have assumed \(\Phi\) to be surjective, it is a Cayley epimorphism.

The fact that for all \(i \in I\) the projections of \(E\) onto \(E_{i}\) is a Cayley
epimorphism is trivial.
\end{proof}

\subsubsection{Cayley-Dickson for algebras}\label{ssc:caydial}
\indent

This well-known procedure is a way around the shortcoming described in \ref{ssc:funshrt}.
Notably, there is no compatibility with the stripping functor of \ref{ssc:caymor} (though
there \emph{is} compatibility if we strip all the way to the module), and only the product
of a Cayley algebra with itself is defined. We present here a version wich differ only in
notation from that of \cite{Bourbaki(A)}. Note that for our intended uses, we must consider
algebras over rings which may happen not to be fields (in contrast to, \emph{e.g.}
\cite{GoodaireJespersMiles(1996)}).

Let \((A,+,.)\) be a commutative and associative ring with unit (with neutral element
\(1\)), and \((E,+,\times,\cdot,\sigma)\) a Cayley algebra over \(A\) (with neutral
elements \(e\) and \(0\) for ``\(\times\)''  and ``\(+\)'' respectively). Let
\(F = E \vartimes E\), \(\varepsilon = (e,0) \in F\), and
\emph{choose} an element \(\zeta \in A\). Define now
\begin{equation*}
\begin{aligned}
\dagger\negthinspace : & &  F \vartimes F & & \rightarrow & & & F \\
 & & ((x,y),(x',y')) & &\mapsto & & & (x+x',y+y') \\
 *\negthinspace : & &  F \vartimes F & & \rightarrow & & & F \\
  & & ((x,y),(x',y')) & & \mapsto & & & (x\times x' - \zeta \cdot (\sigma(y')\times y),y \times \sigma(x') + y'\times x) \\
  \bullet\negthinspace : & &  A \vartimes F & & \rightarrow & & & F \\
  & & (\lambda,(x,y)) & & \mapsto & & & (\lambda \cdot x, \lambda \cdot y) \\
  s\negthinspace : & &  F & & \rightarrow & & & F \\
  & & (x,y) & & \mapsto & & & (\sigma(x),-y)
\end{aligned}
\end{equation*}

The following proposition is entirely classical (\cite{Bourbaki(A)} and \cite{Klein(2002)})

\begin{proposition}[Structure]\label{thm:struct}
\((F,\dagger,*,\bullet,s)\) is a Cayley algebra over \(A\), whose neutral element is
\(\varepsilon\); \(F\) is commutative if and only if \(E\) is commutative and
\(\sigma\) is the identity on \(E\); \(F\) is associative if and only \(E\) is both
associative and commutative; \(F\) is alternative if and only if \(E\) is associative.

Furthermore \(\CayleyTrace{F}((x,y)) = (\CayleyTrace{E}(x), 0)\)  and
\(\CayleyNorm{F}((x,y)) = (\CayleyNorm{E}(x)+\zeta \cdot \CayleyNorm{E}(y),0)\).
\end{proposition}

An especially important use case of the Cayley-Dickson process is when we start with some
associative and commutative ring with unit \((A,+,.)\), and repeatedly iterate the
process, starting with \((E,+,\times,\cdot,\sigma) = (A,+,.,Id_{A})\), yielding
Cayley algebra structures on \(A\), \(A \vartimes A\),...
When we do thus with \(A = \bbR\), and choose \(\zeta = 1\) at each step, we build
\(\bbC\), \(\bbH\), \(\bbO\), \(\bbX\) (see \ref{sec:fourtout}). We will see other
applications of this in the
remainder of this text. In such cases, it is convenient to identify \(E\) with
\(E \vartimes \{0\}\), or \(F\) with a superset of \(E\), at each step. This is one case
where one sees at once that the Cayley norm and trace actually takes their values in the first
rung, \emph{i.e.} \(A\), so we can drop the reference to precisely which algebra we
consider when computing them (which makes for far more readable
\(\CayleyTrace{}((x,y)) = \CayleyTrace{}(x)\)  and
\(\CayleyNorm{}((x,y)) = \CayleyNorm{}(x)+\zeta \cdot \CayleyNorm{}(y)\)).

We will refer to the Cayley-Dickson process as ``doubling'' when no confusion is to be feared.

\subsubsection{Cayley-Dickson for morphisms, and consequences}
\indent

\begin{theorem}\label{thm:mordbl}
Let \((E,+,\times,\cdot,\sigma)\) and \((F,\dagger,*,\bullet,s)\) be two Cayley
algebras over the commutative and associative ring with unit \((A,+,.)\), and
\(\phi : E \rightarrow F\) an \(A\)-Cayley morphism. Let
\(\Phi : E \vartimes E \rightarrow F \vartimes F; (x,y) \mapsto (\phi(x), \phi(y))\).
For any \(\alpha \in A\), denote by \((E_{\alpha},+_{\alpha},\times_{\alpha},
\cdot_{\alpha},\sigma_{\alpha})\) (respectively \((F_{\alpha},\dagger_{\alpha},
*_{\alpha},\bullet_{\alpha},s_{\alpha})\)) the result of applying the Cayley-Dickson
process to \((E,+,\times,\cdot,\sigma)\) (respectively \((F,\dagger,*,\bullet,s)\))
using \(\alpha\); then \(\Phi\) is also an \(A\)-Cayley morphism from
\(E_{\alpha}\) to \(F_{\alpha}\).
\end{theorem}

\begin{proof}
The \(A\)-linearity of \(\Phi\) is trivial.

\(\Phi(\sigma_{\alpha}(x,y)) = \Phi((\sigma(x), -y)) = (\phi(\sigma(x)), \phi(-y))\),
and as \(\phi\) is an \(A\)-Cayley morphism, \(\phi(\sigma(x)) = s(\phi(x))\) and
\(\phi(-y) = -\phi(y)\). Thus \(\Phi(\sigma_{\alpha}(x,y)) = (s(\phi(x)),-\phi(y)) =
s_{\alpha}(\phi(x),\phi(y)) = s_{\alpha}(\Phi((x,y)))\).

\begin{equation*}
\begin{aligned}
\Phi((x,y) \times_{\alpha} (z,t)) &=
\Phi((x \times z - \alpha \cdot (\sigma(t) \times y) , y \times \sigma(z)+t \times x))\\
 &=  (\phi(x) * \phi(z) - \alpha \bullet (s(\phi(t)) * \phi(y)) ,
 \phi(y) * s(\phi(z)) + \phi(t) * \phi(x))\\
 &= (\phi(x),\phi(y)) *_{\alpha} (\phi(z),\phi(t))\\
 &= \Phi((x,y)) *_{\alpha} \Phi((z,t))
 \end{aligned}
 \end{equation*}
 
 This proves that \(\Phi\) is a Cayley morphism.
\end{proof}

\begin{remark}
 If \(\phi\) is injective, so is \(\Phi\); if \(\phi\) is surjective, so is \(\Phi\).
\end{remark}

\begin{theorem}\label{thm:decdbl}
Let \(I\) be a set, \(((E_{i},+_{i},\times_{i},\cdot_{i},\sigma_{i}))_{i \in I}\) a
family of Cayley algebras on a same commutative and associative ring with unit
\((A,+,.)\), such that \((E,+,\times,\cdot,\sigma)\), the product algebra of that family
together with \linebreak \(\sigma : E \rightarrow E;
(x_{i})_{i \in I} \mapsto (\sigma_{i}(x_{i}))_{i \in I}\) form a Cayley algebra and that
every projection of \(\pi_{i} : E \rightarrow E_{i}\) is a Cayley (epi)morphism.

Let \(\alpha \in A\), and consider \((E_{\alpha},+_{\alpha},\times_{\alpha},
\cdot_{\alpha},\sigma_{\alpha})\) (respectively \((E_{i,\alpha},+_{i,\alpha},
\times_{i,\alpha},\cdot_{i,\alpha},\sigma_{i,\alpha})\)) the result of applying the
Cayley-Dickson process to \((E,+,\times,\cdot,\sigma)\) (respectively 
\((E_{i},+_{i},\times_{i},\cdot_{i},\sigma_{i})\)) using \(\alpha\),
\((F,\dagger,*,\bullet)\) the product algebra of the family \((E_{i,\alpha},
+_{i,\alpha},\times_{i,\alpha}, \cdot_{i,\alpha})_{i \in I}\) and
\(s : F \rightarrow F ; ((x_{i},y_{i}))_{i \in I} \mapsto
(\sigma_{i,\alpha}(x_{i},y_{i}))_{i \in I}\).
Then \(s\) is a conjugation on \(F\), \((F,\dagger,*,\bullet,s)\) is a Cayley algebra,
and \(\Phi : E_{\alpha} \rightarrow F ; ((x_{i})_{i \in I},(y_{i})_{i \in I})
\mapsto ((x_{i},y_{i})_{i \in I})\) is a Cayley isomorphism from
\((E_{\alpha},+_{\alpha},\times_{\alpha}, \cdot_{\alpha}, \sigma_{\alpha})\)
to  \((F,\dagger,*,\bullet,s)\).
\end{theorem}

\begin{remark}
\thref{thm:decdbl} simply states that, provided the product is dominated, the doubling of a
product of Cayley algebras is Cayley isomorphic to the product of the doubling of the
algebras.
\end{remark}

\begin{proof}
\thref{thm:struct} proves that \((E_{\alpha},+_{\alpha},\times_{\alpha},
\cdot_{\alpha}, \sigma_{\alpha})\) and the  \((E_{i,\alpha},+_{i,\alpha},
\times_{i,\alpha},\cdot_{i,\alpha},\sigma_{i,\alpha})\) are Cayley algebras.
\thref{thm:mordbl} proves that all the \(\Pi_{i} : E_{\alpha} \rightarrow
E_{i,\alpha} ; (x,y) \mapsto (\pi_{i}(x), \pi_{i}(y))\) are Cayley epimorphisms.

\(\Phi\) as defined above is trivially an isomorphism of algebra with unit. Since it is, in
particular, surjective, \thref{thm:domprod} proves \(s\) is a conjugation,
\((F,\dagger,*,\bullet,s)\) is a Cayley algebra and \(\Phi\) is a Cayley epimorphism.
And since \(\Phi\) is actually bijective, it is a Cayley isomorphism, as noted earlier.
\end{proof}

\subsection{Unimodulars and Quadratic Residues}
\indent

Let \((A,+,.)\) be a commutative and associative ring with unit
, \((E,+,\times,\cdot,\sigma)\) an \emph{alternative} Cayley algebra over \(A\) (with
neutral elements \(e\) and \(0\) for ``\(\times\)''  and ``\(+\)'' respectively), and
\(\rsfsA = A \cdot e\). We know that \((\rsfsA,+,\times)\) is an associative and 
commutative ring with
unit (with neutral element \(e\)), and denoting by \(\rsfsA^{*}\) the elements of
\(\rsfsA\) which are invertible, we see at once that \((\rsfsA^{*},\times)\) is an
abelian group. Furthermore, the Cayley norm takes its values in \(\rsfsA\), and
the Cayley norm, restricted to \(E^{*}\) takes its values in \(\rsfsA^{*}\).

\begin{proposition}\label{prp;normisloophomo}
\(\CayleyNorm{E}|_{E^{*}}^{\rsfsA^{*}}\) is a loop homomorphism.
\end{proposition}

\begin{proof}
This is merely a restatement of part of \thref{prp:altcay}.
\end{proof}

The \emph{unimodulars} of \(E\), which we will denote by \(\mathfrak{U}_{E}\), are
defined to be the kernel of \(\CayleyNorm{E}|_{E^{*}}^{\rsfsA^{*}}\), \emph{i.e.}
\allowbreak\(\mathfrak{U}_{E} = \{x \in E^{*} \thinspace|\thinspace
\CayleyNorm{E}(x) = e\}\).
This definition is of course consistent with what is defined under the same name in \(\bbC\)
(and \(\bbH\) and \(\bbO\)). It is known (\cite{Pflugfelder(1990)}) that the kernel of a
loop homomorphism is a normal subloop. Among other consequences, this entails that left and
right coset decompositions modulo \(\mathfrak{U}_{E}\) exist, that they coincide, and
that all cosets are equipotent (even when the loop is not finite, as is readily seen).
Furthermore, if \(E^{*}\) happens to be \emph{finite}, then \(\cardinal(E^{*}) =
\cardinal(\mathfrak{U}_{E}) \cardinal(E^{*} / \mathfrak{U}_{E})\) (where
\(\cardinal(X)\) denotes the cardinal of the set \(X\)).
Also, as \(E^{*}\) is a di-associative I.P. loop, then so is \(\mathfrak{U}_{E}\), and thus
the material of \ref{ssc:ordind} apply (note that if \(\omega \in \mathfrak{U}_{E}\) then
\(P_{\omega} \subset \mathfrak{U}_{E}\)). In particular, if
\(\cardinal(\mathfrak{U}_{E})\) is finite we have
\((\forall\thinspace \omega \in \mathfrak{U}_{E})\thickspace
\cardinal(\mathfrak{U}_{E}) =
[\mathfrak{U}_{E}:P_{\omega}] \cardinal(P_{\omega})\).

The \emph{quadratic residues} of \(E\), which we will denote by \(\mathfrak{R}_{E}\),
are defined to be the range of \(\CayleyNorm{E}|_{E^{*}}^{\rsfsA^{*}}\), \emph{i.e.}
\(\mathfrak{R}_{E} = \{y \in \rsfsA^{*} \thinspace|\thinspace
(\exists\thinspace x \in E^{*})\thickspace \CayleyNorm{E}(x) = y\} =
\{y \in E^{*} \thinspace|\thinspace
(\exists\thinspace x \in E^{*})\thickspace \CayleyNorm{E}(x) = y\}\).
This definition is of course consistent with the classical definition of quadratic residues
(\cite[...]{IrelandRosen(1982), JonesJones(1998)}), where one has \(E = A = \bbZ/n\bbZ\)
for some non-zero integer \(n\), and the conjugation is the identity. It is easy to
verify\footnote{The homomorphic image of a quasigroup into an associative quasigroup is
a quasigroup (\cite[page 29]{Pflugfelder(1990)}); the rest is immediate.}
that \((\mathfrak{R}_{E},\times)\) is an abelian subgroup of \((\rsfsA^{*},\times)\).
Since \(\mathfrak{U}_{E}\) is the kernel of
\(\CayleyNorm{E}|_{E^{*}}^{\rsfsA^{*}}\) and \(\mathfrak{R}_{E}\) is its range,
we also know (\cite{Pflugfelder(1990)}) that \(\mathfrak{R}_{E}\) is
loop-isomorphic (and hence group-isomorphic) to \(E^{*} / \mathfrak{U}_{E}\).
We have therefore proved the following:

\begin{proposition}\label{prp:unires}
If \(E\) is an alternative Cayley algebra and \(E^{*}\) is finite, then
\(\mathfrak{R}_{E}\) is group-isomorphic to \(E^{*} / \mathfrak{U}_{E}\) and
\begin{gather*}
\cardinal(E^{*}) = \cardinal(\mathfrak{U}_{E}) \cardinal(\mathfrak{R}_{E})\\
\cardinal(\mathfrak{R}_{E}) \thinspace|\thinspace \rsfsA^{*}
\end{gather*}
\end{proposition}

The following is also useful.

\begin{proposition}\label{prp:intouni}
If \(E\) is an alternative Cayley algebra,
\([\bbZ \rightarrow \mathfrak{R}_{E}; n \mapsto \CayleyNorm{E}(x^{n})]\) is a
group homomorphism; if in addition \(E^{*}\) is finite, then
\((\forall\thinspace x \in E^{*})\thickspace
\cardinal(\CayleyNorm{E}(P_{x})) \thinspace|\thinspace
\cardinal(\mathfrak{R}_{E})\).
\end{proposition}

\begin{proof}
Obviously, the range of the group homomorphism is \(\CayleyNorm{E}(P_{x})\).
\end{proof}

Note that, in particular, if \(E\) is a Galois Field, seen as an algebra over itself and with the
identity for conjugation, the quadratic residues are simply the non-zero squares, and we have
additional results.

\begin{proposition}\label{prp:evengal}
Let \(E\) be a Galois Field such that \(\cardinal(E)\) is \emph{even}, then
\(\cardinal(\mathfrak{U}_{E}) = 1\),
\(\cardinal(\mathfrak{R}_{E}) = \cardinal(E^{*}) = \cardinal(E)-1\);
in particular, every non-zero element is a square.
\end{proposition}

\begin{proof}
If \(\cardinal(E)\) is even then \([E \rightarrow E; x \mapsto x^2]\) is a bijection
as follows immediately from \ref{ssc:galeq}, so \(\mathfrak{U}_{E} = \{1_{F}\}\);
\end{proof}

To deal with Galois fields of odd cardinal, we will first state a simple lemma.

\begin{lemma}\label{lem:oddgal}
Let \(E\) be a Galois Field such that \(\cardinal(E)\) is \emph{odd}, let \(\gamma\)
be a generator of \(E^{*}\), and let \(\alpha \in E^{*}\). 

The set \(\{n \in \bbZ \thinspace|\thinspace \alpha = \gamma^{n}\}\) is always
non-empty. It contains either only even numbers or only odd numbers. \(\alpha\) is a square
if and only if \(\{n \in \bbZ \thinspace|\thinspace \alpha = \gamma^{n}\}\) contains
only even numbers.
\end{lemma}

\begin{proof}
The first statement is a direct consequence of the existence of generators for \(E^{*}\),
and does not depend upon the parity of \(\cardinal(E)\).

Let \(n_{0} \in \bbZ\) such that \(\alpha = \gamma^{n_{0}}\); we then have
\(\{n \in \bbN \thinspace|\thinspace \alpha = \gamma^{n}\} =
n_{0} +\cardinal(F^{*})\thinspace\bbZ\), also irrespective of the parity of
\(\cardinal(E)\), as if \(n\) verifies \(\alpha = \gamma^{n}\) then
\(\gamma^{n-n_{0}} = 1_{F}\), and \(\gamma\) is a generator of \(E^{*}\).

As \(\cardinal(E^{*}) = \cardinal(E)-1\), the statement about parity of the elements of
\(\{n \in \bbZ \thinspace|\thinspace \alpha = \gamma^{n}\}\) follows.

If \(\alpha = \gamma^{n_{0}}\) with \(n_{0} = 2m_{0}\), then
\(\alpha = (\gamma^{m_{0}})^2\) so is
a square. Conversely, if \(\alpha\) is a square, there exists \(\beta \in E^{*}\) such
that \(\alpha = \beta^{2}\). But there exists \(m_{0} \in \bbZ\) such that
\(\beta = \gamma^{m_{0}}\), so \(\alpha = \gamma^{2m_{0}}\), and the set
\(\{n \in \bbZ \thinspace|\thinspace \alpha = \gamma^{n}\}\) contains the even
number \(2m_{0}\).
\end{proof}

\begin{proposition}\label{prp:gasqnsq}
Let \(E\) be a Galois Field such that \(\cardinal(E)\) is \emph{odd}, then the product of
two non-zero squares or two non-zero non-squares is a square and the product of a non-zero
square by a non-zero non-square is a non-square.

There are exactly \((\cardinal(E)-1)/2\) non-zero squares and exactly as many non-zero
non-squares.
\end{proposition}

\begin{proof}
The first statement is a trivial consequence of \thref{lem:oddgal}. The second statement
is simple the constatation that \((\cardinal(E)-1)\) is even, therefore that there are as many
even numbers that there are odd numbers in \(\{1,\ldots,(\cardinal(E)-1)\}\), and that
each number in that set gives rise to a different element of \(E^{*}\).
\end{proof}

\section{Applications of finite alternative Cayley algebras}
\indent

The main goal of this section is to build concrete examples of finite alternative Cayley
algebras, and explicitly compute some important numbers related to them (such as the
number of invertible elements and the number of unimodular elements).

As announced in \ref{ssc:caydial}, we will first consider an associative and commutative
ring with unit \(A\); in \ref{sbs:cagaco} \(A\) will be a Galois field, and in
\ref{sbs:caylint} it will be a congruence ring. We then choose three elements in \(A\),
\(\alpha\), \(\beta\) and \(\gamma\), and repeatedly apply the Cayley-Dickson
procedure using these constants. More precisely, we start with \(E = A\), seen as a Cayley
algebra over \(A\), and using the identity as conjugations. Applying the Cayley-Dickson
process to \(E\) using \(\alpha\) yields a Cayley algebra we will denote by
\(E_{\alpha}\). Then we apply the Cayley-Dickson process to \(E_{\alpha}\) using
\(\beta\) and yielding \(E_{\alpha,\beta}\). Finally, we apply the Cayley-Dickson
process one last time, to \(E_{\alpha,\beta}\) using \(\gamma\) and yielding
\(E_{\alpha,\beta,\gamma}\). We could keep doing this indefinitely, of course, but
beyond \(E_{\alpha,\beta,\gamma}\) the algebraic properties are usually so poor that there
is no real incentive to do so, and no outside uses have been identified which would make it
worthwhile either; technically, the fact that we could no longer rely on \thref{prp:altcay}
would also prove bothersome. As indicated in \ref{ssc:caydial}, this is a setting in which it
is sometimes advantageous to consider the successive algebras as supersets of each others, and
at any rate most convenient to consider the various Cayley norms as taking their values in
\(A\) (the identification of \(A\) with \(\rsfsA = A \cdot e\), where \(e\) is the
neutral element of one of the successive algebras is possible in this case, as the mapping
\([x \mapsto x \cdot e]\) is bijective). In particular, this leads to the following simple
expressions:
\begin{align*}
\CayleyNorm{E}(x) &= x^{2}\\
\CayleyNorm{E_{\alpha}}((x,y)) &= x^{2} + \alpha x^{2}\\
\CayleyNorm{E_{\alpha,\beta}}((x,y,z,t))& =
(x^{2}+ \alpha y^{2}) + \beta (z^{2}+ \alpha t^{2})\\
\CayleyNorm{E_{\alpha,\beta,\gamma}}((x,y,z,t,u,v,w,s)) &=
((x^{2}+ \alpha y^{2}) + \beta (z^{2}+ \alpha t^{2})) +
\gamma ((u^{2}+ \alpha v^{2}) + \beta (w^{2}+ \alpha s^{2}))
\end{align*}

\subsection{Some common properties}
\indent

We collect here some properties which are scattered in the previous section, and give some
immediate consequences.

We consider a commutative and associative ring with unit
 \((A,+,.)\) and \((E,+,\times,\cdot,\sigma)\) a \emph{finite} and \emph{alternative}
 Cayley algebra over \(A\) (with neutral elements \(e\) and \(0\)  for ``\(\times\)''
and ``\(+\)'' respectively), and \(\rsfsA = A \cdot e\).

Then we know that \((E^{*},\times)\) and \((\mathfrak{U}_{E},\times)\) are finite
di-associative I.P. loops, that \((\rsfsA^{*},\times)\) and
\((\mathfrak{R}_{E},\times)\) are finite abelian groups, that for all \(x \in E^{*}\),
\((P_{x},\times)\) is an abelian group, and that we have
\begin{gather}
(\forall\thinspace \omega \in E^{*})\thickspace
\cardinal(E^{*}) = [E^{*}:P_{\omega}] \cardinal(P_{\omega})\tag{P1}\label{equ:P1}\\
(\forall\thinspace \omega \in \mathfrak{U}_{E})\thickspace
\cardinal(\mathfrak{U}_{E}) =
[\mathfrak{U}_{E}:P_{\omega}] \cardinal(P_{\omega})\tag{P2}\label{equ:P2}\\
\cardinal(E^{*}) =
\cardinal(\mathfrak{U}_{E}) \cardinal(\mathfrak{R}_{E})\tag{P3}\label{equ:P3}\\
 \cardinal(\mathfrak{R}_{E}) \thinspace|\thinspace \rsfsA^{*}\tag{P4}\label{equ:P4}\\
 (\forall\thinspace x \in E^{*})\thickspace
\cardinal(\CayleyNorm{E}(P_{x})) \thinspace|\thinspace
\cardinal(\mathfrak{R}_{E})\tag{P5}\label{equ:P5}
\end{gather}

\begin{proposition}
\((\forall\thinspace x \in E^{*})
(\exists m \thinspace|\thinspace \cardinal(\mathfrak{R}_{E}))\thickspace
\{n \in \bbZ \thinspace|\thinspace x^{n} \in \mathfrak{U}_{E}\} = m \bbZ\);
in particular
\((\forall\thinspace x \in E^{*})\thickspace
x^{\cardinal(\mathfrak{R}_{E})} \in \mathfrak{U}_{E}\)
and
\((\forall\thinspace x \in E^{*})\thickspace
x^{\cardinal(\rsfsA^{*})} \in \mathfrak{U}_{E}\).

\((\forall\thinspace x \in E^{*})
(\exists m \thinspace|\thinspace \cardinal(E^{*}))\thickspace
\{n \in \bbZ \thinspace|\thinspace x^{n} = e\} = m \bbZ\);
in particular
\((\forall\thinspace x \in E^{*})\thickspace
x^{\cardinal(E^{*})} = e\).

\((\forall\thinspace x \in \mathfrak{U}_{E})
(\exists m \thinspace|\thinspace \cardinal(\mathfrak{U}_{E}))\thickspace
\{n \in \bbZ \thinspace|\thinspace x^{n} = e\} = m \bbZ\);
in particular
\((\forall\thinspace x \in \mathfrak{U}_{E})\thickspace
x^{\cardinal(\mathfrak{U}_{E})} = e\).
\end{proposition}

\begin{proof}
The first statement is a consequence \thref{prp:intouni}. The second statement uses
(\ref{equ:P5}) in addition to the first
statement. The third statement uses (\ref{equ:P4}) in addition.

The fourth statement is a consequence of the fact that \((P_{x},\times)\) is an abelian
group, and the fifth is an application of (\ref{equ:P1}).

The sixth statement is also a consequence of the fact that \((P_{x},\times)\) is an abelian
group, and the seventh is an application of (\ref{equ:P2}).
\end{proof}

\subsection{Preliminary computations}
\indent

We will use the notations of Appendix~\ref{apx:gjs}, and the results therein, which are
merely  rewritten from \cite{IrelandRosen(1982)}. This subsection essentially builds upon
the equivalent to \cite[page~106, Exercises 19 and 20]{IrelandRosen(1982)} for any
finite field. We will write here \(q = \cardinal(F)\).

Given a finite field \(F\), \(a_{1} \in F^{*}\), \ldots, \(a_{r} \in F^{*}\), and
\(b \in F\), we will denote by \(N(a_{1}x_{1}^2+\cdots+a_{r}x_{r}^2 = b)\),
or just \(N\) as no confusion is to be feared here, the
number of solutions of the equation \(a_{1}x_{1}^2+\cdots+a_{r}x_{r}^2 = b\) in
\(F\), and by \(N(a_{i}x^2=c)\) the number of solutions of \(a_{i}x^2=c\) in \(F\).
We impose \(r \geqslant 1\) here. Clearly,
\(N(a_{1}x_{1}^2+\cdots+a_{r}x_{r}^2 = b) =
\sum_{\alpha_{1}+\cdots+\alpha_{r} = b}
N(a_{1}x_{1}^2 = \alpha_{1}) \cdots N(a_{r}x_{r}^2 = \alpha_{r})\).

If \(\cardinal(F)\) is even,  \([F \rightarrow F; x \mapsto x^2]\) is a bijection
and thus \(N = \sum_{\alpha_{1}+\cdots+\alpha_{r} = b} 1 = q^{r-1}\).

If \(\cardinal(F)\) is odd, then \(N =
\sum_{\alpha_{1}+\cdots+\alpha_{r} = b}\thinspace \prod_{i=1}^{r}
\sum_{\chi_{i}^{2} = \epsilon_{F}} \check{\chi}_{i}(\frac{\alpha_{i}}{a_{i}})\),
because \(2 \thinspace|\thinspace (q-1)\),
and we therefore need to determine the number of multiplicative characters of \(F\) whose
order divide \(2\). As the group of multiplicative characters of \(F\) is cyclic, of cardinal
\(q-1\), which is even, the the number of multiplicative characters of \(F\) whose order
divide \(2\) is exactly \(2\): the trivial character \(\epsilon_{F}\) and one we shall call
\(\rsfsX_{2}\) and which is merely a generalization of the Legendre symbol as shown
below.

\begin{proposition}
Let \(F\) be a finite field such that \(q = \cardinal(F)\) is \emph{odd}; then
\begin{equation*}
\check{\rsfsX}_{2}(\alpha) =
\begin{cases}
+1 & \text{if \((\exists\thinspace \beta \in F^{*})\thickspace \alpha = \beta^{2}\)} \\
0 & \text{if \(\alpha = 0\)} \\
-1 & \text{otherwise.}
\end{cases}
\end{equation*}
 is a mutiplicative character of order two and
 \(\rsfsX_{2}(-1_{F}) = (-1)^{\frac{q-1}{2}}\); in particular
\(\rsfsX_{2} \neq \epsilon_{F}\).
\end{proposition}

\begin{proof}
It is immediate to check (using \thref{lem:oddgal}) that \(\rsfsX_{2}\) is a multiplicative
character,  and since it is trivial that \((\forall\thinspace \alpha \in F^{*})\thickspace
(\rsfsX_{2}(\alpha))^2 = 1\), its order divides \(2\). Since by \thref{prp:gasqnsq}
there exists non-zero non-squares, \(\rsfsX_{2}\) does take the value \(-1\), and thus
its order can't be one, so is exactly two. At the same time,
\(\rsfsX_{2} \neq \epsilon_{F}\).

The equation \(x^2 = -1_{F}\) in \(F\) has solutions if and only if
\((-1_{F})^{\frac{q-1}{2}} = 1_{F}\) (according to \ref{ssc:galeq}; in this case
\(d = \PGCD(2, q-1) = 2\)). Hence \(-1_{F}\) is square if and only if
\((-1)^{\frac{q-1}{2}} = 1\) (which is equivalent to \(q \equiv 1 \modulus 4\)),
and thus \(\rsfsX_{2}(-1_{F}) = (-1)^{\frac{q-1}{2}}\).
\end{proof}

We can now see that
\begin{equation*}
\begin{aligned}
N &=
\sum_{\alpha_{1}+\cdots+\alpha_{r} = b}\thinspace \prod_{i=1}^{r}
\left[\check{\epsilon}_{F}\left(\frac{\alpha_{i}}{a_{i}}\right)+
\check{\rsfsX}_{2}\left(\frac{\alpha_{i}}{a_{i}}\right)\right] \\
&= \sum_{\alpha_{1}+\cdots+\alpha_{r} = b}\thinspace
\sum_{\omega \in \rsfsP(\{1,\ldots,r\})}\thinspace
\check{\rsfsX}_{2}^{\omega(1)}\left(\frac{\alpha_{1}}{a_{1}}\right)\cdots
\check{\rsfsX}_{2}^{\omega(r)}\left(\frac{\alpha_{r}}{a_{r}}\right)
\end{aligned}
\end{equation*}
where \(\rsfsP(X)\) denotes the power set of \(X\), and we have
\(\rsfsX_{2}^{0} = \epsilon_{F}\) as in any (multiplicative) group.

Therefore
\begin{equation*}
\begin{aligned}
N &=
\sum_{\omega \in \rsfsP(\{1,\ldots,r\})}\thinspace
\sum_{\alpha_{1}+\cdots+\alpha_{r} = b}\thinspace
\check{\rsfsX}_{2}^{\omega(1)}\left(\frac{\alpha_{1}}{a_{1}}\right)\cdots
\check{\rsfsX}_{2}^{\omega(r)}\left(\frac{\alpha_{r}}{a_{r}}\right) \\
&= \sum_{\omega \in \rsfsP(\{1,\ldots,r\})}
\left[\check{\rsfsX}_{2}^{\omega(1)}(a_{1})\cdots\check{\rsfsX}_{2}^{\omega(r)}(a_{r})\right]
\sum_{\alpha_{1}+\cdots+\alpha_{r} = b}\thinspace
\check{\rsfsX}_{2}^{\omega(1)}(\alpha_{1})\cdots\check{\rsfsX}_{2}^{\omega(r)}(\alpha_{r})
\end{aligned}
\end{equation*}
as both \(\epsilon_{F}\) and \(\rsfsX_{2}\) are multiplicative characters (so
\(\check{\epsilon}_{F}\left(\frac{\alpha_{i}}{a_{i}}\right) =
\check{\epsilon}_{F}(\alpha_{i}) \thinspace \check{\epsilon}_{F}(a_{i}^{-1})\)
and \(\check{\rsfsX}_{2}\left(\frac{\alpha_{i}}{a_{i}}\right) =
\check{\rsfsX}_{2}(\alpha_{i}) \thinspace \check{\rsfsX}_{2}(a_{i}^{-1})\)) and
are of an order dividing \(2\) (so \(\epsilon_{F}(a_{i}^{-1}) = \epsilon_{F}(a_{i})\)
and \(\rsfsX_{2}(a_{i}^{-1}) = \rsfsX_{2}(a_{i})\)).

We now come to a point where we must distinguish between \(b = 0\) and \(b \neq 0\).
\begin{equation*}
N =
\begin{cases}
\sum_{\omega \in \rsfsP(\{1,\ldots,r\})}
\left[\check{\rsfsX}_{2}^{\omega(1)}(a_{1})\cdots\check{\rsfsX}_{2}^{\omega(r)}(a_{r})\right]
J_{0}(\rsfsX_{2}^{\omega(1)},\ldots,\rsfsX_{2}^{\omega(r)}) & \text{if \(b = 0\)} \\
\sum_{\omega \in \rsfsP(\{1,\ldots,r\})}
\left[\left(\rsfsX_{2}^{\omega(1)}\cdots\rsfsX_{2}^{\omega(r)}\right)(b)\right]\thinspace
\left[\check{\rsfsX}_{2}^{\omega(1)}(a_{1})\cdots\check{\rsfsX}_{2}^{\omega(r)}(a_{r})\right]
J(\rsfsX_{2}^{\omega(1)},\ldots,\rsfsX_{2}^{\omega(r)}) & \text{if \(b \neq 0\)}
\end{cases}
\end{equation*}

Using the results in \ref{sbs:jacosum} we find that
\begin{equation*}
N =
\begin{cases}
q^{r-1} + \rsfsX_{2}(a_{1} \cdots a_{r})\thinspace J_{0}(\overbrace{\rsfsX_{2},\ldots,\rsfsX_{2}}^{\text{\(r\) terms}}) & \text{if \(b = 0\)} \\
q^{r-1} + \rsfsX_{2}^{r}(b) \thinspace \rsfsX_{2}(a_{1} \cdots a_{r})\thinspace J(\overbrace{\rsfsX_{2},\ldots,\rsfsX_{2}}^{\text{\(r\) terms}}) & \text{if \(b \neq 0\)}
\end{cases}
\end{equation*}

We now must also distinguish between whether \(r\) is even or odd.

If \(r\) is odd, then \(\rsfsX_{2}^{r-1} = \epsilon_{F}\) and
\(\rsfsX_{2}^{r} = \rsfsX_{2} \neq \epsilon_{F}\), so
\(J_{0}(\overbrace{\rsfsX_{2},\ldots,\rsfsX_{2}}^{\text{\(r\) terms}}) = 0\) and
\(g(\rsfsX_{2})^{r} =
J(\overbrace{\rsfsX_{2},\ldots,\rsfsX_{2}}^{\text{\(r\) terms}})\thinspace
g(\rsfsX_{2}^{r})\), and since \(\rsfsX_{2}^{r} = \rsfsX_{2}\),
\(J(\overbrace{\rsfsX_{2},\ldots,\rsfsX_{2}}^{\text{\(r\) terms}} =
g(\rsfsX_{2})^{r-1}\). But \(\rsfsX_{2}^{-1} = \overline{\rsfsX}_{2} = \rsfsX_{2}\)
so \(g(\rsfsX_{2})^{2} = \rsfsX_{2}(-1_{F}) \thinspace q\) and
\(J(\overbrace{\rsfsX_{2},\ldots,\rsfsX_{2}}^{\text{\(r\) terms}}) =
[\rsfsX_{2}(-1_{F}) \thinspace q]^{\frac{r-1}{2}}\). Combining this with the fact
that  \(\rsfsX_{2}(-1_{F}) = (-1)^{\frac{q-1}{2}}\) we finally find that
\(N = q^{r-1}\) if \(b = 0\), and \(N = q^{r-1} + \rsfsX_{2}(b)\thinspace
\rsfsX_{2}(a_{1} \cdots a_{r})\thinspace (-1)^{\frac{q-1}{2}\thinspace \frac{r-1}{2}}\thinspace
q^{\frac{r-1}{2}}\) otherwise.

If \(r\) is even, \(\rsfsX_{2}^{r-1} = \rsfsX_{2} \neq \epsilon_{F}\)
and \(\rsfsX_{2}^{r} = \epsilon_{F}\), so
\(J_{0}(\overbrace{\rsfsX_{2},\ldots,\rsfsX_{2}}^{\text{\(r\) terms}}) =
\rsfsX_{2}(-1_{F})\thinspace (q-1)\thinspace
J(\overbrace{\rsfsX_{2},\ldots,\rsfsX_{2}}^{\text{\((r-1)\) terms}})\) and
\(J(\overbrace{\rsfsX_{2},\ldots,\rsfsX_{2}}^{\text{\((r-1)\) terms}})\thinspace
g(\rsfsX_{2}^{r-1}) = g(\rsfsX_{2})^{r-1}\). Together these give
\(J_{0}(\overbrace{\rsfsX_{2},\ldots,\rsfsX_{2}}^{\text{\(r\) terms}}) =
\rsfsX_{2}(-1_{F})\thinspace (q-1)\thinspace g(\rsfsX_{2})^{r-2} =
\rsfsX_{2}(-1_{F})\thinspace (q-1)\thinspace
[\rsfsX_{2}(-1_{F})\thinspace q]^{\frac{r}{2}-1}\). Furthermore
\(J(\overbrace{\rsfsX_{2},\ldots,\rsfsX_{2}}^{\text{\(r\) terms}}) =
\rsfsX_{2}(-1_{F})\thinspace
J(\overbrace{\rsfsX_{2},\ldots,\rsfsX_{2}}^{\text{\((r-1)\) terms}})\) and
\(g(\rsfsX_{2})^{r-1} = g(\rsfsX_{2}^{r-1})\thinspace
J(\overbrace{\rsfsX_{2},\ldots,\rsfsX_{2}}^{\text{\((r-1)\) terms}})\), so
\(J(\overbrace{\rsfsX_{2},\ldots,\rsfsX_{2}}^{\text{\((r-1)\) terms}}) =
g(\rsfsX_{2})^{r-2}\) and
\(J(\overbrace{\rsfsX_{2},\ldots,\rsfsX_{2}}^{\text{\(r\) terms}}) =
-\rsfsX_{2}(-1_{F})\thinspace g(\rsfsX_{2})^{r-2} =
-\rsfsX_{2}(-1_{F}\thinspace [\rsfsX_{2}(-1_{F}\thinspace q)]^{\frac{r}{2}-1})\)
so finally \(N = q^{r-1} + \rsfsX_{2}(a_{1} \cdots a_{r})\thinspace
(-1)^{\frac{q-1}{2}\thinspace\frac{r}{2}}\thinspace
(q-1)\thinspace q^{\frac{r}{2}-1}\) if \(b = 0\), and
\(N = q^{r-1} - \rsfsX_{2}(b)\thinspace \rsfsX_{2}(a_{1} \cdots a_{r})\thinspace
(-1)^{\frac{q-1}{2}\thinspace\frac{r}{2}}\thinspace q^{\frac{r}{2}-1}\)
otherwise.

We summarize the value of \(N(a_{1}x_{1}^2+\cdots+a_{r}x_{r}^2 = b)\) as follows:

%
\begin{center}
{\renewcommand{\arraystretch}{1.5}
\begin{tabular}{||c|c|c||c||}
\hhline{|t:===:t:=:t|}
\multicolumn{3}{||l||}{\(q\) even} & \(q^{r-1}\) \\[1mm]
\hhline{||----||}
\(q\) odd & \(b = 0\) & \(r\) odd & \(q^{r-1}\) \\[1mm]
\hhline{||~~--||}
&& \(r\) even & \(q^{r-1} + \rsfsX_{2}(a_{1} \cdots a_{r})\thinspace
(-1)^{\frac{q-1}{2}\thinspace\frac{r}{2}}\thinspace
(q-1)\thinspace q^{\frac{r}{2}-1}\) \\[1mm]
\hhline{||~---||}
& \(b \neq 0\) & \(r\) odd & \(q^{r-1} + \rsfsX_{2}(b)\thinspace
\rsfsX_{2}(a_{1} \cdots a_{r})\thinspace \thinspace (-1)^{\frac{q-1}{2}\thinspace \frac{r-1}{2}}\thinspace
q^{\frac{r-1}{2}}\) \\[1mm]
\hhline{||~~--||}
&& \(r\) even & \(q^{r-1} - \rsfsX_{2}(b)\thinspace \rsfsX_{2}(a_{1} \cdots a_{r})\thinspace
(-1)^{\frac{q-1}{2}\thinspace\frac{r}{2}}\thinspace q^{\frac{r}{2}-1}\) \\[1mm]
\hhline{|b:===:b:=:b|}
\end{tabular}
}
\end{center}

\subsection{Cayley-Galois constructs}\label{sbs:cagaco}
\indent

As announced at the begining of this section, we consider here \(A = F_{q}\), with
\(F_{q}\) be a Galois field of cardinal \(q\). We will designate by \(F_{q;\alpha}\)
\emph{etc.} what was denoted by \(E_{\alpha}\) \emph{etc.} in the introduction to this
section.

We now assume \(\alpha\beta\gamma\neq 0\), for brevity's sake. The preliminary
computations above yield the table below.

Note that \(F_{q;\alpha,\beta}\) and \(F_{q;\alpha,\beta,\gamma}\) never are fields,
but that when \(q\) is even, \(F_{q;\alpha}\) is never a field whereas when \(q\) is odd,
\(F_{q;\alpha}\) sometimes \emph{is} a field (it is a field if and only if
\(\rsfsX_{2}(\alpha)\thinspace(-1)^{\frac{q-1}{2}} = -1\)).

It is interesting to note that, when \(q\) is even, \(F_{q}\) is of characteristic \(2\), and
thus we always have \(-x = +x\), so, as a consequence of \thref{thm:struct} and the form
of the conjugation resulting from the Cayley-Dickson process, \(F_{q}\),
\(F_{q;\alpha}\), \(F_{q;\alpha,\beta}\) and \(F_{q;\alpha,\beta,\gamma}\), and
actually any successor in the doubling scheme, are commutative and associative (irrespective
of the value of \(\alpha\), \(\beta\), \(\gamma\) or their successors), and the
conjugation at each step is the identity.

\begin{center}
{\renewcommand{\arraystretch}{1.5}
\begin{tabular}{||l||c|c||}
\hhline{~|t:==:t|}
\multicolumn{1}{c||}{(\(\alpha\beta\gamma\neq 0\))} & \(q\) odd & \(q\) even \\[1mm]
\hhline{|t:=::=|=||}
\(\cardinal(\rsfsA^{*})\) & \(q-1\) & \(q-1\) \\[1mm]
\hhline{||---||}
\(\cardinal(F_{q}^{*})\) & \(q-1\) & \(q-1\) \\[1mm]
\hhline{||---||}
\(\cardinal(F_{q;\alpha}^{*})\) & \(q^{2} - \left[q+ \rsfsX_{2}(\alpha)\thinspace(q-1)\thinspace(-1)^{\frac{q-1}{2}}\right]\) & \(q^{2}-q\) \\[1mm]
\hhline{||---||}
\(\cardinal(F_{q;\alpha,\beta}^{*})\) & \(q^{4} - \left[q^{3}+(q-1)\thinspace q\right]\) & \(q^{4}-q^{3}\) \\[1mm]
\hhline{||---||}
\(\cardinal(F_{q;\alpha,\beta,\gamma}^{*})\) & \(q^{8} - \left[q^{7}+(q-1)\thinspace q^{3}\right]\) & \(q^{8}-q^{7}\) \\[1mm]
\hhline{||---||}
\(\cardinal(\mathfrak{U}_{F_{q}})\) & \(2\) & \(1\) \\[1mm]
\hhline{||---||}
\(\cardinal(\mathfrak{U}_{F_{q;\alpha}})\) & \(q - \rsfsX_{2}(\alpha)\thinspace(-1)^{\frac{q-1}{2}}\) & \(q\) \\[1mm]
\hhline{||---||}
\(\cardinal(\mathfrak{U}_{F_{q;\alpha,\beta}})\) &\(q^{3}-q\)& \(q^{3}\) \\[1mm]
\hhline{||---||}
\(\cardinal(\mathfrak{U}_{F_{q;\alpha,\beta,\gamma}})\) & \(q^{7}-q^{3}\) & \(q^{7}\) \\[1mm]
\hhline{|b:=:b:==:b|}
\end{tabular}
}
\end{center}

\pagebreak
\subsection{Cayley integers}\label{sbs:caylint}
\indent

As announced at the beginning of this section, we consider here
\(A = \bbZ_{n} = \bbZ/n\bbZ\), with \(n \in \bbN\), not necessarily prime (but still
\(n \geqslant 1\), of course). We will designate by \(\bbZ_{n;\alpha}\)
\emph{etc.} what was denoted by \(E_{\alpha}\) \emph{etc.} in the introduction to this
section. Elements of these Cayley algebras will be dubbed \emph{Cayley integers}.

The situation here is somewhat more complicated than that of \ref{sbs:cagaco}.

Let \(n = \prod_{i = 1}^{k}n_{i}\) where the \(n_{i} \in \bbN^{*}\) are relatively
prime, and let \(\pi_{n;n_{i}}:\bbZ_{n} \rightarrow \bbZ_{n_{i}}\) such that if
\(X \in \bbZ_{n}\) and \(x \in \bbZ\), \(x \in X\), then \(pi_{n;n_{i}}(X) = Y\)
such that if \(y \in \bbZ\), \(y \in Y\) then \(y \equiv x \modulus n_{i}\).
The Chinese Remainder Theorem informs us that, as algebras with unit over
\(A = \bbZ_{n}\), \(\bbZ_{n}\) and
\(\bbZ_{n_{1}} \vartimes\cdots\vartimes \bbZ_{n_{k}}\) are
isomorphic, \emph{via}
\(\pi_{n} : \bbZ_{n} \rightarrow \bbZ_{n_{1}} \vartimes\cdots\vartimes \bbZ_{n_{k}};
x \mapsto (\pi_{n;n_{1}}(x),\ldots,\pi_{n;n_{k}}(x))\).
By using the identity on \(\bbZ\) and the \(\bbZ_{n_{i}}\) as
conjugations, \(\bbZ_{n}\) and the \(\bbZ_{n_{i}}\) are Cayley algebras on the same
ring \(A = \bbZ_{n}\), the \(\pi_{n;n_{i}}\) are Cayley epimorphisms.
\thref{thm:domprod} now informs us
that, using the identity on \(\bbZ_{n_{1}} \vartimes\cdots\vartimes \bbZ_{n_{k}}\)
as conjugations, \(\bbZ_{n_{1}} \vartimes\cdots\vartimes \bbZ_{n_{k}}\) is also
a Cayley algebra over \(A = \bbZ_{n}\). In fact we see that
\(\bbZ_{n_{1}} \vartimes\cdots\vartimes \bbZ_{n_{k}}\) is
Cayley-isomorphic to \(\bbZ\) (it is isomorphic as an algebra with unit by the Chinese
remainder theorem, and the conjugation is the identity on both sides).
Therefore, \thref{thm:decdbl} proves that, given \(\alpha \in \bbZ_{n}\),
\(\bbZ_{n_{1};\alpha} \vartimes\cdots\vartimes \bbZ_{n_{k};\alpha}\) is a
Cayley algebra which is Cayley-isomorphic with the result of applying the Cayley-Dickson
procedure to \(\bbZ_{n_{1}} \vartimes\cdots\vartimes \bbZ_{n_{k}}\) using
\(\alpha\), \emph{i.e.} Cayley-isomorphic to \(\bbZ_{n;\alpha}\),
as detailed in \thref{thm:decdbl} (in somewhat more details: to
\((x,y) \in \bbZ_{n;\alpha}\) we associate
\(((\pi_{n;n_{1}}(x),\pi_{n;n_{1}}(y)),\ldots,
(\pi_{n;n_{k}}(x),\pi_{n;n_{k}}(y))) \in 
\bbZ_{n_{1};\alpha} \vartimes\cdots\vartimes \bbZ_{n_{k};\alpha}\)).
We can iterate the above reasoning, and we find that
\(\bbZ_{n_{1};\alpha,\beta} \vartimes\cdots\vartimes
\bbZ_{n_{k};\alpha,\beta}\) (respectively
\(\bbZ_{n_{1};\alpha,\beta,\gamma} \vartimes\cdots\vartimes
\bbZ_{n_{k};\alpha,\beta,\gamma}\)),
is a Cayley algebra isomorphic to \(\bbZ_{n;\alpha,\beta}\) (respectively
\(\bbZ_{n;\alpha,\beta,\gamma}\)), in a way similar to that detailed above.

Using the above isomorphisms we see that an elements of \(\bbZ_{n}\) (respectively
\(\bbZ_{n;\alpha}\), \(\bbZ_{n;\alpha,\beta}\), \(\bbZ_{n;\alpha,\beta,\gamma}\))
is invertible if and only if \(\pi_{n;n_{i}}(x)\) (respectively
\(\pi_{n;n_{i}}^{\alpha}(x)\), \(\pi_{n;n_{i}}^{\alpha,\beta}(x)\),
 \(\pi_{n;n_{i}}^{\alpha,\beta,\gamma}(x)\), where \(\pi_{n;n_{i}}^{\alpha}\)
 is the result of applying the Cayley-Dickson procedure to \(\pi_{n;n_{i}}\) using
 \(\alpha\), \emph{etc.}) is invertible for all \(i\).
 Therefore \(\bbZ_{n}^{*}\) (respectively
\(\bbZ_{n;\alpha}^{*}\), \(\bbZ_{n;\alpha,\beta}^{*}\),
\(\bbZ_{n;\alpha,\beta,\gamma}^{*}\)) is equipotent to
\(\prod_{i = 1}^{k}\bbZ_{n_{i}}^{*}\) (respectively
\(\prod_{i = 1}^{k}\bbZ_{n_{i};\alpha}^{*}\),
\(\prod_{i = 1}^{k}\bbZ_{n_{i};\alpha,\beta}^{*}\),
\(\prod_{i = 1}^{k}\bbZ_{n_{i};\alpha,\beta,\gamma}^{*}\)).
We find a similar statement for unimodulars by noting that \(\pi_{n}(1) = (1,\ldots,1)\).
In short:
\begin{align*}
&\cardinal(\bbZ_{n}^{*}) = \prod_{i = 1}^{k} \cardinal(\bbZ_{n_{i}}^{*}) &&\cardinal(\mathfrak{U}_{\bbZ_{n}}) = \prod_{i = 1}^{k} \cardinal(\mathfrak{U}_{\bbZ_{n_{i}}})\\
&\cardinal(\bbZ_{n;\alpha}^{*}) = \prod_{i = 1}^{k} \cardinal(\bbZ_{n_{i}\alpha}^{*}); &&\cardinal(\mathfrak{U}_{\bbZ_{n;\alpha}}) = \prod_{i = 1}^{k} \cardinal(\mathfrak{U}_{\bbZ_{n_{i};\alpha}})\\
&\cardinal(\bbZ_{n;\alpha,\beta}^{*}) = \prod_{i = 1}^{k} \cardinal(\bbZ_{n_{i};\alpha,\beta}^{*}) &&\cardinal(\mathfrak{U}_{\bbZ_{n;\alpha,\beta}}) = \prod_{i = 1}^{k} \cardinal(\mathfrak{U}_{\bbZ_{n_{i};\alpha,\beta}})\\
&\cardinal(\bbZ_{n;\alpha,\beta,\gamma}^{*}) = \prod_{i = 1}^{k} \cardinal(\bbZ_{n_{i};\alpha,\beta,\gamma}^{*}) &&\cardinal(\mathfrak{U}_{\bbZ_{n;\alpha,\beta,\gamma}}) = \prod_{i = 1}^{k} \cardinal(\mathfrak{U}_{\bbZ_{n_{i};\alpha,\beta,\gamma}})
\end{align*}

A last touch is to notice that the number of invertible (or unimodular) elements in
\(\bbZ_{n_{i};\alpha}\) (\emph{etc.}) seen as a Cayley algebra over \(\bbZ_{n}\) is the
same as when seen as a Cayley algebra over \(\bbZ_{n_{i}}\).

Finally, to compute the cardinals of interest to us in this segment, we simply have to do so
for integers of the form \(n = p^{s}\), with \(p\) prime and \(s \geqslant 1\).

If \(s = 1\), then we are face to a special case of the Cayley-Galois constructs we have just
studied.

The elements of \(\bbZ_{p^{s}}\) which are not invertible, are, as is well-known, exactly
those of the form \(\bar{x} = k \bar{p}\), for some integer \(k\), where \(\bar{a}\)
represents the congruence class of the integer \(a\) modulo \(p^{s}\), a convention we will
use throughout this segment (and not to be confused with the notation for the conjugation).
Hence there are \(p^{s-1}\) elements which are not invertible, and thus
\(\cardinal(\bbZ_{p^{s}}^{*}) = \cardinal(\rsfsA^{*}) = p^{s-1}\thinspace(p-1)\).

Assume \(p \neq 2\). We will further assume, for brevity's sake, that
\(\alpha \in \bbZ_{p^{s}}^{*}\), \(\beta \in \bbZ_{p^{s}}^{*}\),
\(\gamma \in \bbZ_{p^{s}}^{*}\).

Let us first determine the number of invertible elements. As we will see, finding the number
of unimodulars will use essentially the same techniques.

More to the point, we will determine the number of elements which are \emph{not}
invertible, and assume \(s > 1\). \thref{prp:altcay} informs us that an elements of
\(\bbZ_{p^{s}}\) (respectively \(\bbZ_{p^{s};\alpha}\),
\(\bbZ_{p^{s};\alpha,\beta}\), \(\bbZ_{p^{s};\alpha,\beta,\gamma}\)) is not
invertible if and only if its Cayley norm is
not invertible in \(\rsfsA\), which here is Cayley-isomorphic to \(\bbZ_{p^{s}}\).
Taking into account the expression of the Cayley norm on \(\bbZ_{p^{s};\alpha}\)
\emph{etc.} given at the beginning of this section, we see that finding the non-invertible
element of \(\bbZ_{p^{s};\alpha}\) amounts to finding couples
\emph{etc.} of elements of \(\bbZ_{p^{s}}\) such that
\(x_{1}^{2}+\alpha\thinspace x_{2}^{2}\) \emph{etc.} is not invertible in
\(\bbZ_{p^{s}}\).
As an element \(x_{i} \in \bbZ_{p^{s}}\) is completely characterized by the integer
\(\bm{x}_{i} \in \{0,\ldots,p^{s}-1\}\) such that \(\bm{x}_{i} \in x_{i}\), we are
(unsurprisingly) led to solving a congruence equation in \(\bbZ\). We will find it convenient
to write \(\bm{x}_{i}\) in the basis \(p\), so that
\(\bm{x}_{i} = x_{i;0}\thinspace p^{0} + x_{i;1}\thinspace p^{1}+\cdots
+x_{i;s-1}\thinspace p^{s-1}\).

Finding the non-invertible elements of \(\bbZ_{p^{s};\alpha}\), for instance, is therefore
equivalent to finding all the \((x_{1;0},\ldots,x_{1;s-1},x_{2;0},\ldots,x_{2;s-1}) \in
\bbZ^{2s}\) such that there exists \(k \in \bbZ\) such that
\begin{multline*}
( x_{1;0}\thinspace p^{0} + x_{1;1}\thinspace p^{1}+\cdots
+x_{1;s-1}\thinspace p^{s-1})^{2}+( \alpha_{0}\thinspace p^{0} +
\alpha_{1}\thinspace p^{1}+\cdots + \alpha_{s-1}\thinspace p^{s-1})
( x_{2;0}\thinspace p^{0} + x_{2;1}\thinspace p^{1}+\cdots
+x_{2;s-1}\thinspace p^{s-1})^{2}\\
\equiv k\thinspace p \modulus p^{s}
\end{multline*}

Solving this kind of equation is classical and elementary (\cite{JonesJones(1998)}).

In the case of \(\bbZ_{p^{s};\alpha}\) it is even especially simple, as we see that in that
particular case our equation is equivalent to \(x_{1;0}^{2}+
\alpha_{0}\thinspace x_{2;0}^{2} \equiv 0 \modulus p\), which we have already
solved in \ref{sbs:cagaco}.

In brief, there are
\(p+\left(\frac{a_{0}}{p}\right)\thinspace(p-1)\thinspace(-1)^{\frac{p-1}{2}}\)
possible choices for \((x_{1;0},x_{2;0})\), with \(\left(\frac{x}{p}\right)\)
denoting the Legendre symbol of \(x\) with respect to \(p\), and \(p\) choices for every
other variables. Hence \(\cardinal(\bbZ_{p^{s};\alpha}^{*}) = p^{2s}-
p^{2(s-1)}\left[p+\left(\frac{a_{0}}{p}\right)\thinspace(p-1)\thinspace(-1)^{\frac{p-1}{2}}\right]\).

In the case of \(\bbZ_{p^{s};\alpha,\beta}\) (and
\(\bbZ_{p^{s};\alpha,\beta,\gamma}\)), the method is identical, and the results are
summarized in the table below.

Let us now turn to the case of the unimodulars. By definition, this means we look for the
elements for which the Cayley norm is equal to \(\bar{1} \in \bbZ_{p^{s}}\).

let us first determine the number of unimodulars in \(\bbZ_{p^{s}}\). By the same token
as above, we are led to solve the equation
\begin{equation}
( x_{1;0}\thinspace p^{0} + x_{1;1}\thinspace p^{1}+\cdots
+x_{1;s-1}\thinspace p^{s-1})^{2} \equiv 1
\modulus p^{s}\tag{U}\label{eq:unimodeq}
\end{equation}

This is of course equivalent to the following expanded equation
\begin{equation}
\sum_{j=0}^{s-1}\left(\sum_{k=0}^{j}x_{1;k}\thinspace x_{1;j-k}\right)\thinspace p^{j}
\equiv 1 \modulus p^{s}\tag{U'}\label{eq:unimodeqex}
\end{equation}

The first thing we note is that (\ref{eq:unimodeqex}) implies
\begin{equation}
x_{1;0}^{2} \equiv 1 \modulus p^{1}\tag{U'1}\label{eq:unimodeqex1}
\end{equation}

 If \(s = 1\) these equations are actually equivalent, of course, but otherwise
 (\ref{eq:unimodeqex}) is equivalent to the conjunction of itself and
 (\ref{eq:unimodeqex1}).

We know (by \ref{ssc:galeq}) that the solutions to (\ref{eq:unimodeqex1}) are
\(x_{1;0} = 1\) and \(x_{1;0} = p-1\), both of which happen to be invertible in
\(\bbZ_{p}\). If \(s = 1\), we are done, as noted above, otherwise,
(\ref{eq:unimodeqex1}) is a constraint which our solutions must verify.

The second thing we note is that (\ref{eq:unimodeqex}) also implies
\begin{equation}
[x_{1;0}^{2}]\thinspace p^{0}+[2\thinspace x_{1;0}\thinspace
x_{1;1}]\thinspace p^{1} \equiv 1
\modulus p^{2}\tag{U'2}\label{eq:unimodeqex2}
\end{equation}

As above, if \(s = 2\), (\ref{eq:unimodeqex}) is equivalent to (\ref{eq:unimodeqex2}),
but even if such is not the case\(\ref{eq:unimodeqex}\) is equivalent to
[(\ref{eq:unimodeqex}) and (\ref{eq:unimodeqex1}) and (\ref{eq:unimodeqex2})].
We note also that (\ref{eq:unimodeqex1}) implies that
\(x_{1;0}^{2} \equiv 1+k_{0}\thinspace p^{1} \modulus p^{2}\), for some integer
\(k_{0}\) dependent only on \(x_{1;0}\). But this means that [(\ref{eq:unimodeqex1})
and (\ref{eq:unimodeqex2})] is equivalent to the conjunction of (\ref{eq:unimodeqex1})
and \(2\thinspace x_{1;0}\thinspace x_{1;1} \equiv k_{0} \modulus p^{1}\).
Since \(x_{1;0}\) always is invertible in \(\bbZ_{p}\) for the solutions of
(\ref{eq:unimodeqex}), this implies a single possible value for \(x_{1;1}\) given a value
for \(x_{1;0}\).

We may iterate this reasoning until we have reached the value of \(s\). At each step save
the first, the value of \(x_{1;i}\) can be seen to be uniquely determined by the values of
\(x_{1;0}\), \ldots, \(x_{1;i-1}\). Finally,
\(\cardinal(\mathfrak{U}_{\bbZ_{p^{s}}}) = 2\).

We briefly sketch the computation of
\(\cardinal(\mathfrak{U}_{\bbZ_{p^{s};\alpha}})\), which is essentially similar to
the above, albeit yet more computationally unpleasant.
\(\cardinal(\mathfrak{U}_{\bbZ_{p^{s};\alpha,\beta}})\) and
\(\cardinal(\mathfrak{U}_{\bbZ_{p^{s};\alpha,\beta,\gamma}})\) appear in the
table below.

To compute that number, we have to solve
\begin{multline}
( x_{1;0}\thinspace p^{0} + x_{1;1}\thinspace p^{1}+\cdots
+x_{1;s-1}\thinspace p^{s-1})^{2}+( \alpha_{0}\thinspace p^{0} +
\alpha_{1}\thinspace p^{1}+\cdots + \alpha_{s-1}\thinspace p^{s-1})
( x_{2;0}\thinspace p^{0} + x_{2;1}\thinspace p^{1}+\cdots
+x_{2;s-1}\thinspace p^{s-1})^{2}\\
\equiv 1 \modulus p^{s}\tag{UU}\label{eq:uunimod}
\end{multline}

The first step, as above, it to note that (\ref{eq:uunimod}) is equivalent to the conjunction
of itself and
\begin{equation}
x_{1;0}^{2}+\alpha_{0}\thinspace x_{2;0}^{2} \equiv 1
\modulus p\tag{UU1}\label{eq:uunimod1}
\end{equation}

We have solved this equation too, in \ref{sbs:cagaco}, and we know it has
\(p-\left(\frac{\alpha_{0}}{p}\right)(-1)^{\frac{p-1}{2}}\) solutions, all of which
evidently differ from \((0,0)\).

The second step, which is only necessary if \( s \geqslant 2\), is to see that
(\ref{eq:uunimod}) is equivalent to the conjunction of itself, (\ref{eq:uunimod1}) and
\begin{equation}
[x_{1;0}^{2}+\alpha_{0}\thinspace x_{2;0}^{2}]\thinspace p^{0}+
[2\thinspace x_{1;0}\thinspace x_{1;1}+
2\alpha_{0}\thinspace x_{2;0}\thinspace x_{2;1}+
\alpha_{1}\thinspace x_{2;0}^{2}]\thinspace p^{1}
\equiv 1 \modulus p^{2}\tag{UU2}\label{eq:uunimod2}
\end{equation}

Since we know that \(x_{1;0}^{2}+\alpha_{0}\thinspace x_{2;0}^{2} \equiv
1 + k_{0}\thinspace p^{1} \modulus p^{2}\) for some integer \(k_{0}\) dependant
only upon \(x_{1;0}\) and \(x_{2;0}\), we see that what we actually have to solve is
\(2\thinspace x_{1;0}\thinspace x_{1;1}+
2\alpha_{0}\thinspace x_{2;0}\thinspace x_{2;1}+
\alpha_{1}\thinspace x_{2;0}^{2} \equiv
k_{0} - \alpha_{1}\thinspace x_{2;0}^{2}\). Hence for each \((x_{1;0},x_{2;0})\)
solution of (\ref{eq:uunimod1}) there are \(p\) choices for \((x_{1;1},x_{2;1})\).

We may iterate this reasoning until we have reached the value of \(s\). At each step save
the first, there are \(p\) values of \((x_{1;i},x_{2;i})\) which are uniquely determined by
the values of \(x_{1;0}\), \ldots, \(x_{1;i-1}\), \(x_{2;0}\), \ldots, \(x_{2;i-1}\).
Finally, \(\cardinal(\mathfrak{U}_{\bbZ_{p^{s};\alpha}}) =
p^{s-1} \left[p-\left(\frac{\alpha_{0}}{p}\right)(-1)^{\frac{p-1}{2}}\right]\).

We now assume \(p = 2\), where things keep getting ugly. We will also assume, for brevity's
 sake, that \(\alpha \in \bbZ_{2^{s}}^{*}\), \(\beta \in \bbZ_{2^{s}}^{*}\),
\(\gamma \in \bbZ_{2^{s}}^{*}\)
. Note that 
\(\alpha \in \bbZ_{2^{s}}^{*}\) is equivalent to \(\alpha_{0} \equiv 1 \modulus 2\),
and likewise for \(\beta\) and \(\gamma\). Recall also we have already shown that
\(\cardinal(\bbZ_{2^{s}}^{*}) = 2^{s-1}\).

We again look for elements of \(\bbZ_{2^{s};\alpha}\) \emph{etc.} which are \emph{not}
invertible. As above, we look for them in by identifying members of their class in
\(\{0,\ldots,2^{s}-1\}\), and writing them in base \(2\). However in this case the
expansion of the equation is different, precisely because the number \(2\) which is the
base of expansion also is a coefficient of elements in the square of a sum!

For \(\bbZ_{2^{s};\alpha}\), we have to solve
\begin{multline*}
( x_{1;0}\thinspace 2^{0} + x_{1;1}\thinspace 2^{1}+\cdots
+x_{1;s-1}\thinspace 2^{s-1})^{2}+( \alpha_{0}\thinspace 2^{0} +
\alpha_{1}\thinspace 2^{1}+\cdots + \alpha_{s-1}\thinspace2p^{s-1})
( x_{2;0}\thinspace 2^{0} + x_{2;1}\thinspace 2^{1}+\cdots
+x_{2;s-1}\thinspace 2^{s-1})^{2}\\
\equiv k\thinspace 2 \modulus 2^{s}
\end{multline*}
for all integers \(k\).

Expanding, we find that this is equivalent to
\(x_{1;0}^{2}+\alpha_{0}\thinspace x_{2;0}^{2} \equiv 0 \modulus 2\), and
we have determined that (since we want to consider only
\(\alpha \in \bbZ_{2^{s}}^{*}\)) \(\alpha_{0} = 1\). Furthermore, we know that
\([x \mapsto x^{2}]\) is a bijection on \(\bbZ_{2} = \GaloisField(2)\). Hence there
just two solutions for \((x_{1;0},x_{2;0})\), and \(2\) solution for every other unknown.
Hence \(\cardinal(\bbZ_{2^{s};\alpha}^{*}) = 2^{2s-1}\).

We proceed likewise for \(\bbZ_{2^{s};\alpha,\beta}^{*}\) and
\(\bbZ_{2^{s};\alpha,\beta,\gamma}^{*}\) and give the relevant values in the table
below.

We now turn to the unimodulars.

In the case of \(\bbZ_{2}\), the relevant numbers have been computed already in
\ref{sbs:cagaco}.

Since \(\bbZ_{2^{2}}\) has just four elements, we list them explicitly and discover
that there are exactly \(2\) unimodulars. We note that the Cayley norm of an element of
\(\bbZ_{2^{2}}\) is always either \(\bar{0}\) or \(\bar{1}\) (and in particular all
invertibles are unimodular), and never either \(\bar{2}\) or \(\bar{3}\) (the only
remaining possibilities).

If \((x,y) \in \bbZ_{2^{2};\alpha}\), we know that (with our above conventions)
\(\CayleyNorm{\bbZ_{2^{2};\alpha}}((x,y))=
\CayleyNorm{\bbZ_{2^{2}}}(x)+
\overline{(\alpha_{0}+\alpha_{1}\thinspace 2)}\CayleyNorm{\bbZ_{2^{2}}}(y)\),
and we have assumed here that \(\alpha_{0} = 1\).

Therefore, if \(\alpha_{1} = 0\)
the Cayley norm of an element of \(\bbZ_{2^{2};\alpha}\) can only take the value
\(\bar{0}\), \(\bar{1}\) or \(\bar{2}\), hence here too unimodulars and invertible
are the same elements (and so \(\cardinal(\mathfrak{U}_{\bbZ_{2^{2};\alpha}}) = 8\)).

On the other hand, if \(\alpha_{1} = 1\) (the only other possible value), then we have
\(\CayleyNorm{\bbZ_{2^{2};\alpha}}((\bar{1},\bar{0})) = \bar{1}\) and 
\(\CayleyNorm{\bbZ_{2^{2};\alpha}}((\bar{0},\bar{1})) = \bar{3}\) and we see that
\(\bbZ_{2^{2}}^{*} \subset \{\bar{1},\bar{3}\} \subset
\mathfrak{R}_{\bbZ_{2^{2};\alpha}} \subset \bbZ_{2^{2}}^{*}\); we now invoke
\thref{prp:unires} (more precisely (\ref{equ:P3})):
\(\cardinal(\bbZ_{2^{2};\alpha}^{*}) =
\cardinal(\mathfrak{U}_{\bbZ_{2^{2};\alpha}})\thinspace
\cardinal(\mathfrak{R}_{\bbZ_{2^{2};\alpha}})\), which here means
\(8 = \cardinal(\mathfrak{U}_{\bbZ_{2^{2};\alpha}})\times 2\).

In a similar way,
\(\CayleyNorm{\bbZ_{2^{2};\alpha,\beta}}((\bar{1},\bar{0},\bar{0},\bar{0})) =
\bar{1}\), so \(\{\bar{1}\} \in \mathfrak{R}_{\bbZ_{2^{2};\alpha,\beta}}\).
If \(\alpha_{1} = \beta_{1} = 0\), then
\(\CayleyNorm{\bbZ_{2^{2};\alpha,\beta}}((\bar{0},\bar{1},\bar{1},\bar{1})) =
\bar{3}\). If \(\alpha_{1} = 0\) and \(\beta_{1} = 1\) then
\(\CayleyNorm{\bbZ_{2^{2};\alpha,\beta}}((\bar{0},\bar{0},\bar{0},\bar{1})) =
\bar{3}\). Finally, if \(\alpha_{1} = 1\) and irrespective of the value of \(\beta_{1}\)
\(\CayleyNorm{\bbZ_{2^{2};\alpha,\beta}}((\bar{0},\bar{1},\bar{0},\bar{0})) =
\bar{3}\). So here too we have
\(\mathfrak{R}_{\bbZ_{2^{2};\alpha,\beta}} = \bbZ_{2^{2}}^{*}\). We reason in
just the same way for \(\mathfrak{U}_{\bbZ_{2^{2};\alpha,\beta,\gamma}}\).

Assume now that \(s \geqslant 3\).

The unimodulars of \(\bbZ_{2^{s}}\) are a well-studied topic, and we find, in \emph{e.g.}
\cite{JonesJones(1998)}, that \(\cardinal(\mathfrak{U}_{\bbZ_{2^{s}}}) = 4\) and
that \(\mathfrak{R}_{\bbZ_{2^{s}}} =
\{\bar{a}\thinspace|\thinspace a \in \bbZ, a \equiv 1 \modulus 8\} =
\{\bar{1}+k\thinspace\bar{8}\thinspace|\thinspace k \in \bbZ\}\).
On the other hand, an element of \(\bbZ_{2^{s}}\) is not invertible if and only if it is of
the form \(k\thinspace\bar{2}\), with \(0 \leqslant k < 2^{s-1}\), 
and these have squares of the form \(k^{2}\thinspace\bar{4}\), and these in turn are all
of the form \(l\thinspace\bar{16}\) (when \(k\) is even) or
\(\bar{4}+l\thinspace\bar{32}\) (when \(k\) is odd), for some
(though perhaps not every) integer \(l\).

Given \(\alpha \in \bbZ_{2^{s}}^{*}\), by considering elements of the form
\((x,\bar{0})\) we see that \(\mathfrak{R}_{\bbZ_{2^{s};\alpha}} \supset
\mathfrak{R}_{\bbZ_{2^{s}}}\).

Assume \(\alpha_{1} = \alpha_{2} = 0\) (and still \(\alpha_{0} = 1\)). We thus
have \(\alpha = \bar{1}+k_{\alpha}\thinspace\bar{8}\) for some integer
\(k_{\alpha}\), and this implies, as stated above, that
\(\alpha \in \mathfrak{R}_{\bbZ_{2^{s}}}\).

Since \(\mathfrak{R}_{\bbZ_{2^{s}}}\) is a group, there exists
\(z_{\alpha} \in \mathfrak{R}_{\bbZ_{2^{s}}}\) such that
\(\alpha \times z_{\alpha} = \bar{1}\), and by definition of
\(\mathfrak{R}_{\bbZ_{2^{s}}}\), there exists \(y_{\alpha} \in \bbZ_{2^{s}}^{*}\)
such that \(y_{\alpha}^{2} = z_{\alpha}\), and thus
\(\alpha \times y_{\alpha}^{2} = \bar{1}\). By considering elements of the form
\((x,2y_{\alpha})\), with \(x \in \bbZ_{2^{s}}^{*}\), we see that
\(\mathfrak{R}_{\bbZ_{2^{s};\alpha}}\) also contains all elements of the form
\(\bar{5}+k\thinspace\bar{8}\).

However, \(\mathfrak{R}_{\bbZ_{2^{s};\alpha}}\) contains no element of the form
\(\bar{3}+k\thinspace\bar{8}\) or \(\bar{7}+k\thinspace\bar{8}\).
Indeed, if we look at the four possible situations for \((x,y) \in \bbZ_{2^{s};\alpha}\),
we see that if \(x \not\in \bbZ_{2^{s}}^{*}\) and \(y \not\in \bbZ_{2^{s}}^{*}\)
or \(x \in \bbZ_{2^{s}}^{*}\) and \(y \in \bbZ_{2^{s}}^{*}\), then
\(\CayleyNorm{\bbZ_{2^{s};\alpha}}((x,y))\) is a multiple of \(\bar{2}\) (and
\((x,y) \not\in \bbZ_{2^{s};\alpha}^{*}\));
if \(x \in \bbZ_{2^{s}}^{*}\) and \(y \not\in \bbZ_{2^{s}}^{*}\) or
\(x \not\in \bbZ_{2^{s}}^{*}\) and \(y \in \bbZ_{2^{s}}^{*}\) then
\(\CayleyNorm{\bbZ_{2^{s};\alpha}}((x,y))\) is of the form
\(\bar{1}+k^{2}\thinspace\bar{4}+l\thinspace\bar{8}\) for some integers \(k\) and
\(l\), which means that they are of the form \(\bar{1}+m\thinspace\bar{8}\) or
\(\bar{5}+m\thinspace\bar{8}\) for some integer \(m\), according to whether \(k\) is
even or odd.

Therefore \(\cardinal(\mathfrak{R}_{\bbZ_{2^{s};\alpha}}) =
\frac{1}{2} \cardinal( \bbZ_{2^{s}}^{*})\), and
\(\cardinal(\mathfrak{U}_{\bbZ_{2^{s};\alpha}}) = 2^{s+1}\).

Assume now \(\alpha_{1} = 0\), \(\alpha_{2} = 1\) (and still \(\alpha_{0} = 1\)). We thus
have \(\alpha = \bar{5}+k_{\alpha}\thinspace\bar{8}\) for some integer
\(k_{\alpha}\).

If \(s = 3\), by listing all the elements explicitly, we see the the Cayley norm can
reach the value \(\bar{5}\) (and \(\bar{1}\), of course), but not the values \(\bar{3}\)
and \(\bar{7}\). Hence here \(\cardinal(\mathfrak{U}_{\bbZ_{2^{3};\alpha}}) =
2^{3+1}\).

Assume therefore that \(s \geqslant 4\).

By using the projection 
\(\pi_{2^{s},2^{s-3}}:\bbZ_{2^{s}} \rightarrow \bbZ_{2^{s-3}}\) we see that
\(\pi_{2^{s},2^{s-3}}(\alpha)\) is invertible in \(\bbZ_{2^{s-3}}\) and that we can,
for any value of \(t\) in \(\{0,\ldots,2^{s-3}-1\}\), solve the equation
\(l(5+k_{\alpha}\thinspace 8) \equiv t - k_{\alpha} \modulus 2^{s-3}\), yielding
\(l_{t}\). For each \(l_{t}\), we know that \(\bar{1}+l_{t}\thinspace\bar{8} \in
\mathfrak{R}_{\bbZ_{2^{s};\alpha}}\), and therefore that there exists \(y_{t} \in
\bbZ_{2^{s}}^{*}\) such that \(y_{t}^{2} = bar{1}+l_{t}\thinspace\bar{8}\).
But then \(\alpha\thinspace\times y_{t}^{2} = \bar{5}+t\thinspace\bar{8}\), so
all elements of the form \(\bar{5}+t\thinspace\bar{8}\) are reached by the Cayley norm.

A discussion similar to the above shows that no value of the form
\(\bar{3}+k\thinspace\bar{8}\) or \(\bar{7}+k\thinspace\bar{8}\) can be reached
by the Cayley norm. Threfore, as above,
\(\cardinal(\mathfrak{U}_{\bbZ_{2^{s};\alpha}}) = 2^{s+1}\).

Assume now \(\alpha_{1} = 1\), \(\alpha_{2} = 0\) (and still \(\alpha_{0} = 1\)). We thus
have \(\alpha = \bar{3}+k_{\alpha}\thinspace\bar{8}\) for some integer
\(k_{\alpha}\).

We know all values of the form \(\bar{1}+k\thinspace\bar{8}\) can be reached by the
Cayley norm, and by considering elements of the form \((x,\bar{2})\) with
\(x \in \bbZ_{2^{s}}\) we see all values of the form  \(\bar{5}+k\thinspace\bar{8}\)
can also be reached. If \(s = 3\) we see that the values \(\bar{3}\) and \(bar{7}\) are
reached (by \((\bar{0},\bar{1})\) and \((\bar{2},\bar{1})\) respectively, for instance).
If \(s \geqslant 4\), solving the congruence equations
\(l(3+k_{\alpha}\thinspace 8) \equiv t - k_{\alpha} \modulus 2^{s-3}\) proves as
above that all values of the form \(\bar{5}+k\thinspace\bar{8}\) are reached too.
Since \(\mathfrak{R}_{\bbZ_{2^{s};\alpha}}\) is a subgroup of \(\bbZ_{2^{s}}^{*}\)
we see that we actually have \(\mathfrak{R}_{\bbZ_{2^{s};\alpha}} =
\bbZ_{2^{s}}^{*}\). Therefore \(\cardinal(\mathfrak{R}_{\bbZ_{2^{s};\alpha}}) =
2^{s}\).

The final case, \(\alpha_{1} = 1\), \(\alpha_{2} = 1\) (and still \(\alpha_{0} = 1\)),
for which we have have \(\alpha = \bar{7}+k_{\alpha}\thinspace\bar{8}\) for some
integer \(k_{\alpha}\) is entirely similar to the above, and we thus find that
\(\cardinal(\mathfrak{R}_{\bbZ_{2^{s};\alpha}}) = 2^{s}\).

We will gloss over the unimodulars of \(\bbZ_{2^{2};\alpha,\beta}\),
\(\bbZ_{2^{2};\alpha,\beta,\gamma}\), \(\bbZ_{2^{3};\alpha,\beta}\) and
\(\bbZ_{2^{3};\alpha,\beta,\gamma}\); to find the values we want we either investigate
the explicit list of the elements or we adapt what follows. We therefore assume
\(s \geqslant 4\).

If either \(\alpha = \bar{3}+k\thinspace\bar{4}\) or
\(\beta = \bar{3}+k\thinspace\bar{4}\), by considering elements of the form
\((x,y,\bar{0},\bar{0})\) or \((x,\bar{0},y,\bar{0})\) and using the methods we have
developed above, we see that actually \(\mathfrak{R}_{\bbZ_{2^{s};\alpha,\beta}} =
\bbZ_{2^{s}}^{*}\). Therefore
\(\cardinal(\mathfrak{R}_{\bbZ_{2^{s};\alpha,\beta}}) = 2^{3s}\).

If both \(\alpha = \bar{1}+k\thinspace\bar{4}\) and
\(\beta = \bar{1}+l\thinspace\bar{4}\), let \(u \in \bbZ_{2^{s}}\) and
\(v \in \bbZ_{2^{s}}\) such that
\begin{align*}
&u =
\begin{cases}
\bar{1} & \text{if \(\beta = 1 + l'\thinspace\bar{8}\)}\\
\bar{5} & \text{if \(\beta = 1 + l'\thinspace\bar{8}\)}
\end{cases}
&v =
\begin{cases}
\bar{1} & \text{if \(\alpha = 1 + k'\thinspace\bar{8}\)}\\
\bar{5} & \text{if \(\alpha = 1 + k'\thinspace\bar{8}\)}
\end{cases}
\end{align*}

Then
\begin{equation*}
\CayleyNorm{\bbZ_{2^{s};\alpha,\beta}}((x,y,u,uv)) = x^{2}+y^{2}+
\begin{cases}
\bar{2}\\
\bar{6}\\
\bar{10}
\end{cases}
\end{equation*}

We can thus see that \(\bar{3}\), \(\bar{7}\) or \(\bar{13}\) is reached, and reasoning
as above, \(\mathfrak{R}_{\bbZ_{2^{s};\alpha,\beta}} = \bbZ_{2^{s}}^{*}\).
Therefore
\(\cardinal(\mathfrak{R}_{\bbZ_{2^{s};\alpha,\beta}}) = 2^{3s}\).

Since by considering elements of the form \((x,y,z,t,\bar{0},\bar{0},\bar{0},\bar{0})\)
the Cayley norm reaches all of \(\bbZ_{2^{s}}^{*}\), we can immediately see that
\(\cardinal(\mathfrak{R}_{\bbZ_{2^{s};\alpha,\beta,\gamma}}) = 2^{7s}\).

\begin{center}
{\renewcommand{\arraystretch}{1.5}
\begin{tabular}{||l||c|c|c|c||}
\hhline{~|t:====:t|}
\multicolumn{1}{c||}{(\(\alpha\beta\gamma \in \bbZ_{p^{s}}^{*}\))} & \(p \neq 2\) & \multicolumn{3}{c||}{\(p = 2\)} \\[1mm]
\hhline{~~|---||}
\multicolumn{1}{c||}{} &\multicolumn{1}{c|}{} & \(s = 1\) & \(s = 2\) & \(s \geqslant 3\)\\[1mm]
\hhline{|t:=::=|===||}
\(\cardinal(\rsfsA^{*})\) & \(p^{s-1}(p-1)\) & \multicolumn{3}{c||}{\(2^{s-1}\)} \\[1mm]
\hhline{||-----||}
\(\cardinal(\bbZ_{p^{s}}^{*})\) & \(p^{s-1}(p-1)\) & \multicolumn{3}{c||}{\(2^{s-1}\)} \\[1mm]
\hhline{||-----||}
\(\cardinal(\bbZ_{p^{s};\alpha}^{*})\) & \(p^{2s}-p^{2(s-1)}\left[p+\left(\frac{a_{0}}{p}\right)\thinspace(p-1)\thinspace(-1)^{\frac{p-1}{2}}\right]\) & \multicolumn{3}{c||}{\(2^{2s-1}\)} \\[1mm]
\hhline{||-----||}
\(\cardinal(\bbZ_{p^{s};\alpha,\beta}^{*})\) & \(p^{4s} - p^{4(s-1)}\left[p^{3}+(p-1)\thinspace p\right]\) & \multicolumn{3}{c||}{\(2^{4s-1}\)} \\[1mm]
\hhline{||-----||}
\(\cardinal(\bbZ_{p^{s};\alpha,\beta,\gamma}^{*})\) & \(p^{8s} - p^{8(s-1)}\left[p^{7}+(p-1)\thinspace p^{3}\right]\) & \multicolumn{3}{c||}{\(2^{8s-1}\)} \\[1mm]
\hhline{||-----||}
\(\cardinal(\mathfrak{U}_{\bbZ_{p^{s}}})\) & \(2\) & \(1\) & \(2\) & \(4\) \\[1mm]
\hhline{||-----||}
\(\cardinal(\mathfrak{U}_{\bbZ_{p^{s};\alpha}})\) & \(p^{s-1} \left[p-\left(\frac{\alpha_{0}}{p}\right)(-1)^{\frac{p-1}{2}}\right]\) & \(2\) & \multicolumn{2}{c||}{\(2^{s+1-\alpha_{1}}\)} \\[1mm]
\hhline{||-----||}
\(\cardinal(\mathfrak{U}_{\bbZ_{p^{s};\alpha,\beta}})\) &\(p^{3(s-1)} \left[p^{3}-p\right]\)& \multicolumn{3}{c||}{\(2^{3s}\)} \\[1mm]
\hhline{||-----||}
\(\cardinal(\mathfrak{U}_{\bbZ_{p^{s};\alpha,\beta,\gamma}})\) & \(p^{7(s-1)} \left[p^{7}-p^{3}\right]\) & \multicolumn{3}{c||}{\(2^{7s}\)} \\[1mm]
\hhline{|b:=:b:====:b|}
\end{tabular}
}
\end{center}

\subsection{Further results}
\indent

We collect here some remarks and results about the objects we have studied in this section and
relatives thereof, which are too few to warrant a section of their own.

\begin{proposition}
Let \(p\) be an odd prime, \(k \geqslant 1\) an integer, \(q = p^{k}\) and
\(\alpha \in F_{q}^{*}\) such that
\(\rsfsX_{2}(\alpha)\thinspace(-1)^{\frac{q-1}{2}} = -1\). Then
\(F_{q;\alpha}^{*}\) is cyclic and in particular there is an element of order
\(q^{2}-1\). Furthermore, if \(z_{0}\) is a generator of \(F_{q;\alpha}^{*}\),
then \(\CayleyNorm{F_{q;\alpha}}(z_{0})\) is a generator of
\(\mathfrak{R}_{F_{q;\alpha}}\).
\end{proposition}

\begin{proof}
Under these hypotheses, \(F_{q;\alpha}\) is the finite field of cardinal \(q^{2}\),
hence its invertible elements are cyclic. Since \(\cardinal(F_{q;\alpha}^{*}) = q^{2}-1\)
this means that there is an element of order \(q^{2}-1\). The last assertion is trivial by
contraposition.
\end{proof}

\begin{proposition}
Let \(p\) be an odd prime, \(k \geqslant 1\) an integer, \(q = p^{k}\) and
\(\alpha \in F_{q}^{*}\) such that
\(\rsfsX_{2}(\alpha)\thinspace(-1)^{\frac{q-1}{2}} = +1\). Then the order of any
element of \(\in F_{q;\alpha}^{*}\) divides \(q-1\), and \(\in F_{q;\alpha}^{*}\)
is not cyclic; however \(\mathfrak{U}_{F_{q;\alpha}}\) is cyclic and isomorphic to
\(F_{q}^{*}\).
\end{proposition}

\begin{proof}
Under these hypotheses, \(F_{q;\alpha}^{*}\) , \(\mathfrak{U}_{F_{q;\alpha}}\) and
\(\mathfrak{R}_{F_{q;\alpha}}\) are abelian groups, and
\(\cardinal(\mathfrak{U}_{F_{q;\alpha}}) =
\cardinal(\mathfrak{R}_{F_{q;\alpha}}) = q-1\). Therefore, given any
\((n,u) \in \cardinal(\mathfrak{U}_{F_{q;\alpha}}) \vartimes
\cardinal(\mathfrak{R}_{F_{q;\alpha}})\), \((n,u)^{q-1} = (\dot{1},\dot{1})\),
with \(\dot{1}\) the neutral element for the multiplication in \(F_{q}\).
The first statement then results from \thref{prp:unires}. The second statement then results
from the fact that \(q-1 < (q-1)^{2} = \cardinal(F_{q;\alpha}^{*})\).

Note that \(\rsfsX_{2}(\alpha)\thinspace(-1)^{\frac{q-1}{2}} = +1\) is equivalent to
the existence of \(\omega \in  F_{q}^{*}\) such that \(\omega^{2} = -\alpha\).
Therefore \((x,y) \in \mathfrak{U}_{F_{q;\alpha}} \Leftrightarrow
(x+\omega y)(x-\omega y) = e\). Let now \(z = x+\omega y \in F_{q}^{*}\).
We compute immediately that \((x,y) = \dot{2}^{-1}(z+z^{-1},\omega^{-1}(z-z^{-1}))\),
with \(\dot{2} = \dot{1}+\dot{1}\) (which is invertible because \(q\) is odd). Better
yet, we immediately verify that
\([F_{q}^{*} \rightarrow \mathfrak{U}_{F_{q;\alpha}};
z \mapsto \dot{2}^{-1}(z+z^{-1},\omega^{-1}(z-z^{-1}))]\) is a group isomorphism.
\end{proof}

\pagebreak
\appendix
\appendixpage
\addappheadtotoc

\section{Some common and not-so-common algebraic structures}
\indent

Algebra is rife with named structures. Unfortunately, authors tend to disagree about
what a given name refers to precisely (always with good reason, but easily resulting in
confusion nonetheless; \cite{Rusin(1998)}), even for the many structures
found in most introductory-level material (\cite[...]{Lelong-FerrandA(1978)}), and even
for reference material (\cite{Bourbaki(A), Lang(1971), Pflugfelder(1990)}).
This appendix is therefore here to describe, precisely, the meaning we give to various
terms used in the main text of this document. We have included in this discussion a few other
familiar structures to serve as ``landmarks'', but this is \emph{not} a comprehensive
catalogue of all algebraic structures (not even of the most common ones).
We will endeavor to follow what we feel is current usage (\cite{Klein(2002)}),
though we take this opportunity to perform a few generalizations.
One particular emphasis of this appendix is the careful
treatment of associativity, in some (again, not \emph{all}) of its various shades.

\subsection{Structures involving only one law}
\subsubsection{A protean foundation: the magma}
\indent

Given a set \(S\), an \emph{internal law} on \(S\) (or simply a \emph{law} on \(S\))
is a function \(f\) defined on all of
\(S \vartimes S\) with values in \(S\). A \emph{magma}
is a structure \((S, f)\) where \(S\) is a (perhaps empty) set and \(f\) is a law on  \(S\).
A \emph{groupoid} is a non-empty magma.
Given a set \(T \subset S\), if \(f\) source-restricted to \(T \vartimes T\) takes its
values in \(T\), we call the structure \((T, f|_{T \vartimes T}^{T})\) a
\emph{submagma} of \((S, f)\). Given any (non-empty) collection
\((T_{i})_{i \in I}\) of subsets
of \(S\) such that \(f\) source-restricted to \(T_{i} \vartimes T_{i}\) takes its
values in \(T_{i}\) for all \(i \in I\), if we name \(U = \bigcap_{i \in I} T_{i}\),
we find that \((U, f|_{U \vartimes U}^{U})\) is a submagma of \((S,f)\) and of every
\((T_{i}, f|_{T_{i} \vartimes T_{i}}^{T_{i}})\). Therefore, given any (perhaps empty)
set \(X \subset S\), Therefore there exist a (perhaps empty) set \(Y \subset S\) such that
\(X \subset Y\), \((Y, f|_{Y \vartimes Y}^{Y})\) is a submagma of \((S,f)\),
and such that if \((T, f|_{T \vartimes T}^{T})\) is a submagma of \((S,f)\)
such that \(X \subset T\), then \(Y \subset T\) (in other words, \(Y\) is the
smallest submagma of \(S\) containing \(X\)); we call this the submagma
\emph{generated} by \(X\). Given a collection \((a)_{j \in J}\) of elements of \(S\),
we call the submagma \emph{generated} by the \(a_{j}\), the submagma generated by
\(\{a_{j}| j \in J\}\).

A magma \((S,f)\) is \emph{monogenic} (or \emph{cyclic}) if and only if there exists some
\(a \in S\) such that the submagma generated by \(a\) is the whole of \(S\).

A law \(f\) on \(S\) is\footnote{This definition generalizes that \cite{Bourbaki(A)}.
A definition equivalent to that of \cite{Bourbaki(A)} involves separating left and right
alternativity, and flexibility (though left and right alternativity together imply flexibility, 
\cite{Schafer(1966), Carmody(1988)}).} \emph{alternative} if and only
if \((\forall\thinspace x \in S)(\forall\thinspace y \in S)
(\forall\thinspace z \in S)\thickspace
[(x = y \text{ or } x = z \text{ or } y = z) \Rightarrow f(f(x,y),z) = f(x,f(y,z)]\).
If the law of a magma is alternative we say that the magma itself is \emph{alternative}.

A law \(f\) on \(S\) is \emph{associative} if and only if
\((\forall\thinspace x \in S)(\forall\thinspace y \in S)
(\forall\thinspace z \in S)\thickspace f(f(x,y),z) = f(x,f(y,z))\).
If the law of a magma is associative we say that the magma itself is \emph{associative}.

A law \(f\) on \(S\) is \emph{power-associative} if and only if for each \(a \in S\)
the submagma generated by \(a\) is associative.
If the law of a magma is power-associative we say that the magma itself is
\emph{power-associative}. An alternative law (or magma) need not be
power-associative\footnote{Consider:
\begin{center}
\begin{tabular}{c|cccc}
f & a & b & c & d \\
\hhline{-----}
a & a & a & a & a \\
b & a & d & a & c \\
c & a & a & c & a \\
d & a & c & a & a
\end{tabular}
\end{center}

In this case, the law is alternative, but the submagma generated by \(b\) is the whole of
\(\{a,b,c,d\}\) and the law is not associative on it because
\(f(b,f(d,c)) = f(b,a) = a \neq c = f(c,c) = f(f(b,d),c)\), hence the law not
power-associative.}, and a power-associative law
(or magma) need not be alternative\footnote{\label{ftn:powass}Consider:
\begin{center}
\begin{tabular}{c|ccc}
f & a & b & c \\
\hhline{----}
a & a & c & b \\
b & c & b & a \\
c & b & a & c
\end{tabular}
\end{center}

In this case the law is power-associative, as the submagma generated by any element
is reduced to that element, but the law is not alternative, as for instance
\(f(a,f(a,b)) = f(a,c) = b \neq c = f(a,b) = f(f(a,a),b)\).}.

It is perhaps surprising that
a monogenic magma need not be power-associative\footnote{\cite{Pflugfelder(1990)}
gives the following example:
\begin{center}
\begin{tabular}{c|ccccc}
f & a & b & c & d & e \\
\hhline{------}
a & a & b & c & d & e \\
b & b & a & e & c & d \\
c & c & e & d & b & a \\
d & d & c & a & e & b \\
e & e & d & b & a & c
\end{tabular}
\end{center}

In this case \(e\) generates the magma (a loop, actually, see further on for the definition),
but \(f(c,f(c,c)) = f(c,d) = b \neq a = f(d,c) = f(f(c,c),c)\).}! A power-associative magma
need not be monogenic, of course\footnote{As the example in footnote \ref{ftn:powass}
attests.}.

A law \(f\) on \(S\) is \emph{di-associative} if and only if for each couple of (not
necessarily distinct) elements of
\(S\), \(a_{1}\) and \(a_{2}\), the submagma generated by \(\{a_{1},a_{2}\}\)
is associative.
If the law of a magma is di-associative we say that the magma itself is
\emph{di-associative}. A di-associative magma is always both alternative and
power-associative, but an alternative magma need not be di-associative\footnote{Consider:
\begin{center}
\begin{tabular}{c|cccc}
f & a & b & c & d \\
\hhline{-----}
a & a & a & a & a \\
b & a & a & a & c \\
c & a & a & c & a \\
d & a & a & a & a
\end{tabular}
\end{center}

In this case, the law is alternative, but if \(a_{1} = b\) and \(a_{2} = d\), then the
generated submagma is \(\{a,b,c,d\}\), which is not associative, since
\(f(f(b,d),c) = a \neq c = f(b,f(d,c))\).}, and a power-associative magma need not be
di-associative\footnote{Consider:
\begin{center}
\begin{tabular}{c|ccc}
f & a & b & c \\
\hhline{----}
a & a & b & a \\
b & b & a & a \\
c & a & a & a
\end{tabular}
\end{center}

In this case the law is power-associative, but not di-associative because if \(a_{1} = a\)
and \(a_{2} = c\) then the submagma generated by these two element is the whole of
\(\{a,b,c\}\), and the law is not associative on it (because
\(f(a,f(c,b)) = f(a,a) = a \neq b = f(a,b) = f(f(a,c),b)\)).}. Likewise, an associative magma is always
di-associative, but a di-associative magma need not be associative\footnote{consider:
\begin{center}
\begin{tabular}{c|ccc}
f & a & b & c \\
\hhline{----}
a & a & a & c \\
b & a & b & b \\
c & a & b & c
\end{tabular}
\end{center}

In this case the law is di-associative, but is not associative because
\(f(a,f(b,c)) = a \neq c = f(f(a,b),c)\).}.

Given a law \(f\) on \(S\), an \(e_{\mathfrak{l}} \in S\) is a
\emph{left identity element} for the law
\(f\) if and only if \((\forall\thinspace x \in S)\thickspace
f(e_{\mathfrak{l}},x) = x\); and  an \(e_{\mathfrak{r}} \in S\) is a
\emph{right identity element} for the law
\(f\) if and only if \((\forall\thinspace x \in S)\thickspace
f(x,e_{\mathfrak{r}}) = x\). An \(e \in S\) is an \emph{identity element}
(or \emph{neutral element}) for the law \(f\) if and only if it is both a left identity element
for \(f\) and a right identity element for \(f\). A neutral element is necessarily unique.
A neutral element is frequently called a ``unit'' for the law, but we will not use this
terminology here due as this would conflict with other possible meanings of this word.
A \emph{magma with unit} is a magma whose law has a neutral element.

A law \(f\) on \(S\) is \emph{commutative} if and only if
\((\forall\thinspace x \in S)(\forall\thinspace y \in S)\thickspace
f(x,y) = f(y,x)\).
If the law of a magma is commutative we say that the magma itself is \emph{commutative}.

Given two magmas \((S,f)\) and \((S',f')\),
a \emph{magma (homo)morphism} is a function
\(\theta : S \rightarrow S'\) such that
\((\forall\thinspace (x,y) \in S^{2})\thickspace
\theta(f(x,y)) = f'(\theta(x),\theta(y))\). If  \((S,f)\) and \((S',f')\) are magmas with
units, whose neutral elements are respectively \(e\) and \(e'\), then a function
\(\theta : S \rightarrow S'\) is a \emph{(homo)morphism of magmas with unit}
if and only if it is a magma morphism and \(\theta(e) = e'\). A surjective morphism is
called an \emph{epimorphism}, and an injective one is called a \emph{monomorphism}.
An \emph{endomorphism} is a morphism from one set to itself. A function
\(\theta : S \rightarrow S'\) is an \emph{isomorphism} if and only if it is a
homomorphism, it is bijective and the reciprocal function
\(\theta^{-1} : S' \rightarrow S\) is a homomorphism. An \emph{automorphism} is
an isomorphism from one set to itself.

If \(I\) is a set and \((S_{i},f_{i})_{i \in I}\) is a (non-empty) family of magmas, we
define a magma structure on the product set \(S\) by \(f : S \vartimes S ; 
((x_{i})_{i \in I},(x'_{i})_{i \in I}) \mapsto (f_{i}(x_{i},x'_{i}))_{i \in I}\).
If all the magmas in the family are commutative (respectively power-associative, alternative,
di-associative, associative, with unit), then the product magma is likewise commutative
(respectively power-associative, alternative, di-associative, associative, with unit). Even if
all the magmas in the family are monogenic, however, the product clearly does not have to
be monogenic.

\subsubsection{Strongly associative subsets}
\indent

We present this interesting variation on the idea of associativity which we will use later on
in this appendix. This is essentially a simplified presentation of notions found in
\cite{Bourbaki(A)}.

Let \((S,f)\) be a magma, and \(T \subset S\). We will say \(T\) is a
\emph{strongly associative subset} of \(S\) if and only if
\begin{equation}
(\forall\thinspace (u,v,w) \in S^{3})\thickspace
[[[(u,v) \in T^{2}] \text{ or } [(u,w) \in T^{2}] \text{ or } [(v,w) \in T^{2}]]
\Rightarrow f(u,f(v,w)) = f(f(u,v),w)] \notag
\end{equation}

Notice that \(\emptyset\) is a strongly
associative subset of \(S\)!

Before coming to the crux of this segment, let us recall that an ordered set
\((E, \preccurlyeq)\), whose order may perhaps only be partial, is said to be
\emph{inductive} if and only if any totally ordered subset of \(E\) has a majorant
in \(E\) for \(\preccurlyeq\). The very important Zorn
theorem\footnote{Given the axiomatic framework of Zermelo-Fr{\ae}nkel, the axiom of
choice and the theorem of Zorn are actually equivalent.} then states that
every inductive set has maximal elements. We can now state the following simple proposition:

\begin{proposition}[Inductivity of the Strongly Associative]\label{prop:indstrass}
Let \((S,f)\) be a magma, and \(\rsfsE\) the set of all strongly associative subsets of
\(S\); then \((\rsfsE, \subset)\) is inductive.
\end{proposition}

\begin{proof}
Let \(\rsfsF \subset \rsfsE\), totally ordered by \(\subset\), and let
\(U = \bigcup_{T \in \rsfsF} T\), then \(U\) is a majorant of \(\rsfsF\) in
\(\rsfsP(S)\) (the set of subsets of \(S\)).

If \(U = \emptyset\), then \(U\) is strongly associative.

If \(U \neq \emptyset\), let \((u,v,w) \in S^{3}\). Assume that \(u \in U\),
\(v \in U\). Then \((\exists\thinspace F_{u} \in \rsfsF)\thickspace u \in F_{u}\),
\((\exists\thinspace F_{v} \in \rsfsF)\thickspace v \in F_{v}\). Since
\((\rsfsF, \subset)\) is totally ordered by hypothesis, then either
\(F_{u} \subset F_{v}\) or  \(F_{v} \subset F_{u}\), so \(u \in F_{\alpha}\) and
\(v \in F_{\alpha}\) with \(F_{\alpha} = F_{u}\) or \(F_{\alpha} = F_{v}\) as
the case may be. But \(\rsfsF \subset \rsfsE\) and \(F_{\alpha} \in \rsfsF\), so
\(F_{\alpha} \in \rsfsE\), and therefore \(f(u,f(v,w)) = f(f(u,v),w)\). The other
possible branches of assumption are dealt with in the same manner. Finally \(U\) is strongly
associative in this case too.
\end{proof}

\subsubsection{Beyond the magma}
\indent

A \emph{quasigroup} is a groupoid \((S,f)\) such that for all \(a \in S\) the functions
\(\mathfrak{l}_{a} = [S \rightarrow S\thinspace;\thinspace x \mapsto f(a,x)]\)
(left translation, or multiplication, by \(a\)) and
\(\mathfrak{r}_{a} = [S \rightarrow S\thinspace;\thinspace x \mapsto f(x,a)]\)
(right translation, or multiplication, by \(a\)) are both
bijections\footnote{\label{ftn:quasigroup}Which \emph{does not} imply the existence of
either left or right inverses, or even that of a neutral element! Consider:
\begin{center}
\begin{tabular}{c|ccc}
f & a & b & c \\
\hhline{----}
a & b & c & a \\
b & a & b & c \\
c & c & a & b
\end{tabular}
\end{center}}.
We define \emph{subquasigroups} as non-empty submagmas of a quasigroup.
An intersection of quasigroups may be empty\footnote{\cite{Pflugfelder(1990)}
gives the following example:
\begin{center}
\begin{tabular}{c|ccc}
f & a & b & c \\
\hhline{----}
a & a & c & b \\
b & c & b & a \\
c & b & a & c
\end{tabular}
\end{center}}, but when it is
not, it is a quasigroup, so we may define the subquasigroup
\emph{generated} by a (non-empty) subset or a (non-empty) collection of elements.
A quasigroup is \emph{commutative}
(respectively \emph{monogenic}, \emph{alternative}, \emph{associative}) if, as a magma, it
is such. It is \emph{power-associative}, respectively \emph{di-associative}, if the
subquasigroup generated respectively by any one or any two (not necessarily distinct)
elements is associative.

A quasigroup \((Q,\star)\), even if devoid of a neutral element, is said to have the
\emph{left inverse property} (or \emph{L.I.P.} for short), or to be a \emph{L.I.P.
quasigroup} if and only if there exists a
bijection \(J_{\lambda}:Q \rightarrow Q; a \mapsto a^{\lambda}\) such that
\((\forall\thinspace a \in Q)(\forall\thinspace x \in Q)\thickspace
a^{\lambda} \star (a \star x) = x\), and is said to have the \emph{right inverse property}
(or \emph{R.I.P} for short), or to be a \emph{R.I.P. quasigroup} if and only if there exists a
bijection \(J_{\rho}:Q \rightarrow Q; a \mapsto a^{\rho}\) such that
\((\forall\thinspace a \in Q)(\forall\thinspace x \in Q)\thickspace
(x \star a) \star a^{\rho}= x\). A quasigroup which has both the L.I.P. property and the
R.I.P. property is said to have the \emph{inverse property} (or \emph{I.P.} for short), or to
be an \emph{I.P. quasigroup}.

A \emph{loop} is a quasigroup whose law has a neutral element. A quasigroup need not
be a loop\footnote{Consider:
\begin{center}
\begin{tabular}{c|ccc}
f & a & b & c \\
\hhline{----}
a & b & a & c \\
b & a & c & b \\
c & c & b & a
\end{tabular}
\end{center}}.
In a loop \((L,\times)\) therefore,
every element \(a \in L\) has both a right inverse \(a^{\lambda}\) and a left inverse
\(a^{\rho}\) (which are necessarily unique), but it may
happen\footnote{\cite{Pflugfelder(1990)} gives the following example:
\begin{center}
\begin{tabular}{c|ccccc}
f & a & b & c & d & e \\
\hhline{------}
a & a & b & c & d & e \\
b & b & a & e & c & d \\
c & c & d & a & e & b \\
d & d & e & b & a & c \\
e & e & c & d & b & a
\end{tabular}
\end{center}

In this case \((\forall\thinspace x \in L)\thickspace
x^{\lambda} = x^{\rho} = x\), but \(f(f(c,d),d) \neq c\) and
\(f(c,f(c,d)) \neq d\).} that there exists
\(b \in L\) such that \(a^{\lambda} \times (a \times b) \neq b\) or
\((b \times a)\times a^{\rho} \neq b\).
We define \emph{subloops} in a like manner as subquasigroups.
Note that one can prove\footnote{Let \((L,\times)\) be a loop with neutral element
\(e\), and \((H,\times)\) be a subloop. Then there exists \(a \in H\) such that
\(a \times a = a\), but as \(H \subset L\), \(a \in L\) and
\(a \times a = a = a \times e\).}
that the neutral element of a loop is included in every subloop
and is  is necessarily identical to that subloop's neutral element.
The intersection of subloops being a loop, we define likewise the subloop
\emph{generated} by a subset or a collection of elements.
A loop is \emph{commutative}
(respectively \emph{monogenic}, \emph{alternative}, \emph{L.I.P.}, \emph{R.I.P},
\emph{I.P.}) if, as a quasigroup, it is such. It is \emph{power-associative}, respectively
\emph{di-associative}, if the subloop generated respectively by any one or any two (not
necessarily distinct) elements is associative.

Note that in an I.P. loop, \(J_{\lambda} = J_{\rho}\), and therefore every element has an
inverse.

Given a loop \((L,\times)\), the following properties are equivalent
(\cite{Pflugfelder(1990)}):
\begin{align}
(\forall\thinspace x \in L)(\forall\thinspace y \in L)
(\forall\thinspace z \in L)\thickspace &(x \times y)\times(z \times x) =
[x \times(y \times z)]\times x \tag{M1}\label{eq:M1}\\
(\forall\thinspace x \in L)(\forall\thinspace y \in L)
(\forall\thinspace z \in L)\thickspace &(x \times y)\times(z \times x) =
x \times[(y \times z)\times x] \tag{M2}\label{eq:M2}\\
(\forall\thinspace x \in L)(\forall\thinspace y \in L)
(\forall\thinspace z \in L)\thickspace &[x \times(z \times x)]\times y =
x \times[z \times(x \times y)] \tag{M3}\label{eq:M3}\\
(\forall\thinspace x \in L)(\forall\thinspace y \in L)
(\forall\thinspace z \in L)\thickspace &[(x \times z)\times x]\times y =
x \times[z \times(x \times y)] \tag{M4}\label{eq:M4}\\
(\forall\thinspace x \in L)(\forall\thinspace y \in L)
(\forall\thinspace z \in L)\thickspace &[(y \times x)\times z]\times x =
y \times[x \times(z \times x)] \tag{M5}\label{eq:M5}
\end{align}

Any loop which verifies one (and therefore all) of the above properties is called a
\emph{Moufang loop}. A Moufang loop is \emph{commutative} (respectively
\emph{monogenic}) if, as a loop, it is such. We define \emph{sub-Moufang loops} in a like
manner as subloops. Any subloop of a Moufang loop is also a Moufang loop.

A Moufang loop always is
(\cite[Theorem~IV.1.4 and Corollary~IV.2.9]{Pflugfelder(1990)}) a di-associative I.P.
loop (which is equivalent to saying that all subloops generated by any two, not
necessarily distinct, elements always is a group; this is a corollary of Moufang's Theorem).

A \emph{monoid} is a magma with unit for which the law \(f\) is associative.
A monoid is \emph{commutative} (respectively \emph{monogenic}) if it is
commutative (respectively monogenic) as a magma.
A monoid need not be a quasigroup\footnote{Consider:
\begin{center}
\begin{tabular}{c|ccc}
f & a & b & c \\
\hhline{----}
a & c & a & a \\
b & a & b & c \\
c & a & c & c
\end{tabular}
\end{center}}, hence need not be a loop, and a quasigroup need not be a
monoid\footnote{See the example in
footnote~\ref{ftn:quasigroup}, which is not even alternative.}; a loop, even a
Moufang loop, need not be associative\footnote{\label{ftn:loopnotassociative}The
multiplication of the classical
octonions is an example of this; see also (\cite[page 89]{Pflugfelder(1990)}
for a reference to the smallest possible example of a non-associative Moufang loop.},
hence need not be a monoid either.

\pagebreak
A \emph{group} is a monoid in which every element has an inverse. One can prove
(\cite{Pflugfelder(1990)}) that being a group is equivalent to being a quasigroup
whose law is associative, or a loop whose law is associative, or a Moufang loop
whose law is associative\footnote{However, since
a loop, even a Moufang loop, needs not be associative, it needs not be a group; see
footnote \ref{ftn:loopnotassociative}.}.
We define \emph{subgroups} in a like manner as submagmas. As with loops, groups and
subgroups share their neutral element.
A group is \emph{commutative} (respectively \emph{monogenic}) if, as a monoid, it is such.
A commutative group is also called an \emph{abelian} group.

Morphisms (including isomorphims, \emph{etc.}) of quasigroups, loops, monoids and groups
are simply morphisms for the underlying magma (or magma with unit). While these notions
are quite interesting for groups, they are not so for poorer
structures (a quasigroup can perfectly be homomorphic to something which is \emph{not}
a quasigroup, and a loop can be homomorphic to something which is \emph{not}
a loop, see \cite[page~28]{Pflugfelder(1990)}), and are too strict a
tool to be very useful in classifying these later (a better tool for this is
\emph{isotopism}; see \cite{Pflugfelder(1990)}).

\subsubsection{Derivation graph for one law}\label{ssc:derone}
\indent

We present here a graph representing the ``is-a'' relationship for some of the most important
structures in this first part of the appendix.  We have abbreviated ``power-associative magma''
into ``P.-A. Magma'', ``di-associative magma'' into ``Di-A. Magma'', ``alternative magma''
into ``Alt. Magma'', ``associative magma'' into ``Ass. Magma'', ``magma with unit'' into
``Magma w. U.'' \emph{etc.}.

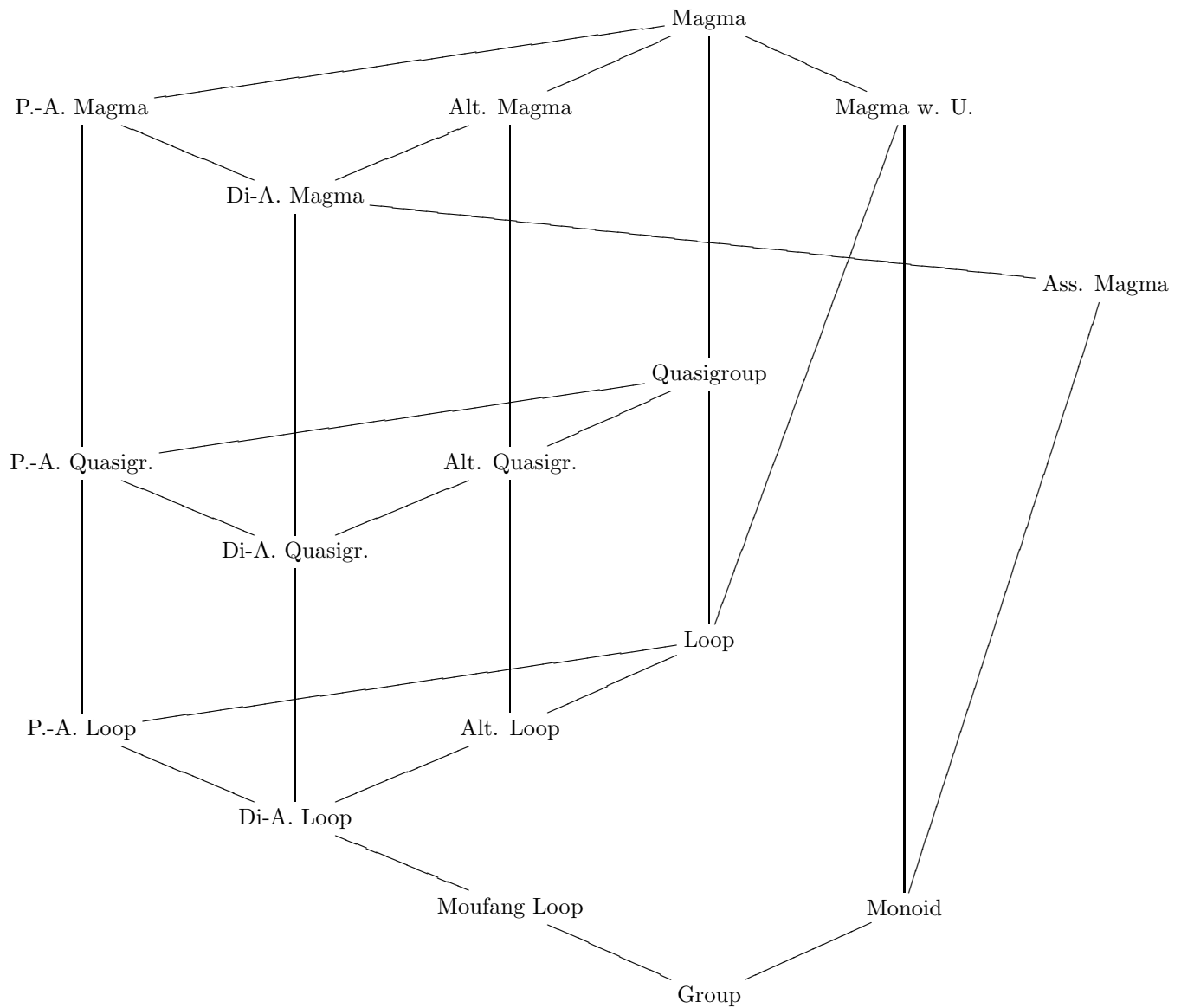
\begin{center}
\begin{math} 
\xymatrix{
 & & & \text{Magma} \ar@{-}[dddd]  \ar@{-}[llld]  \ar@{-}[ld]  \ar@{-}[rd] & & \\
\text{P.-A. Magma} \ar@{-}[dddd]  \ar@{-}[rd] & & \text{Alt. Magma} \ar@{-}[dddd]  \ar@{-}[ld] & & \text{Magma w. U.} \ar@{-}[ldddddd] \ar@{-}[ddddddddd] & \\
 & \text{Di-A. Magma} \ar@{-}[dddd] \ar@{-}[rrrrd] & & & & \\
 & & & & & \text{Ass. Magma} \ar@{-}[lddddddd] \\
 & & & \text{Quasigroup} \ar@{-}[ddd] \ar@{-}[llld]  \ar@{-}[ld]  & & \\
 \text{P.-A. Quasigr.} \ar@{-}[ddd]  \ar@{-}[rd] & & \text{Alt. Quasigr.} \ar@{-}[ddd]  \ar@{-}[ld] & & & \\
 & \text{Di-A. Quasigr.} \ar@{-}[ddd] & & & & \\
 & & & \text{Loop} \ar@{-}[llld]  \ar@{-}[ld]  & & \\
 \text{P.-A. Loop} \ar@{-}[rd] & & \text{Alt. Loop}  \ar@{-}[ld] & & & \\
 & \text{Di-A. Loop} \ar@{-}[rd] & & & & \\
 & & \text{Moufang Loop} \ar@{-}[rd] & & \text{Monoid} \ar@{-}[ld] & \\
 & & & \text{Group} & & 
}
\end{math}
\end{center}

\pagebreak
\subsubsection{A festival of pitfalls: powers (or multiples) of an element}\label{ssc:pitpow}
\indent

We will now take an interlude with an eye towards applications (however remote). Several
successful methods in cryptography involve taking powers (or multiples) of an element
belonging to rather rich structures, groups at the least, usually commutative
(\cite{Zemor(2000)}). One can't help but wonder what can be salvaged of these
methods in far poorer settings. Of course, the crux of these techniques usually is a variation
on Fermat's little theorem, and we will not investigate this here, but the simple idea of
taking powers bears some investigating.

Given any magma  \((S,f)\), for any \(a \in S\) and any non-zero integer
\(n\), we define \(a^n\) by induction, \emph{i.e.}, \(a^1 = 1\), and having
defined \(a^p\) for some non-zero integer \(p\), we define
\(a^{p+1} = f(a,a^p) = f(a^p,a)\).

While this is a perfectly reasonable definition, which might even be useful in some
contexts, without any additional constraints on the magma there not much that can be
proved about these objects. For instance, we might want to have the same values if we
repeatedly multiply on the right by \(a\) rather than multiply on the left as we have done.
Unfortunately, without any assumption of commutativity or associativy, in whatever form,
it may happen that this is not the case\footnote{Consider:
\begin{center}
\begin{tabular}{c|ccc}
f & a & b & c \\
\hhline{----}
a & b & c & a \\
b & a & a & a \\
c & a & a & a
\end{tabular}
\end{center}

In this case \(a^{2} = f(a,a) = b\) and
\(a^{3} = f(a,b) = f(a,f(a,a)) = c \neq a = f(b,a) = f(f(a,a),a)\).}.

Commutativity will avoid this complication, but we will in this document consider
very interesting cases which are not commutative. Furthermore, commutativity on its own
will not be sufficient to insure some nice properties such as:
\begin{equation}
(\forall\thinspace n \in \bbN-\{0\})
(\forall\thinspace m \in \bbN-\{0\})\thickspace
a^{n+m} = f(a^n,a^m) = f(a^m,a^n)\tag{A}\label{eqn:AdditiveNP}
\end{equation}

\noindent
which some forms of associativity will provide.

Alternativity is a rather poor choice for a shade of associativity, as then
(\ref{eqn:AdditiveNP}) is true if \(m\) is a power of \(2\) (as can be proved by
induction), but may otherwise fail\footnote{Consider:
\begin{center}
\begin{tabular}{c|ccccc}
f & a & b & c & d & e \\
\hhline{------}
a & b & c & d & e & e \\
b & c & d & e & e & e \\
c & d & e & d & e & e \\
d & e & e & e & e & e \\
e & e & e & e & e & e
\end{tabular}
\end{center}

In this case,
\(a^1 = a\), \(a^2 = b\), \(a^3 = c\), \(a^4 = d\), \(a^5 = e\), but
\(a^6 = e \neq d = f(a^3,a^3)\), so (\ref{eqn:AdditiveNP}) does not hold here for
\(m = 3\).}. A power-associative magma, on the other hand
does verify (\ref{eqn:AdditiveNP}) for every element.

If the magma under consideration is not only power-associative, but is also with unit,
with neutral element \(e\), we define \(a^0 = e\), in addition to the above, and now have
an improvement over  (\ref{eqn:AdditiveNP}):
\begin{equation}
(\forall\thinspace n \in \bbN)
(\forall\thinspace m \in \bbN)\thickspace
a^{n+m} = f(a^n,a^m) = f(a^m,a^n)\tag{A'}\label{eqn:AdditiveN}
\end{equation}

If, in addition, the element \(a\) happens to have an inverse \(\alpha\), we define for any
integer \(n\), \(a^{-n}\) as \(\alpha^n\). Unfortunately, if the magma is no more
than power-associative, it may happen\footnote{Consider:
\begin{center}
\begin{tabular}{c|cccc}
f & a & b & c & d \\
\hhline{-----}
a & a & b & c & d \\
b & b & c & b & a \\
c & c & b & c & a \\
d & d & a & a & a
\end{tabular}
\end{center}

In this case, the law is power-associative, has a unit, \(a\), but the element \(b\)
has an inverse, \(d\), while \(f(b^{-1},b^2) = f(d,c) = a \neq b\).}
that \(f(a^{-1},a^{2}) \neq a\). If, however, the magma is actually di-associative,
or is a power-associative loop (which is more constraining than being power-associative as a
magma), then we can improve upon (\ref{eqn:AdditiveN}):
\begin{equation}
(\forall\thinspace n \in \bbZ)
(\forall\thinspace m \in \bbZ)\thickspace
a^{n+m} = f(a^n,a^m) = f(a^m,a^n)\tag{A''}\label{eqn:AdditiveZ}
\end{equation}

We may want to investigate what can be said for products of powers of two elements which
may not be inverses of each other. More
precisely, even if the magma is without unit, we would like to have the following property:
\begin{multline}
(\forall\thinspace n \in \bbN-\{0\})
(\forall\thinspace m \in \bbN-\{0\})(\forall\thinspace p \in \bbN-\{0\})
(\forall\thinspace q \in \bbN-\{0\}) \\
f(a^{n+m},b^{p+q}) = f(a^m,f(a^n,b^{p+q})) = f(f(a^{n+m},b^p),b^q)\tag{B}\label{eqn:AdditiveDNP}
\end{multline}

If the magma is di-associative, then indeed, we do have (\ref{eqn:AdditiveDNP}).
Furthermore, if the magma also has a unit we can extend (\ref{eqn:AdditiveNP}) into:
\begin{equation}
(\forall\thinspace n \in \bbN)
(\forall\thinspace m \in \bbN)(\forall\thinspace p \in \bbN)
(\forall\thinspace q \in \bbN) \thickspace
f(a^{n+m},b^{p+q}) = f(a^m,f(a^n,b^{p+q})) = f(f(a^{n+m},b^p),b^q)\tag{B'}\label{eqn:AdditiveDN}
\end{equation}

About the best we can hope for, if the magma has a unit and \(a\) and \(b\) have inverses
would be:
\begin{equation}
(\forall\thinspace n \in \bbZ)
(\forall\thinspace m \in \bbZ)(\forall\thinspace p \in \bbZ)
(\forall\thinspace q \in \bbZ) \thickspace
f(a^{n+m},b^{p+q}) = f(a^m,f(a^n,b^{p+q})) = f(f(a^{n+m},b^p),b^q)\tag{B''}\label{eqn:AdditiveDZ}
\end{equation}

As one might come to expect, if the only shade of associativity which can be mustered is
di-associativity, this might fail\footnote{Consider the set \(S = \{2^{p} | p \in \bbZ\}
\cup \{3^{q} | q \in \bbZ\} \cup \{5^{r} | r \in \bbZ\}\) along with the following law
``\( \star \)'':
\(2^{p} \star 2^{p'} = 2^{p+p'}\), \(3^{q} \star 3^{q'} = 3^{q+q'}\),
\(5^{r} \star 5^{r'} = 5^{r+r'}\),
\(2^{p} \star 3^{q} = 3^{q} \star 2^{p} = 5^{p+q}\),
\(2^{p} \star 5^{r} = 5^{r} \star 2^{p} = 5^{p+r}\),
\(3^{q} \star 5^{r} = 5^{r} \star 3^{q} = 5^{q+r}\).

One verifies that \((S, \star)\) is indeed di-associative, has a neutral element (\(1\)),
and that every element has inverses: an element of the form \(2^p\) has \(2^{-p}\) for
only inverse, an element of the form \(3^q\) has \(3^{-q}\) for only inverse, and an
element of the form \(5^r\) have exactly \emph{three} inverses (\(2^{-r}\), \(3^{-r}\)
and \(5^{-r}\)) if \(r \neq 0\) (if \(r = 0\), of course it only has \(1\) for inverse).

We now see that \(2^{-1} \star (2^{+1} \star 3^{+1}) = 2^{-1} \star 5^{2} = 5^{+1}
\neq 3^{+1} = (2^{-1} \star 2^{+1}) \star 3^{+1}\).}. There are, however, very
interesting cases where (\ref{eqn:AdditiveDZ}) is true without imposing (full) associativity
(which is of course enough to guaranty (\ref{eqn:AdditiveDZ})).

\begin{proposition}
If \((S,f)\) is a di-associative magma with unit, and if \(S\) is \emph{finite}, then
(\ref{eqn:AdditiveDZ}) holds for every \(a\) and \(b\) which have inverses.
\end{proposition}

\begin{proof}
Let \(e\) be the neutral element of \((S,f)\), and let \(a \in S\) having \(a^{-1}\) as
inverse. Since we have assumed \((S,f)\) to be di-associative, we immediately see that the
sub-magma generated by \(a\) and \(a^{-1}\) is actually a group, and since it is a subset
of \(S\) which is finite, it is a finite group. Hence
\((\exists\thinspace n_{0} \in \bbN)\thickspace a^{-1} = a^{n_{0}}\), and
(\ref{eqn:AdditiveDZ}) reduces to (\ref{eqn:AdditiveDN}).
\end{proof}

It is also worth mentioning, irrespective of whether \(S\) is finite or not,
that (\ref{eqn:AdditiveDZ}) holds if \((S,f)\) is a
di-associative loop (which is a more stringent condition than being di-associative as a
magma), in particular a Moufang loop.

\subsection{Structures involving more than one law}\label{sbs:manylaws}
\subsubsection{Familiar structures}
\indent

A \emph{ring} is a structure \((A,+,\times)\) such that \((A,+)\) is a commutative group
(whose neutral element we will write \(0\)),
and that the multiplication (the law ``\(\times\)'') is distributive over the
addition\footnote{Which means that \((\forall\thinspace x \in A)
(\forall\thinspace y \in A)(\forall\thinspace z \in A)\thickspace
[(x+y)\times z = (x\times z) + (y\times z) \text{ and }
z\times(x+y) = (z\times x) + (z\times y)]\)}.
A ring is said to be \emph{with unit} (respectively \emph{commutative}, 
\emph{monogenic}, \emph{power-associative}, \emph{alternative}, \emph{associative})
if and only if the multiplication is with unit (respectively commutative, monogenic,
power-associative, alternative, associative). As we will see later on, if the multiplication is
alternative, it will turn out to be di-associative as well, contrary to what happens for magmas
(in other words, not all magma laws are fit to play the role of multiplication in rings). Of
course, di-associative multiplications are still alternative!

Given two rings \((A,+,\times)\) and \((B,\dagger,*)\), a \emph{ring (homo)marphism}
is a function \(\theta: A \rightarrow B\) such that it is a group morphism from 
\((A,+)\) to \((B,\dagger)\), and a magma morphism from \((A,\times)\) to
\((B,*)\). If \(A\) and \(B\) are with unit, whose neutral elements are respectively
\(e\) and \(\varepsilon\), then a \emph{morphism of ring with unit} is a morphism
of rings such that \(\theta(e) = \varepsilon\). We define, epi-, mono-, endo-, iso- and
automorphisms for rings as we have for the other structures, from ring homomorphisms.

\pagebreak
Note that if \(A\) is a ring with unit, whose neutral element we will write \(1_{A}\),
then \(1_{A} \neq 0 \Leftrightarrow A \neq \{0\}\).
Given a ring with unit \(A\) with at least two elements, we will denote by
\(A^*\) the set of elements of \(A\) which have an inverse (for the multiplication
in \(A\)); of course \(A^* \subset A-\{0\}\)
but the inclusion may be strict. The elements of \(A^*\) are also
(very unfortunately) called the \emph{units} of \(A\).
The magma with unit \((A^{*},\times)\), has a rather poor structure in general, but
when \(A\) is associative, it is a group  (\cite{Lang(1971)}); when \(A\) is only
alternative, it is a Moufang loop (\cite{GoodaireJespersMiles(1996)}; we will show a special
case of this result in this document).

Given a ring \((A,+,\times)\) and a \(T \subset A\), we will say that
\(T\) is a  \emph{subring} if and only if \((T,+)\) is a subgroup of \((A,+)\), and \(T\)
is stable for the multiplication of \(A\) (\emph{i.e.}
\((\forall (x,y) \in T^{2})\thickspace x \times y \in T\)). If \((A,+,\times)\) is a
ring \emph{with unit}, whose neutral element for ``\(\times\)'' we will denote by \(e\),
then \(T\) will be said to be a \emph{subring of a ring with unit } if and only if it is a
subring and furthermore that \(e \in T\)\footnote{It is well-known that
undesirable things may happen otherwise, and that
\((T,+|_{T \vartimes T}^{T},\times|_{T \vartimes T}^{T})\) may happen to be a ring
with unit, without being a ``subring of a ring with unit'' as per our definition; consider the
square matrices of order two on \(\bbR\) for \(A\), and the matrices of the form
\(\left(
\begin{smallmatrix}
x & 0\\
0 & 0
\end{smallmatrix}
\right)\) with \(x \in \bbR\) for T.}.

A \emph{semifield} is a structure \((K,+,\times)\) which is a ring with unit such that,
writing \(K^* = K-\{0\}\) (with \(0\) being the neutral element for the addition,
\emph{i.e.} the law \(+\)), \((K^*,\times)\) is a loop, whose neutral element which
we will call \(1_{K}\) verifies \(0 \neq 1_{K}\).
A semifield is \emph{commutative}
(respectively \emph{monogenic}, \emph{power-associative}, \emph{alternative},
\emph{associative}) if
it is such as a ring.

A \emph{field} is\footnote{\cite{Lang(1971)} requires the multiplication to be
commutative; we do not, however, in accordance with \cite{Bourbaki(A)}.}
an associative semifield.
Denoting by \(F^*\) the set of elements of \(F\) which have an inverse
(for the multiplication in \(F\)) we therefore have \(F^* = F-\{0\}\); \((F^*,\times)\)
is a group. A field is said to be \emph{commutative} (respectively \emph{monogenic})
if, as a semifield, it is such.

Morphisms for semifields and fields are those of the underlying rings.

Given an associative and commutative ring with unit \(A\), whose neutral element we will
denote by \(1_{A}\), a \emph{module} over \(A\) (also called an
\(A\)\emph{-module}) is a structure \((M,+,\cdot)\) where \((M,+)\) is a
commutative group, and such that the function (called an \emph{external law}) ``\(\cdot\)''
is defined on all of \(A \vartimes M\) with values in \(M\), is distributive both on
the addition in \(A\) and the addition in \(M\), and verify
\((\forall x \in M)\thickspace 1_{A}\thinspace \cdot \thinspace x = x\) and
\((\forall x \in M)(\forall \alpha \in A)(\forall \beta \in A)\thickspace
\alpha \cdot (\beta \cdot x) = (\alpha\thinspace . \thinspace\beta) \cdot x\).

Given an associative ring with unit \((A,+,.)\), an \(A\)-module \((M,+,\cdot)\) and a
\(T \subset M\), we will say \(T\) is a \emph{submodule} if and only if \((T,+)\) is a
subgroup of \((M,+)\), and \(T\) is stable for the exterior law (\emph{i.e.}
\((\forall\thinspace \alpha \in A)(\forall\thinspace x \in T)\thickspace
\alpha \cdot x \in T\)).

A \emph{vector space} is a module over a commutative field. A \emph{sub-vector space}
is a submodule over a commutative field.

Given two modules  \((M,+,.)\) and \((N,\dagger,\bullet)\) over the \emph{same}
commutative ring \(A\), a \emph{module (homo)morphism} is a function
\(\theta: M \rightarrow N\) such that it is a group homomorphism from \((M,+)\)
to \((N,\dagger)\), and such that \((\forall\thinspace \alpha \in A)
(\forall\thinspace m \in M)\thickspace \theta(\alpha . m) =
\alpha \bullet \theta(m)\). Morphisms for vector fields are those of the underlying
module. The usual morphism declinations apply here as well.

Given a commutative and associative ring with unit \((A,+,.)\), an \emph{algebra} over
\(A\) (also called an \(A\)\emph{-algebra}) is a structure
\((E,+,\times,\cdot)\) where \((E,+,\cdot)\) is an \(A\)-module\footnote{We do
not limit ourselves to the case where \(A\) is a field, which is the most common situation
(\cite[\ldots]{Schafer(1966)}), but which would be too restrictive for our needs in this
document.}, and the law
 ``\(\times\)'' (the multiplication) \(E\) is \(A\)-bilinear, \emph{i.e.} it is distributive
 over ``+'' and \((\forall\thinspace \alpha \in A)(\forall\thinspace \beta \in A)
(\forall\thinspace x \in E)(\forall\thinspace y \in E)\thickspace
(\alpha \cdot x) \times (\beta \cdot y) = (\alpha . \beta) \cdot (x \times y)\).
An \emph{algebra with unit} is an algebra whose multiplication has a neutral element.
In particular, an algebra is a ring, and an algebra with unit is a ring with unit.
An algebra is said to be \emph{commutative}
(respectively \emph{power-associative}, \emph{alternative}, \emph{associative})
if and only if the multiplication is
commutative (respectively power-associative, alternative, associative).

Given an associative ring with unit \((A,+,.)\), an \(A\)-algebra \((E,+,\times,\cdot)\)
and a \(T \subset E\), we will say \(T\) is a \emph{sub-algebra} if and only if 
\((T,+,\cdot)\) is a submodule of \((E,+,\cdot)\), and \(T\) is stable for
``\(\times\)'' (\emph{i.e.} \((\forall\thinspace (x,y) \in T^{2})\thickspace
(x \times y) \in T\)). If \((E,+,\times,\cdot)\) is an \(A\)-algebra with unit, with
neutral element \(e\), \(T\) will be a \emph{sub-algebra of an algebra with unit} is and
only if is is a sub-algebra of E and \(e \in T\).

If \(E\) is an \(A\)-algebra with unit, we will also denote by \(E^*\) the set of elements
of \(E\) which have an inverse (for the multiplication in \(E\)); of course
\(E^* \subset E-\{0\}\) but the inclusion may be strict. As with rings in general, the
magma with unit \((E^*,\times)\) usually has a rather poor structure.

\pagebreak
If \(E\) is an \(A\)-algebra with unit, whose neutral element we will denote by \(e\),
there is a natural epimorphism of rings with unit defined by
\([A \rightarrow E; a \mapsto a \cdot e]\). It should be noted that this homomorphism
need not be injective\footnote{The following example was communicated to me by
M.~Daniel~R.~Grayson:  \(A = \bbZ\) and \(E = \bbZ/2\bbZ\); in this case \(e\) is
\emph{not} linearly independent on \(A\).}!

Given two algebras \((E,+,\times,\cdot)\) and \((F,\dagger,*,\bullet)\) over the
\emph{same} ring \(A\), an \emph{algebra (homo)morphism} is a function
\(\theta : E \rightarrow F\) such that it is a morphism of modules from \((E,+,\cdot)\)
to \((F,\dagger,\bullet)\) and a morphism of ring from \((E,+,\times)\) to
\((F,\dagger,*)\). If \(E\) and \(F\) are with unit, a \emph{morphism of algebra with
unit} is a morphism of algebra such that it is also a morphism of ring with unit from
\((E,+,\times)\) to \((F,\dagger,*)\). The usual offshoots of the morphism family are
dealt with as usual.

\subsubsection{Alternative rings}\label{ssc:altrng}
\indent

We recall here some important properties of alternative rings that we need in this document.
In particular, we present here a rather elementary proof of a classical, which is a special case
of a theorem of Bruck and Kleinfeld (\cite{GoodaireJespersMiles(1996)}), itself an
extension of a result due to Artin, classically known for algebras
(\cite{Bourbaki(A), Schafer(1966)}). This is mostly for self-sufficiency's sake, but we
believe it fits quite well with this appendix' emphasis on shades of associativity; it also is an
opportunity to remedy what is most likely an unfortunate ellipsis in the presentation found in
\cite{Bourbaki(A)}.

Given a ring \((A,+,\times)\) and a function \(f: A^{n} \rightarrow A\), we will say
that \(f\) is \emph{additive} with respect to the \(i^{\text{th}}\) variable if and only if:
\begin{multline}
(\forall\thinspace (x_{1},\ldots,x_{i-1},x_{i+1},\ldots,x_{n}) \in A^{n-1})
(\forall\thinspace x \in A)(\forall\thinspace x' \in A)\thickspace
f(x_{1},\ldots,x_{i-1},x+x',x_{i+1},\ldots,x_{n}) = \\
f(x_{1},\ldots,x_{i-1},x,x_{i+1},\ldots,x_{n})+
f(x_{1},\ldots,x_{i-1},x',x_{i+1},\ldots,x_{n}) \notag
\end{multline}

We trivially verify that a function \(f\) additive with respect to the \(i^{\text{th}}\)
variable takes the value \(0\) when the \(i^{\text{th}}\) variable takes the value \(0\),
and that furthermore
\begin{multline}
(\forall\thinspace (x_{1},\ldots,x_{i-1},x_{i+1},\ldots,x_{n}) \in A^{n-1})
(\forall\thinspace x \in A)\thickspace \\
f(x_{1},\ldots,x_{i-1},-x,x_{i+1},\ldots,x_{n}) =
-f(x_{1},\ldots,x_{i-1},x,x_{i+1},\ldots,x_{n})\notag
\end{multline}

We will say that a function \(f\), additive with respect to each of its variables, is
\emph{alternated} if and only if:
\begin{equation}
(\forall\thinspace i \in \{1,\ldots,n-1\})
(\forall\thinspace (x_{1},\ldots,x_{n}) \in A^{n})\thickspace
[x_{i} = x_{i+1} \Rightarrow f(x_{1},\ldots,x_{n}) = 0]\notag
\end{equation}

Let us denote by \(\mathfrak{S}_{n}\) the group of permutations on \(\{1,\ldots,n\}\)
and, for \(\sigma \in \mathfrak{S}_{n}\), by \(\epsilon_{\sigma}\) the signature of
\(\sigma\). We trivially prove that if \(f\) is alternated then
\(f(x_{\sigma(1)},\ldots,x_{\sigma(n)}) =
\epsilon_{\sigma}\thinspace f(x_{1},\ldots,x_{n})\), and that a function is
alternated if and only if it is zero whenever any two variables have the same value, whether
these variables are consecutive or not.

Given any ring \((A,+,\times)\), we can always define the \emph{alternator} by:
\begin{equation*}
\mathfrak{a}\negthinspace : \negthinspace
\begin{array}[t]{rcl}
A^{3} & \rightarrow & A\\
(x,y,z) & \mapsto & x \times (y \times z) - (x \times y) \times z
\end{array}
\end{equation*}
\noindent
and it is always true that the alternator is additive with respect to each variable. It is also
always true that:
\begin{equation}
(\forall\thinspace (p,q,r,s) \in A^{4})\thickspace
\mathfrak{a}(p \times q,r,s) - \mathfrak{a}(p,q \times r,s) +
\mathfrak{a}(p,q,r \times s) =
p \times \mathfrak{a}(q,r,s) + \mathfrak{a}(p,q,r) \times s
\tag{T}\label{eqn:foursome}
\end{equation}
\noindent
(which is proved by simply evaluating both sides, and is known as the Teichm{\"u}ller
identity).

Another interesting function is the Kleinfeld function (which differ only in sign from the
definition in \cite{ZhevlakovSlin'koShestakovShirshov(1982), GoodaireJespersMiles(1996)}):
\begin{equation*}
\mathfrak{k}(p,q,r,s) =
q\times\mathfrak{a}(p,r,s)+\mathfrak{a}(q,r,s)\times p-\mathfrak{a}(p\times q,r,s)
\end{equation*}
which also is an additive function with respect to each variable.

The alternator
is interesting as we trivially prove that \(A\) is alternative if and only if \(\mathfrak{a}\)
is alternated, and that \(A\) is associative if and only if  \(\mathfrak{a}\) is the zero
function. Likewise, a subset \(T\) of \(A\) is strongly associative (regardless of
whether \(A\) is alternative or not) if and only if
\begin{equation}
(\forall\thinspace (u,v,w) \in A^{3})\thickspace
[[[(u,v) \in T^{2}] \text{ or } [(u,w) \in T^{2}] \text{ or } [(v,w) \in T^{2}]]
\Rightarrow \mathfrak{a}(u,v,w) = 0] \notag
\end{equation}
and if \(A\) is indeed alternative, then \(T\) is strongly associative if and only if
\begin{equation}
(\forall\thinspace (u,v) \in T^{2})(\forall\thinspace w \in A)\thickspace
\mathfrak{a}(u,v,w) = 0 \notag
\end{equation}

As well we have (\cite{GoodaireJespersMiles(1996)}):

\begin{proposition}\label{prp:kleinalt}
If \((A,+,\times)\) is \emph{alternative} then \(\mathfrak{k}\) is alternated.
\end{proposition}

\begin{proof}
We first compute \(\mathfrak{k}(s,p,q,r) = p\times\mathfrak{a}(s,q,r)+
\mathfrak{a}(p,q,r)\times s-\mathfrak{a}(s\times p,q,r)\).

Using (\ref{eqn:foursome}) we find that \(\mathfrak{k}(s,p,q,r) =
\mathfrak{a}(p\times q,r,s)-\mathfrak{a}(p,q\times r,s)
+\mathfrak{a}(p,q,r\times s)-\mathfrak{a}(s\times p,q,r)\).

However, writing (\ref{eqn:foursome}) for \((q,r,s,p)\) we see that
\(\mathfrak{a}(q \times r,s,p) - \mathfrak{a}(q,r \times s,p) +
\mathfrak{a}(q,r,s \times p) =
q \times \mathfrak{a}(r,s,p) + \mathfrak{a}(q,r,s) \times p\) and since \(A\) is
alternated, \(\mathfrak{a}(r,s,p) = \mathfrak{a}(p,r,s)\), so
\(\mathfrak{k}(p,q,r,s) =
\mathfrak{a}(q\times r,s,p)-\mathfrak{a}(p\times q,r,s)
+\mathfrak{a}(q,r,s\times p)-\mathfrak{a}(q,r\times s,p)\).

Since \(A\) is alternative, the associator is alternated, so
\(\mathfrak{a}(p,q\times r,s) = \mathfrak{a}(q\times r,s,p)\),
\(\mathfrak{a}(p,q,r\times s) = \mathfrak{a}(q,r\times s,p)\) and
\(\mathfrak{a}(s\times p,q,r) = \mathfrak{a}(q,r,s\times p)\) and thus
\(\mathfrak{k}(s,p,q,r) = -\mathfrak{k}(p,q,r,s)\).

Since the alternator is alternated because the ring is alternative, we have
\(\mathfrak{k}(p,q,r+s,r+s) = \mathfrak{k}(p,q,r,r) = \mathfrak{k}(p,q,s,s) = 0\)
and so \(\mathfrak{k}(p,q,s,r) = -\mathfrak{k}(p,q,r,s)\).

Finally,
\(
\left(
\begin{smallmatrix}
1&2&3&4\\
4&1&2&3
\end{smallmatrix}
\right)
\) and
\(
\left(
\begin{smallmatrix}
1&2&3&4\\
1&2&4&3
\end{smallmatrix}
\right)
\) generate the symmetric group of order \(4\).
\end{proof}

The alternator is fundamental for proving the next theorem, which contrasts
strongly with the case of magmas. 

\begin{theorem}[Di-associativity of alternative rings]\label{thm:diassaltrng}
A ring is alternative if and only if it is di-associative, and if and only if all subrings generated
by any two (not necessarily distinct) elements are associative. A ring with unit is alternative
if and only if it is di-associative, and if and only if all subrings of a ring with unit generated
by any two (not necessarily distinct) elements are associative.
\end{theorem}

We will mostly, but not quite completely, paraphrase the proof in \cite{Bourbaki(A)},
which requires a few lemmas.

We start with a technical result which shows that the presence of an additional law on which
the product is distributive imposes some regularity on the product.

\begin{lemma}\label{lem:strassgrp}
Let \((A,+,\times)\) be any ring, and let \(\rsfsE_{\times}\) the set of all subsets of
\(A\) which are strongly associative for ``\(\times\)''. Let \(G\) be a maximal
element (for the inclusion) of \(\rsfsE_{\times}\); then if \(G \neq \emptyset\),
\((G,+)\) is a group.
\end{lemma}

Note that we know there are maximal elements in \(\rsfsE_{\times}\) thanks to
\thref{prop:indstrass}, and the Zorn theorem.

\begin{proof}
Assume that \(G \neq \emptyset\).

Let's first prove that \(0 \in G\). Consider \(K = G \cup \{0\}\). Of course
\(G \subset K\), so as \(G\) is maximal, it is enough to prove that \(K\) is
strongly associative. Let \((u,v,w) \in A^{3}\); if two of the three are in \(G\), then
\(\mathfrak{a}(u,v,w) = 0\), because \(G\) is strongly associative, and if at least one is
zero, then \(\mathfrak{a}(u,v,w) = 0\) as well because the alternator is additive with
respect to each of its variables. Hence \(K\) is strongly additive.

Let \(u \in G\); considering \(L = G \cup \{-u\}\), and invoking the additivity of the
alternator and the maximality of \(G\), we see as well that \(-u \in G\).

Let \(u \in G\), \(v \in G\); considering \(M = G \cup \{u+v\}\), and invoking the
additivity of the alternator and the maximality of \(G\), we also see that \((u+v) \in G\).
\end{proof}

\begin{lemma}\label{lem:genrngstrass}
Let \((A,+,\times)\) be an alternative ring, \(H\) a subset of \(A\) strongly
associative for ``\(\times\)'' and \(F\) the subring of \(A\) generated by \(H\);
then \(F\) is strongly associative.
\end{lemma}

\begin{proof}
If \(F = \emptyset\), there is nothing to prove. We therefore assume that
\(F \neq \emptyset\).

As the subsets of \(F\) which are strongly associative for ``\(\times\)'' are inductive for
set inclusion, by \thref{prop:indstrass}, we may consider \(G\) the greatest (for set
inclusion) strongly associative subset of \(F\) which contains \(H\). to prove the lemma
it therefore suffices to prove that \(G = F\).

Using \thref{lem:strassgrp}, we know that \((G,+)\) is a group, so all that is left to prove
is that it is also stable for ``\(\times\)'' (as then \(G\) will be a subring of \(A\)
containing \(H\), \(F\) is by definition the smallest subring of \(A\) containing \(H\),
and \(G \subset F\)).

Let \(u \in G\), \(v \in G\), and consider \(L = G \cup \{u \times v\}\).
Of course \(H \subset L\) (since \(H \subset G\)), and \(L \subset F\) (as \(F\)
is a ring). Consider \(x \in L\), \(y \in L\), \(z \in A\) and
\(\alpha = \mathfrak{a}(x,y,z)\).
If \(x \in G\) and \(y \in G\) then \(\alpha = 0\) because \(G\) is strongly associative.
Assume therefore that \(x = u \times v\). If we also have \(y = u \times v\) then
\(\alpha = 0\) because the alternator is alternated. Assume therefore only that \(y \in G\).
Using (\ref{eqn:foursome}) we find that
\begin{equation}
-\alpha = \mathfrak{a}(x,z,y) = \mathfrak{a}(u \times v,z,y) =
\mathfrak{a}(u,v \times z,y) - \mathfrak{a}(u,v,z \times y) +
\mathfrak{a}(u,v,z) \times y + u \times \mathfrak{a}(v,z,y)\notag
\end{equation}
However, since \(u \in G\), \(v \in G\) and \(y \in G\) and \(G\) is strongly
associative, \(\mathfrak{a}(v,z,y) = 0\), \(\mathfrak{a}(u,v,z) = 0\),
\(\mathfrak{a}(u,v,z \times y) = \mathfrak{a}(u,v,z') = 0\) and
\(\mathfrak{a}(u,v \times z,y) = \mathfrak{a}(u,z'',y) = 0\). Hence \(\alpha = 0\).
The case where we assume \(x \in G\) and \(y = u \times v\) is deduced from the above
by \(\mathfrak{a}(x,y,z) = -\mathfrak{a}(y,x,z)\).
\end{proof}

\begin{lemma}\label{lem:monostrass}
Let \((A,+,\times)\) be an alternative ring; for any \(x \in A\), the subring generated by
\(x\) is strongly associative.
\end{lemma}

\begin{proof}
Let \(H = \{x\}\); \(H\) is strongly associative since the alternator is alternated. We
conclude by using \thref{lem:genrngstrass}.
\end{proof}

\begin{lemma}\label{lem:stblsubrng}
Let \((A,+,\times)\) be an alternative ring, \(T\) a subgroup of \((A,+)\), \(S\) a
subset of \(T\), and \(U\) the subring generated by \(S\). If
\((\forall\thinspace x \in S)
[\thinspace x \times T \subset T \text{ and } T \times x \subset T \thinspace]\)
then \(T \supset U\).
\end{lemma}

\begin{proof}
Let \(V = \{ x \in U \thinspace | \thinspace x \times T \subset T \text{ and }
T \times x \subset T \}\). By hypothesis \(S \subset V\). Let \(x \in V\),
\(y \in V\) and \(t \in T\).

Let \(P = (x \times y) \times t\). Obviously, \(P = x \times (y \times t) -
[\thinspace x \times (y \times t) - (x \times y) \times t \thinspace]\), which by
definition of the alternator means that \(P = x \times (y \times t) - \mathfrak{a}(x,y,t) \),
or, as the alternator is alternated, \(P = x \times (y \times t) + \mathfrak{a}(x,t,y) \),
which evaluates into \(P = x \times (y \times t) + x \times (t \times y) -
(x \times t) \times y\).

Since \(y \in V\) and \(t \in T\), we have both \(y \times t \in T\) and
\(t \times y \in T\). But \(y \times t \in T\) and \(x \in V\) together imply
\(x \times (y \times t) \in T\), and likewise \(x \times (t \times y)\),
\((x \times t) \times y \in T\). Since \((T,+)\) is a group, this proves that
\((x \times y) \times t = P \in T\).

Considering \(Q = t \times (x \times y)\), we likewise prove that
\( t \times (x \times y) \in T\).

Hence \(V\) is a subring of \(A\), and since \(S \subset V\) and \(U\) is the subring
generated by \(S\), then \(V \supset U\). As \(V \supset U\) by definition, this means
that \(V = U\), and therefore that \((\forall\thinspace x \in U)
[\thinspace x \times T \subset T \text{ and } T \times x \subset T \thinspace]\).

As \(T \subset U\), this means that \emph{a fortiori}, \(T \times T \subset T\), and
so that \(T\) is a subring of \(A\). As \(S  \subset T\) by definition, and \(U\) is the
subring generated by \(S\), \(T \supset U\).
\end{proof}

\begin{remark}
If in \thref{lem:stblsubrng} \(T \subset U\) then of course \(T = U\)...
\end{remark}

\begin{lemma}\label{lem:unionstrass}
Let \((A,+,\times)\) be an alternative ring and \(X\) and \(Y\) two strongly associative
\emph{subrings} of \(A\). The subring generated by \(X \cup Y\) is associative.
\end{lemma}

\begin{proof}
Let \(W\) be the subring generated by \(X \cup Y\).

Consider \(Z = \{z \in W \thinspace |
(\forall\thinspace x \in X)(\forall\thinspace y \in Y)\thickspace
\mathfrak{a}(x,y,z) = 0\}\).

We trivially see that \((Z,+)\) is a group. Since \(X\) is
strongly associative, we have \(X \subset Z\) and likewise \(Y \subset Z\). By
definition we have \(Z \subset W\). So \(W\) is also the subring generated by \(Z\).

Let \(x \in X\), \(x' \in X\), \(y \in Y\), \(z \in Z\).

Using (\ref{eqn:foursome}) we find that
\begin{equation}
\mathfrak{a}(x' \times x,z,y)-\mathfrak{a}(x',x \times z,y)+
\mathfrak{a}(x',x,z \times y) = \mathfrak{a}(x',x,z) \times y + 
x' \times \mathfrak{a}(x,z,y)\notag
\end{equation}
but \(\mathfrak{a}(x',x,z) = 0\) and \(\mathfrak{a}(x',x,z \times y) = 0\) because
\(x \in X\) and \(x' \in X\) and \(X\) is strongly associative,
\(\mathfrak{a}(x,z,y) = -\mathfrak{a}(x,y,z) = 0\) because \(x \in X\), \(y \in Y\)
and by definition of \(Z\). As \(X\) is a subring of \(A\), \(x'' = x' \times x \in X\),
so \(\mathfrak{a}(x' \times x, z,y) = -\mathfrak{a}(x'',y,z) = 0\).
Hence \(\mathfrak{a}(x',y,x \times z) = 0\).
This proves that \((\forall\thinspace x \in X)(\forall\thinspace z \in Z)\thickspace
x \times z \in Z\).

Using (\ref{eqn:foursome}) with \((p,q,r,s) = (y,z,x,x')\) we likewise prove that
\((\forall\thinspace x \in X)(\forall\thinspace z \in Z)\thickspace
z \times x \in Z\).

Since the alternator is alternated, \(Z\) can also be defined as
\(\{z \in W \thinspace |
(\forall\thinspace x \in X)(\forall\thinspace y \in Y)\thickspace
\mathfrak{a}(y,x,z) = 0\}\) and what we just proved also implies that
\((\forall\thinspace y \in Y)(\forall\thinspace z \in Z)\thickspace
[y \times z \in Z \text{ and } z \times y \in Z]\).

Using \thref{lem:stblsubrng} we conclude that \(Z = W\); also
\((\forall\thinspace x \in X)(\forall\thinspace y \in Y)
(\forall\thinspace w \in W)\thickspace
\mathfrak{a}(x,y,w) = 0\).

Consider now \(Z' = \{z \in W \thinspace |
(\forall\thinspace x \in X)(\forall\thinspace y \in Y)
(\forall\thinspace w \in W)\thickspace
\mathfrak{a}(z,y,w) = \mathfrak{a}(x,z,w) = 0\}\).

Again, \((Z',+)\) is trivially a group. We have just seen that
\((\forall\thinspace y \in Y)(\forall\thinspace w \in W)
[z \in X \Rightarrow \mathfrak{a}(z,y,w) = 0]\), and as \(X\) is strongly associative,
\((\forall\thinspace x \in X)(\forall\thinspace w \in W)
[z \in X \Rightarrow \mathfrak{a}(x,z,w) = 0]\), so \(X \subset Z'\).
Likewise \(Y \subset Z'\). By definition, \(Z' \subset W\). So \(W\) is also
the subring generated by \(Z'\).

Let \(x \in X\), \(y \in Y\), \(y' \in Y\), \(w \in W\), \(z \in Z'\).

Using (\ref{eqn:foursome}) we find that
\begin{equation}
\mathfrak{a}(z \times y,y',w) - \mathfrak{a}(z,y \times y',w) +
\mathfrak{a}(z,y,y' \times w) =
\mathfrak{a}(z,y,y') \times w + z \times \mathfrak{a}(y,y',w)\notag
\end{equation}
but \(\mathfrak{a}(z,y,y') = 0\) and \(\mathfrak{a}(y,y',w) = 0\) because \(Y\) is
strongly associative, \(\mathfrak{a}(z,y,y' \times w) = 0\) because \(W\) is a subring
and by definition of \(Z'\), and \(\mathfrak{a}(z,y \times y',w) = 0\) because
\(Y\) is a subring and by definition of \(Z'\). Hence
\(\mathfrak{a}(z \times y,y',w) = 0\).

Using (\ref{eqn:foursome}) once more, we find that
\begin{equation}
\mathfrak{a}(z \times y,w,x) - \mathfrak{a}(z,y \times w,x) +
\mathfrak{a}(z,y,w \times x) =
\mathfrak{a}(z,y,w) \times x + z \times \mathfrak{a}(y,w,x)\notag
\end{equation}
but \( \mathfrak{a}(y,w,x) = 0\) because \(Z = W\), \(\mathfrak{a}(z,y,w) = 0\)
by definition of \(Z'\), \(\mathfrak{a}(z,y,w \times x) = 0\) and
\(\mathfrak{a}(z,y \times w,x) = 0\) because \(W\) is a subring and by definition of
\(Z'\). Hence \(\mathfrak{a}(z \times y,w,x) = 0\).

These two results prove that \((\forall\thinspace y \in Y)
(\forall\thinspace z \in Z')\thickspace z \times y \in Z'\).

Using (\ref{eqn:foursome}) again, we find that
\begin{equation}
\mathfrak{a}(y \times z,w,y') - \mathfrak{a}(y,z \times w,y') +
\mathfrak{a}(y,z,w \times y') =
\mathfrak{a}(y,z,w) \times y' + y \times \mathfrak{a}(z,w,y')\notag
\end{equation}
but \(\mathfrak{a}(z,w,y') = 0\) and \(\mathfrak{a}(y,z,w) = 0\) by definition of
\(Z'\), \(\mathfrak{a}(y,z,w \times y') = 0\) because \(W\) is a subring and by
definition of \(Z'\), and \(\mathfrak{a}(y,z \times w,y') = 0\) because \(W\) is a
subring and \(Y\) is strongly associative. Hence \(\mathfrak{a}(y \times z,w,y') = 0\).

Using (\ref{eqn:foursome}) yet again, we find that
\begin{equation}
\mathfrak{a}(y \times z,w,x) - \mathfrak{a}(y,z \times w,x) +
\mathfrak{a}(y,z,w \times x) =
\mathfrak{a}(y,z,w) \times x + y \times \mathfrak{a}(z,w,x)\notag
\end{equation}
but \(\mathfrak{a}(z,w,x) = 0\) and \(\mathfrak{a}(y,z,w) = 0\) by definition of
\(Z'\), \(\mathfrak{a}(y,z \times w,x) = 0\) because \(W\) is a subring and by
definition of \(Z'\), and \(\mathfrak{a}(y,z \times w,x) = 0\) because \(W\) is a
subring and \(Z = W\). Hence \(\mathfrak{a}(y \times z,w,x) = 0\).

These two results prove that \((\forall\thinspace y \in Y)
(\forall\thinspace z \in Z')\thickspace y \times z \in Z'\).

As \(X\) and \(Y\) play symmetrical roles in the definition of \(Z'\), we have therefore
also proved that
\((\forall\thinspace x \in X)(\forall\thinspace z \in Z')\thickspace
[x \times z \in Z' \text{ and } z \times x \in Z']\).

Using \thref{lem:stblsubrng} we conclude that \(Z' = W\); also
\((\forall\thinspace x \in X)(\forall\thinspace y \in Y)
(\forall\thinspace v \in W)(\forall\thinspace w \in W)\thickspace
\mathfrak{a}(x,v,w) = \mathfrak{a}(y,v,w) = 0\).

Consider finally \(Z'' = \{z \in W \thinspace |
(\forall\thinspace v \in W)(\forall\thinspace w \in W)\thickspace
\mathfrak{a}(z,v,w) = 0\}\).

Once again, \((Z'',+)\) is trivially a group. We have also just proved that \(X \subset Z''\)
and \(T \subset Z''\). By definition, \(Z'' \subset W\). So \(W\) is also
the subring generated by \(Z''\).

Let \(x \in X\), \(y \in Y\), \(v \in W\), \(w \in W\), \(z \in Z''\).

Using (\ref{eqn:foursome}) as usual, we find that
\begin{equation}
\mathfrak{a}(x \times z,v,w) - \mathfrak{a}(x,z \times v,w) +
\mathfrak{a}(x,z,v \times w) =
\mathfrak{a}(x,z,v) \times w + x \times \mathfrak{a}(z,v,w)\notag
\end{equation}
but \(\mathfrak{a}(z,v,w) = 0\) and \(\mathfrak{a}(x,z,v) = 0\) by definition of
\(Z''\), and \(\mathfrak{a}(x,z,v \times w) = 0\) and
\(\mathfrak{a}(x,z \times v,w)\) because \(W\) is a subring and because \(Z' = W\).
Hence \(\mathfrak{a}(x \times z,v,w) = 0\). Likewise
\(\mathfrak{a}(y \times z,v,w) = 0\).

Using (\ref{eqn:foursome}) one last time, we find that
\begin{equation}
\mathfrak{a}(z \times x,v,w) - \mathfrak{a}(z,x \times v,w) +
\mathfrak{a}(z,x,v \times w) =
\mathfrak{a}(z,x,v) \times w + z \times \mathfrak{a}(x,v,w)\notag
\end{equation}
but \(\mathfrak{a}(x,v,w) = 0\) and \(\mathfrak{a}(z,x,v) = 0\) by definition of
\(Z''\), and \(\mathfrak{a}(z,x,v \times w) = 0\) and
\(\mathfrak{a}(z,x \times v,w) = 0\) because \(W\) is a subring and by definition of
\(Z''\). Hence \(\mathfrak{a}(z \times x,v,w) = 0\). Likewise
\(\mathfrak{a}(z \times y,v,w) = 0\).

Using \thref{lem:stblsubrng} we conclude that \(Z'' = W\); also
\((\forall\thinspace u \in W)(\forall\thinspace v \in W)
(\forall\thinspace w \in W)\thickspace \mathfrak{a}(u,v,w) = 0\), which means
that ``\(\times\)'' is associative on \(W\).
\end{proof}

\pagebreak
\begin{proof}[Di-associativity of alternative rings]
Let \((A,+,\times)\) be an alternative ring. Quite obviously, if all subrings generated by
any two (not necessarily distinct) elements are associative, then all submagmas generated by
any two (not necessarily distinct) elements, being subsets of the subrings generated by the
same two (not necessarily distinct) elements, are associative as well, and we
already know that di-associative magmas are alternative. It therefore suffices to show that if
a ring is alternative, then all subrings generated by any two (not necessarily distinct) elements
are associative.

\thref{lem:monostrass} tells us that given
\(x \in A\) and \(y \in A\), if \(X\) is the subring generated by \(x\) and \(Y\) the
subring generated by \(Y\), then \(X\) and \(Y\) are strongly associative subrings.
We use \thref{lem:unionstrass} to conclude that the subring generated by \(\{x,y\}\)
is associative.

If now \((A,+,\times)\) is an alternative ring \emph{with unit}, the same reasoning as
above shows that it suffices to prove that all subrings of rings with units generated by two
(not necessarily distinct) elements are associative.

Given any two (not necessarily distinct) elements, \(x\) and \(y\), of \(A\), we know that
the subring generated by \(\{x, y\}\), which we will call \(S\), is associative. But if \(e\)
is the neutral element for ``\(\times\)'', then the subring of a ring with unit generated by
\(\{x, y\}\) is merely \(S \cup \bbZ\cdot e\) (where we have denoted by
\(\bbZ\cdot e\) the integer multiples, both positive and negative, of \(e\)), which is
clearly also associative.
\end{proof}

We immediately deduce from that the classical result that an algebra is alternative if and only
if it is di-assiociative, and if and only if all sub-algebras generated by any two (not
necessarily distinct) elements are associative, and that an algebra with unit is alternative if and
only if it is di-associative and if and only if all sub-algebras of algebras with unit, generated
by any two (not necessarily distinct) elements are associative.

We close this section with another property we will need
(\cite{GoodaireJespersMiles(1996)}).

\begin{proposition}\label{prp:altrngmouf}
Let \((A,+,\times)\) be an alternative ring with unit; then \((A,\times)\) verifies
\begin{gather}
(\forall\thinspace (x,y,z) \in A^{3})\thickspace \mathfrak{a}(x^{2},y,z) =
x\times\mathfrak{a}(x,y,z)-\mathfrak{a}(x,y,z)\times x\tag{R1}\label{eq:R1}\\
(\forall\thinspace (x,y,z) \in A^{3})\thickspace
\mathfrak{a}(x\times y,x,z) = \mathfrak{a}(x,y,z)\times x\tag{R2}\label{eq:R2}\\
(\forall\thinspace (x,y,z) \in A^{3})\thickspace
\mathfrak{a}(y\times x,x,z) = x\times \mathfrak{a}(x,y,z)\tag{R3}\label{eq:R3}
\end{gather}
and also verifies (\ref{eq:M1}), (\ref{eq:M2}), (\ref{eq:M3}), (\ref{eq:M4}) and
(\ref{eq:M5}).
\end{proposition}

\begin{proof}
The ring being alternative, \thref{prp:kleinalt} ensures that \(\mathfrak{k}(x,x,y,z) = 0\),
wich is (\ref{eq:R1}). Likewise \(\mathfrak{k}(x,y,x,z) = 0\) is (\ref{eq:R2}) and
\(\mathfrak{k}(y,x,x,z) = 0\) is (\ref{eq:R3}).

For (\ref{eq:M4}), we compute
\begin{multline*}
x\times[z\times(x\times y)] - [(x\times z)\times x]\times y =\\
\{x\times[z\times(x\times y)] - (x\times z)\times(x\times y)\}+
\{(x\times z)\times(x\times y) - [(x\times z)\times x]\times y\} =\\
\mathfrak{a}(x, z, x\times y)+\mathfrak{a}(x\times z, x, y)
\end{multline*}
Using (\ref{eqn:foursome}) we thus find that
\(x\times[z\times(x\times y)] - [(x\times z)\times x]\times y =
\mathfrak{a}(x,z\times x,y)+\mathfrak{a}(x,z,x)\times y+
x\times\mathfrak{a}(z,x,y)\), and since the ring is alternative,
\(\mathfrak{a}(x,z,x) = 0\).  Furthermore, since the ring is alternative,
\(\mathfrak{a}(x,z\times x,y) = -\mathfrak{a}(z\times x,x,y)\) and we use
(\ref{eq:R2}) (applied to \((x,z,y)\)) to find that
\(x\times[z\times(x\times y)] - [(x\times z)\times x]\times y = 0\), which is
(\ref{eq:M4}).

For (\ref{eq:M5}), we compute
\begin{multline*}
y\times[x\times(z\times x)] - [(y\times x)\times z]\times x =\\
\{y\times[x\times(z\times x)] - (y\times x)\times(z\times x)\}+
\{(y\times x)\times(z\times x) - [(y\times x)\times z]\times x\} =\\
\mathfrak{a}(y, x, z\times x)+\mathfrak{a}(y\times x, z, x)
\end{multline*}

For (\ref{eq:M1}), we compute
\begin{multline*}
(x\times y)\times(z\times x) - [x\times(y\times z)]\times x =\\
\{(x\times y)\times(z\times x) - [(x\times y)\times z]\times x\}+
\{[(x\times y)\times z]\times x - [x\times(y\times z)]\times x\} =\\
\mathfrak{a}(x\times y,z,x)-\mathfrak{a}(x, y, z)\times x
\end{multline*}

(\ref{eq:M2}) and (\ref{eq:M3}) are dealt with in a similar way.
\end{proof}

\subsubsection{Derivation graph for several laws}
\indent

We present here a graph representing the ``is-a'' relationship for some of the most important
structures in this second part of the appendix. We have used abbreviations akin to those of
\ref{ssc:derone} (\emph{i.e.} ``P.-A.'' stands for ``Power-Associative'', ``Alt.'' stands for
``Alternative'', ``Ass.'' stands for ``Associative'' and ``w. U.'' stands for ``with Unit'').
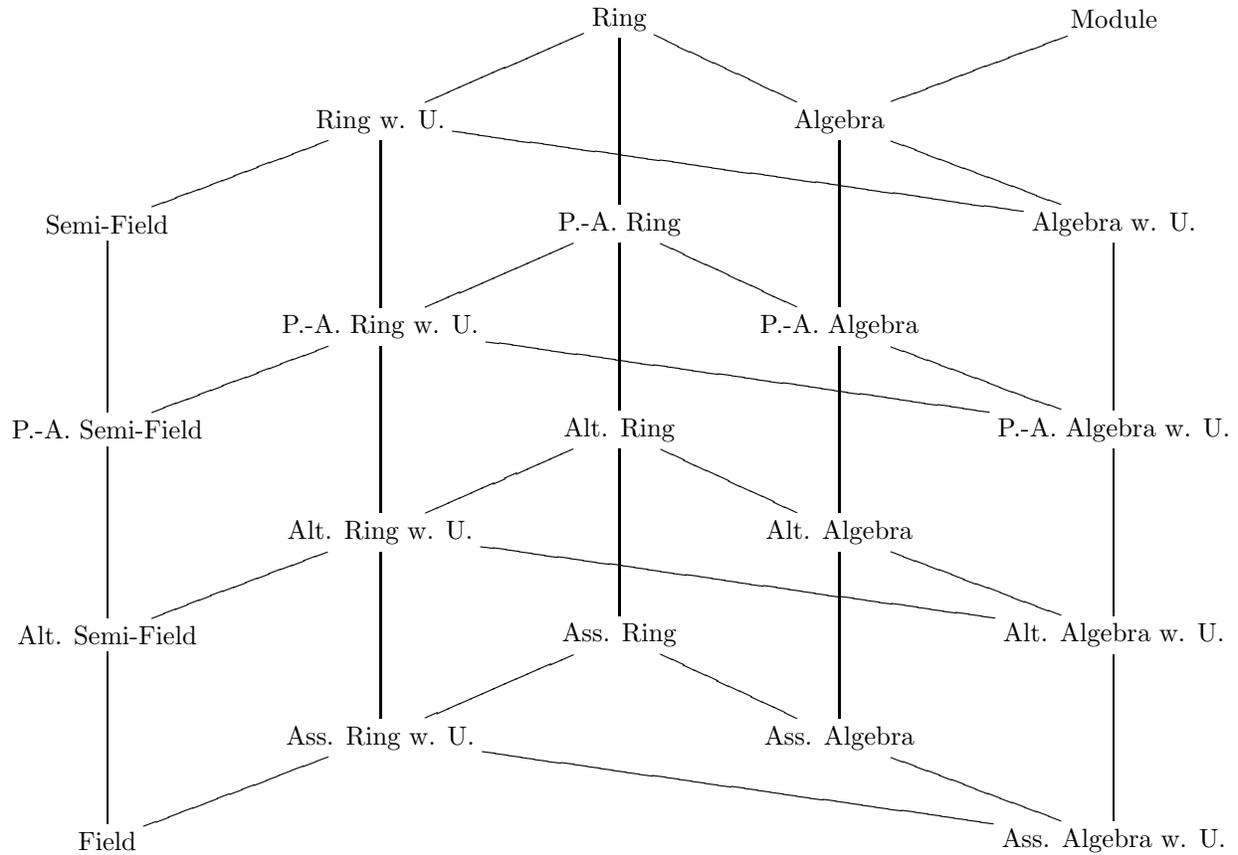
\begin{center}

\begin{math} 
\xymatrix{
 & & \text{Ring} \ar@{-}[dd] \ar@{-}[ld] \ar@{-}[rd] & & \text{Module} \ar@{-}[ld] \\
 & \text{Ring w. U.} \ar@{-}[dd] \ar@{-}[ld] \ar@{-}[rrrd] & & \text{Algebra} \ar@{-}[dd] \ar@{-}[rd] \\
\text{Semi-Field} \ar@{-}[dd] & & \text{P.-A. Ring} \ar@{-}[dd] \ar@{-}[ld] \ar@{-}[rd] & & \text{Algebra w. U.} \ar@{-}[dd] \\
 & \text{P.-A. Ring w. U.} \ar@{-}[dd] \ar@{-}[ld] \ar@{-}[rrrd] & & \text{P.-A. Algebra} \ar@{-}[dd] \ar@{-}[rd] \\
\text{P.-A. Semi-Field} \ar@{-}[dd] & & \text{Alt. Ring} \ar@{-}[dd] \ar@{-}[ld] \ar@{-}[rd] & & \text{P.-A. Algebra w. U.} \ar@{-}[dd]\\
 & \text{Alt. Ring w. U.} \ar@{-}[dd] \ar@{-}[ld] \ar@{-}[rrrd] & & \text{Alt. Algebra} \ar@{-}[dd] \ar@{-}[rd] \\
\text{Alt. Semi-Field} \ar@{-}[dd] & & \text{Ass. Ring} \ar@{-}[ld] \ar@{-}[rd] & & \text{Alt. Algebra w. U.} \ar@{-}[dd]\\
 & \text{Ass. Ring w. U.} \ar@{-}[ld] \ar@{-}[rrrd] & & \text{Ass. Algebra} \ar@{-}[rd] & \\
\text{Field} & & & & \text{Ass. Algebra w. U.} 
}
\end{math}

\end{center}

\pagebreak
\section{A brief reminder about Gauss and Jacobi sums}\label{apx:gjs}
\indent

We spell out here what is briefly outlined in \cite[page 147]{IrelandRosen(1982)}.

\subsection{Basic definitions and results}
\subsubsection{Multiplicative characters}
\indent

Let \(p\) be a prime number, \(k \geqslant 1\) an integer, \(q = p^k\) and
\(F = \GaloisField (q)\) the finite field of cardinal \(q\); as usual \(F^* = F-\{0\}\),
and we will denote by \(1_{F}\) the neutral element for the multiplication in \(F\). Recall
that \(F\) is of characteristic \(p\).

A multiplicative character on \(F\) is a function \(\chi : F^* \rightarrow \bbC^*\)
such that \((\forall\thinspace (a,b) \in F^{*2})\thickspace
\chi(a\cdot b) = \chi(a) \chi(b)\). The character
\(\epsilon_{F}:[F^* \rightarrow \bbC^* , a \mapsto 1]\) is extended to all of \(F\)
as \({\check{\epsilon}}_{F}(0) = 1\), and any character
\(\chi \neq \epsilon_{F}\) is extended to all of \(F\) as \(\check{\chi}(0) = 0\).
Therefore, \((\forall\thinspace (a,b) \in F^{2})\thickspace [
\check{\chi}(a\cdot b) = \check{\chi}(a) \check{\chi}(b) \text{ and }
{\check{\epsilon}}_{F}(a\cdot b) = {\check{\epsilon}}_{F}(a)
{\check{\epsilon}}_{F}(b)]\).

One then proves (as an analogue to \cite[proposition 8.1.1]{IrelandRosen(1982)})
that if \(\chi\) is a multiplicative character then \(\chi(1_{F}) = 1\),
and given \(a \in F^{*}\), \(\chi(a)\) is a \((q-1)^{\text{th}}\) root of unity and
\(\chi(a^{-1}) = {\chi(a)}^{-1} = \overline{\chi(a)}\).

Furthermore (mirroring \cite[proposition 8.1.2]{IrelandRosen(1982)}), given a
multiplicative character \(\chi\), we have:

\begin{center} 

\begin{math}
\displaystyle
\sum_{t \in F} \check{\chi}(t) =
\begin{cases}
q & \text{if \(\chi = \epsilon_{F}\)} \\
0 & \text{otherwise.}
\end{cases}
\end{math}

\end{center}

We also verify that, as in \cite[proposition 8.1.3]{IrelandRosen(1982)}, the group of
multiplicative characters is cyclic of cardinal \(q-1\), and that for any \(a \in F^{*}\)
there is a multiplicative character \(\chi\) such that \(\chi(a ) \neq 1\), with the same
corollary that if \(a \in F^{*}\) and \(a \neq 1_{F}\) then
\(\sum_{\chi \in \text{multiplicative characters of } F} \chi(a) = 0\).

\subsubsection{The \protect\(x^{n} = a\protect\) equation in \protect\(\GaloisField(q)\protect\)}\label{ssc:galeq}
\indent

As well, one sees that (just like in \cite[proposition 7.1.2]{IrelandRosen(1982)}), given
\(a \in F^*\), the equation \(x^n = a\) has solutions in \(F\) if and only if
\(a^{\frac{q-1}{d}} = 1_{F}\), with \(d = \PGCD(n, q-1)\) (\emph{i.e.} \(d\)
is the smallest common multiple of \(n\) and \(q-1\)), and if there are solutions, then there
are exactly \(d\) solutions.

From that, one deduces (as \cite[proposition 8.1.4]{IrelandRosen(1982)}) that given
\(a \in F^{*}\) and \(n\thinspace|\thinspace(q-1)\) (\emph{i.e.} a \(n\) which
divides \(q-1\)), if \(x^n = a\) is not solvable in \(F^{*}\), then there exists a
multiplicative character character \(\chi\) such that \(\chi^n = \epsilon_{F}\)
and \(\chi(a) \neq 1\).

This, of course, leads to the result (corresponding to
\cite[proposition 8.1.5]{IrelandRosen(1982)})
that the number of solutions of the equation \(x^n = a\)
in \(F\), which we will write \(N(x^n = a)\), is given by
\(N(x^n = a) = \sum_{\chi^n = \epsilon_{F}} \check{\chi}(a)\), if
\(n\thinspace|\thinspace(q-1)\).

\subsubsection{Galois trace}
\indent

Recall that given \(p\) a prime number and \(k \geqslant 1\)
an integer, then \(\GaloisField(p)\) is isomorphic, through some isomorphism
\(\rsfsJ_{p,k} : \GaloisField(p) \rightarrow\GaloisField(p^k)\),
to the subfield of \(\GaloisField(p^k)\) defined by \(\{\alpha \in \GaloisField(p^k)
\thickspace|\thickspace \alpha^{p} = \alpha\}\). Furthermore, given any
\(\alpha \in \GaloisField(p^k)\), if we compute
\(\beta = \alpha + \alpha^p + \alpha^{p^2}+ \cdots + \alpha^{p^{k-1}}\) then
we find that \(\beta^p = \beta\). Hence, given \(p\) and \(k\) as above,
we define\footnote{The unfortunate verbosity is to
ensure there is no confusion with the Cayley trace, as we will be using both notions in this
document.} the \emph{Galois trace} on \(\GaloisField(p^k)\) as the function

\begin{center}

\begin{math} 
\displaystyle
\GaloisTrace_{p^k}\negthinspace : \negthinspace
\begin{array}[t]{rcl}
\GaloisField(p^k) & \rightarrow & \GaloisField(p)\\
\alpha & \mapsto & \rsfsJ_{p,k}^{-1}(\alpha + \alpha^p + \alpha^{p^2}+ \cdots + \alpha^{p^{k-1}})
\end{array}
\end{math}

\end{center}

Recall that \(\GaloisField(p^k)\) is an algebra over \(\GaloisField(p)\), and that
\(\GaloisTrace_{p^k}\) is a surjective \(\GaloisField(p)\)-linear operator.

\pagebreak
\subsection{Gauss sums}
\indent

Let us define in this appendix, for convenience's sake,
\(\zeta_{p} = \NeperBase^{\frac{2 \pi i}{p}}\). Since, given
\(\bm{\alpha} \in {\bbZ / {p \bbZ}}\) and \(a \in \bm{\alpha}\) then
\((\forall\thinspace b \in \bm{\alpha})\thickspace
\zeta_{p}^{b} = \zeta_{p}^{a}\), we can unambiguously define
\(\zeta_{p}^{\bm{\alpha}}\) as the common value of \( \zeta_{p}^{a}\)
for all \(a \in \bm{\beta}\). We can therefor (identifying
\(\GaloisField(p)\) with \(\bbZ / {p \bbZ}\)) define the function

\begin{center}

\begin{math} 
\displaystyle
\psi\negthinspace : \negthinspace
\begin{array}[t]{rcl}
F & \rightarrow & \bbC\\
\alpha & \mapsto & \zeta_{p}^{\GaloisTrace(\alpha)}
\end{array}
\end{math}

\end{center}

The salient properties of \(\psi\) are that
\((\forall\thinspace (\alpha, \beta) \in F^2)\thickspace
\psi(\alpha+\beta) = \psi(\alpha)\psi(\beta)\),
\((\exists\thinspace \alpha_{0} \in F)\thickspace \psi(\alpha_{0}) \neq 1\),
and \(\sum_{\alpha \in F} \psi(\alpha) = 0\). Another crucial property is:

\begin{center} 

\begin{math}
\displaystyle
\sum_{\alpha \in F} \psi(\alpha\cdot\beta) =
\begin{cases}
q & \text{if \(\beta = 0\)} \\
0 & \text{otherwise.}
\end{cases}
\end{math}

\end{center}

We finally define
\(g_{\alpha}(\chi) = \sum_{\beta \in F}\check{\chi}(\beta)\psi(\alpha\cdot\beta)\)
to be a \emph{Gauss sum} belonging to the character \(\chi\)
(and relative to an \(\alpha \in F\)). We will also define
\(g(\chi) = g_{1_{F}}(\chi)\).

One can then prove (as for \cite[proposition 8.2.1]{IrelandRosen(1982)}) that
\(g_{\alpha}(\chi)\) takes the following values:
\begin{center}
\begin{tabular}{||c||c|c||}
\hhline{|t:=:t:==:t|}
\(g_{\alpha}(\chi)\) & \(\alpha = 0\) & \(\alpha \neq 0\) \\
\hhline{|:=::=|=:|}
\(\chi = \epsilon_{F}\) & q & 0 \\
\hhline{||---||}
\(\chi \neq \epsilon_{F}\) & 0 & \(\chi(\alpha^{-1})\thinspace g(\chi)\) \\
\hhline{|b:=:b:==:b|}
\end{tabular}
\end{center}

As well (mirroring \cite[proposition 8.2.2]{IrelandRosen(1982)}),
\(g(\chi^{-1}) = g(\bar{\chi}) = \chi(-1_{F}) \overline{g(\chi)}\)
and if furthermore \(\chi \neq \epsilon_{F}\) then
\(g(\chi)\thinspace g(\chi^{-1}) = \chi(-1_{F})\thinspace q\), which together imply
that \(|g(\chi)| = \sqrt{q}\).

\subsection{Jacobi sums}\label{sbs:jacosum}
\indent

Given \(\ell\) (at least two) multiplicative characters on \(F = \GaloisField(q)\),
\(\chi_{1}\), \(\chi_{2}\), \ldots, \(\chi_{\ell}\), the \emph{Jacobi sums} of these
characters are defined to be:
\begin{align*}
J(\chi_{1}, \chi_{2},\ldots,\chi_{\ell}) &= \sum_{t_{1}+t_{2}\cdots+t_{\ell} = 1_{F}} \check{\chi}_{1}(t_{1})\check{\chi}_{2}(t_{2})\ldots\check{\chi}_{\ell}(t_{\ell}) \\
J_{0}(\chi_{1}, \chi_{2},\ldots,\chi_{\ell}) &= \sum_{t_{1}+t_{2}\cdots+t_{\ell} = 0} \check{\chi}_{1}(t_{1})\check{\chi}_{2}(t_{2})\ldots\check{\chi}_{\ell}(t_{\ell})
\end{align*}

Additionally we will define \(J(\chi) = 1_{F}\) for any multiplicative character \(\chi\).

The Jacobi sums verify the following properties (rewriting of
\cite[proposition~8.5.1, theorem~3, corollary~1 and corollary~2]{IrelandRosen(1982)}):
\begin{description}
\item[] if \((\chi_{1},\ldots,\chi_{\ell-1}) = (\epsilon_{F},\ldots,\epsilon_{F})\) then
    \begin{description}
    \item[] \(J(\epsilon_{F}, \epsilon_{F},\ldots,\epsilon_{F}) = J_{0}(\epsilon_{F}, \epsilon_{F},\ldots,\epsilon_{F}) = q^{\ell-1}\)
    \end{description}
\item[] else (\emph{i.e.}: \((\chi_{1},\ldots,\chi_{\ell-1}) \neq (\epsilon_{F},\ldots,\epsilon_{F})\))
    \begin{description}
    \item[] if \(\chi_{1} = \epsilon_{F}\) or \(\chi_{2} = \epsilon_{F}\) or \ldots or \(\chi_{\ell} = \epsilon_{F}\) then
        \begin{description}
        \item[] \(J(\epsilon_{F}, \epsilon_{F},\ldots,\epsilon_{F}) = J_{0}(\epsilon_{F}, \epsilon_{F},\ldots,\epsilon_{F}) = 0\)
        \end{description}
    \item[] else (\emph{i.e.}: \(\chi_{1} \neq \epsilon_{F}\) and \(\chi_{2} \neq \epsilon_{F}\) and \ldots and \(\chi_{\ell} \neq \epsilon_{F}\))
        \begin{description}
        \item[] if \(\chi_{1}\chi_{2}\cdots\chi_{\ell} = \epsilon_{F}\) then
            \begin{description}
            \item[] \( J_{0}(\epsilon_{F}, \epsilon_{F},\ldots,\epsilon_{F}) = (q-1)\thinspace\chi_{\ell}(-1_{F})\thinspace J(\chi_{1},\chi_{2},\ldots,\chi_{\ell-1})\)
            \item[] \(J(\chi_{1},\chi_{2},\ldots,\chi_{\ell}) = -\chi_{\ell}(-1_{F})\thinspace J(\chi_{1},\chi_{2},\ldots,\chi_{\ell-1})\)
            \item[] \(g(\chi_{1})g(\chi_{2}) \cdots g(\chi_{\ell}) = q\thinspace\chi_{\ell}(-1_{F})\thinspace J(\chi_{1},\chi_{2},\ldots,\chi_{\ell-1})\)
            \end{description}
        \item[] else
            \begin{description}
            \item[] \( J_{0}(\epsilon_{F}, \epsilon_{F},\ldots,\epsilon_{F}) = 0\)
            \item[] \(g(\chi_{1})g(\chi_{2}) \cdots g(\chi_{\ell}) = g(\chi_{1}\chi_{2}\cdots\chi_{\ell})\thinspace J(\chi_{1},\chi_{2},\ldots,\chi_{\ell})\)
            \end{description}
        \end{description}
    \end{description}
\end{description}

\pagebreak
\section{Cayley algebra \emph{versus} Quadratic algebra and others}\label{sec:fourtout}
\indent

There are unfortunately several rather different meanings to the term ``quadratic algebra''.

One of these (\cite{Schafer(1966), ZhevlakovSlin'koShestakovShirshov(1982),
GoodaireJespersMiles(1996)}) is to say that every
\(x\) in the algebra with unit \(E\) (with neutral element \(e\)) over the commutative and
associative ring with unit \(A\), satisfies an equation of the form \(x^{2}-t(x)x+n(x)e = 0\)
for some \(t(x) \in A\) and \(n(x) \in A\). For this meaning of the term ``quadratic
algebra'', all Cayley algebras are quadratic algebras.

Another meaning to the term  ``quadratic algebra'' is found in \cite{Bourbaki(A)}.
For this other meaning of ``quadratic algebra'' however, the main thing to keep in mind
is that these two kinds of structure, while related, and sometimes overlapping (some algebras
are both quadratic algebras and Cayley algebras), \emph{are different things}.
We will highlight here some of the similarities and difference between them.

Let \((A,+,.)\) be an associative and commutative ring with unit, whose neutral
element we will denote by \(1_{A}\). An \(A\)-algebra \((E,+,\times,\cdot)\) is a
\emph{quadratic algebra} over \(A\) if and only if it admits, as an \(A\)-module, a
generating\footnote{We do not impose the family to be linearly independent over \(A\),
so as to accept perhaps pathological cases such as \(A = \bbZ/2\bbZ\), which can be seen
as a quadratic algebra of type \((\bar{1}, \bar{0})\) over itself, using the family
 \((e_{1},e_{2}) = (\bar{1}, \bar{1})\).} family \((e_{1},e_{2})\) such that the
 multiplication (``\(\times\)'') verifies
\begin{gather*}
e_{1}^{2} = e_{1} \\
e_{1}\times e_{2} = e_{2}\times e_{1} = e_{2} \\
e_{2}^{2} = \alpha\cdot e_{1}+\beta\cdot e_{2}
\end{gather*}
for some \((\alpha,\beta) \in A^{2}\), and we then say the quadratic algebra is 
\emph{of type} \((\alpha,\beta)\), though this is not an intrinsic property of the algebra,
\emph{i.e.} the algebra can happen\footnote{It should also be born in mind that, where
Cayley algebras are concerned, applying the Cayley-Dickson procedure with different
constants may happen to yield isomorphic structures.} to also be of type 
\((\alpha',\beta')\) with
 \((\alpha',\beta') \neq (\alpha,\beta)\) (using another family).

As such, quadratic algebras are necessarily commutative, associative and with unit (with
neutral element \(e_{1}\)), and this is of course not true of every Cayley algebra (take
the classical quaternions or octonions for instance) which means there are structures which
are Cayley algebra but are not quadratic algebras (as we give examples of below,). As well,
this also means any quadratic \(A\)-algebra can be made into a Cayley algebra,
over itself if perhaps not over \(A\), by using the identity as conjugation.

It is worth noting that if one makes a Cayley algebra out of an associative and commutative
ring with unit \(A\), considering it as an algebra over itself and using the identity as
conjugation, then applying the Cayley-Dickson process with a constant \(\zeta\), then one
gets a structure which is both a Cayley algebra \emph{and} a quadratic algebra of type
\((-\zeta,0)\), when considering the basis \(((1_{A},0),(0,1_{A}))\), of \(A^{2}\) over
\(A\). Of course, then the conjugation over \(A^{2}\) is usually no longer the
identity\footnote{Though it may still be, if for instance \(A\) is of characteristic \(2\)!}
and reiterating the Cayley-Dickson process over this new structure yields another Cayley
algebra which is usually \emph{not} a quadratic algebra over either \(A\) or \(A^{2}\)
(this is the case with \(A = \bbR\) yielding first \(\bbC\) and then \(\bbH\), for instance).

Like Cayley algebras, quadratic algebras have a ``natural'' doubling process, which, of course,
produces structures which always are quadratic algebras, but are usually not Cayley algebras.
If \((A,+,.)\) is an associative and commutative ring with unit (with neutral
element \(1_{A}\)), \((E,+,\times,\cdot)\) is a quadratic \(A\)-algebra (with neutral
element \(1_{E}\)), and
\(\alpha\) and \(\beta\) any two elements of \(E\), then on
\(F = E \vartimes E\) we consider the following laws:
\begin{gather*}
\begin{aligned}
\dagger\negthinspace : & &  F \vartimes F & & \rightarrow & & & F \\
 & & ((x,y),(x',y')) & &\mapsto & & & (x+x',y+y') \\
 *\negthinspace : & &  F \vartimes F & & \rightarrow & & & F \\
 & & ((x,y),(x',y')) & & \mapsto & & & (x\times x' + \alpha \times (y'\times y),y \times x' + y'\times x+\beta\times(y \times y')) \\
  \bullet\negthinspace : & &  E \vartimes F & & \rightarrow & & & F \\
  & & (\lambda,(x,y)) & & \mapsto & & & (\lambda \times x, \lambda \times y)
\end{aligned}
\end{gather*}
it is then trivial to check that \((F,\dagger,*,\bullet)\) is a quadratic \(E\)-algebra, of
type \((\alpha,\beta)\) when considering the basis \(((1_{E},0),(0,1_{E}))\) of \(F\)
over \(E\). Indeed, when one starts with \(E = A\), apply first the above doubling using
constants \(\alpha\) and \(\beta\) in \(E = A\) and \emph{then} apply the
Cayley-Dickson process using a constant \(\zeta \in E = A\) and the identity on
\(E = A\) as conjugation, the result is known in
literature (\cite{Bourbaki(A)}) as a quaternion algebra of type \((\alpha,\beta,-\zeta)\).
This, of course, is a Cayley algebra, but usually not a quadratic algebra (either over
\(E = A\), or \(E^{2}\)). One can therefore usually not indefinitely ``knit'' the two kinds
of doubling, which makes this construction rather \emph{ad'hoc}.

There are also other known ways to build useful algebras in a way reminiscent of the
Cayley-Dickson process. One of the better known (\cite{GoodaireJespersMiles(1996)})
is to consider a commutative and associative ring with unit \(A\), with neutral element
\(1_{A}\), and build the algebra
\(H(A) = A + A\thinspace I + A\thinspace J + A\thinspace K\) with
\(I^{2} = J^{2} = K^{2} = IJK = -1_{A}\). With \(A = \bbR\) one obtains \(\bbH\),
the usual quaternions, and with \(A = \bbC\) one obtains what is known as the
\emph{complex quaternions} (or \emph{biquaternions}) (\cite{Imaeda(1986),
Kravchenko(2003)}).

There is at least one more construction which deserves consideration, and which does not seem
to fit well with those we have considered so far: that of the sedenions.

In \cite{Carmody(1988), Carmody(1997)} Carmody presents constructions originally
due to Mus{\`e}s, some of which (the counter-complex, counter-quaternions and
counter-octonions) are easily seen as examples of the Cayley-Dickson
process in action. Indeed, we can draw the following filliation graph (where \(\bbR\),
\(\bbC\), \(\bbH\) and \(\bbO\) stand for the set of real numbers, the set of complex
numbers, the quaternions and the octonions respectively; \(\bbX\) we have dubbed the
\emph{hexadecimalions}, \cite{Holin(1999)}, but see below):
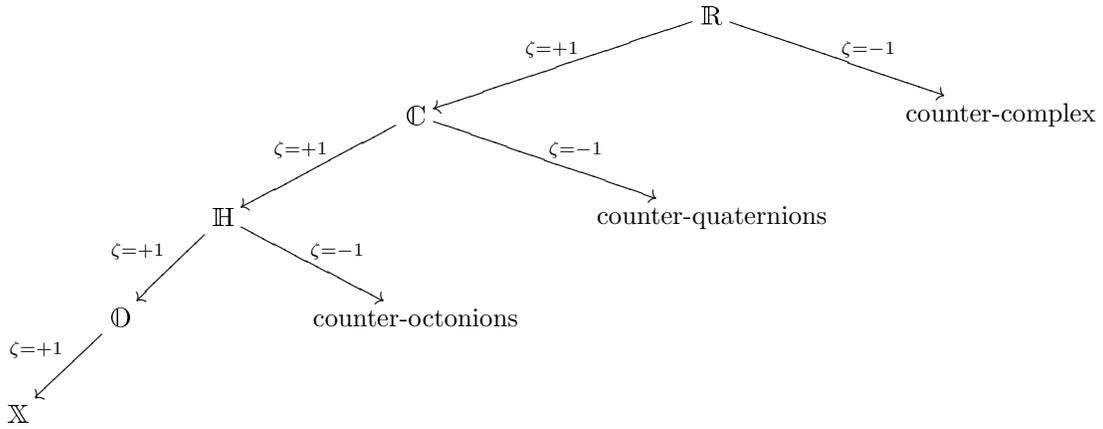
\begin{center}
\begin{math} 
\xymatrix{
&&&&\bbR\ar@{->}[ld]_{\zeta = +1}\ar@{->}[rd]^{\zeta = -1}&\\
&&&\bbC\ar@{->}[ld]_{\zeta = +1}\ar@{->}[rd]^{\zeta = -1}&&\text{counter-complex}\\
&&\bbH\ar@{->}[ld]_{\zeta = +1}\ar@{->}[rd]^{\zeta = -1}&&\text{counter-quaternions}&\\
&\bbO\ar@{->}[ld]_{\zeta = +1}&&\text{counter-octonions}&&\\
\bbX&&&&&
}
\end{math}
\end{center}

There is also a structure on \(\bbR^{16}\), which has been called by  Mus{\`e}s and
Carmody, the \emph{sedenion}. That structure is different from that of \(\bbX\), as we shall
see, but unfortunately the name ``sedenion'' has also (posteriorly) been used to designate the
hexadecimalions (\cite{ImaedaImaeda(2000)} and, at the time of this writing,
\cite{SedenionWiki}).

We will keep the name ``sedenion'' to
designate the construct of Mus{\`e}s and Carmody. It has the following Cayley
(multiplication) table:
\begin{center}
\begin{tabular}{c|*{16}{c}}
&\(1\)&\(i_{1}\)&\(i_{2}\)&\(i_{3}\)&\(i_{4}\)&\(i_{5}\)&\(i_{6}\)&\(i_{7}\)&\(e_{1}\)&\(e_{2}\)&\(e_{3}\)&\(e_{4}\)&\(e_{5}\)&\(e_{6}\)&\(e_{7}\)&\(\omega\)\\
\hhline{-*{16}{-}}
\(1\)&\(+1\)&\(+i_{1}\)&\(+i_{2}\)&\(+i_{3}\)&\(+i_{4}\)&\(+i_{5}\)&\(+i_{6}\)&\(+i_{7}\)&\(+e_{1}\)&\(+e_{2}\)&\(+e_{3}\)&\(+e_{4}\)&\(+e_{5}\)&\(+e_{6}\)&\(+e_{7}\)&\(+\omega\)\\
\(i_{1}\)&\(+i_{1}\)&\(-1\)&\(+i_{3}\)&\(-i_{2}\)&\(+i_{5}\)&\(-i_{4}\)&\(-i_{7}\)&\(+i_{6}\)&\(+\omega\)&\(+e_{3}\)&\(-e_{2}\)&\(+e_{5}\)&\(-e_{4}\)&\(-e_{7}\)&\(+e_{6}\)&\(-e_{1}\)\\
\(i_{2}\)&\(+i_{2}\)&\(-i_{3}\)&\(-1\)&\(+i_{1}\)&\(+i_{6}\)&\(+i_{7}\)&\(-i_{4}\)&\(-i_{5}\)&\(-e_{3}\)&\(+\omega\)&\(+e_{1}\)&\(+e_{6}\)&\(+e_{7}\)&\(-e_{4}\)&\(-e_{5}\)&\(-e_{2}\)\\
\(i_{3}\)&\(+i_{3}\)&\(+i_{2}\)&\(-i_{1}\)&\(-1\)&\(+i_{7}\)&\(-i_{6}\)&\(+i_{5}\)&\(-i_{4}\)&\(+e_{2}\)&\(-e_{1}\)&\(+\omega\)&\(+e_{7}\)&\(-e_{6}\)&\(+e_{5}\)&\(-e_{4}\)&\(-e_{3}\)\\
\(i_{4}\)&\(+i_{4}\)&\(-i_{5}\)&\(-i_{6}\)&\(-i_{7}\)&\(-1\)&\(+i_{1}\)&\(+i_{2}\)&\(+i_{3}\)&\(-e_{5}\)&\(-e_{6}\)&\(-e_{7}\)&\(+\omega\)&\(+e_{1}\)&\(+e_{2}\)&\(+e_{3}\)&\(-e_{4}\)\\
\(i_{5}\)&\(+i_{5}\)&\(+i_{4}\)&\(-i_{7}\)&\(+i_{6}\)&\(-i_{1}\)&\(-1\)&\(-i_{3}\)&\(+i_{2}\)&\(+e_{4}\)&\(-e_{7}\)&\(+e_{6}\)&\(-e_{1}\)&\(+\omega\)&\(-e_{3}\)&\(+e_{2}\)&\(-e_{5}\)\\
\(i_{6}\)&\(+i_{6}\)&\(+i_{7}\)&\(+i_{4}\)&\(-i_{5}\)&\(-i_{2}\)&\(+i_{3}\)&\(-1\)&\(-i_{1}\)&\(+e_{7}\)&\(+e_{4}\)&\(-e_{5}\)&\(-e_{2}\)&\(+e_{3}\)&\(+\omega\)&\(-e_{1}\)&\(-e_{6}\)\\
\(i_{7}\)&\(+i_{7}\)&\(-i_{6}\)&\(+i_{5}\)&\(+i_{4}\)&\(-i_{3}\)&\(-i_{2}\)&\(+i_{1}\)&\(-1\)&\(-e_{6}\)&\(+e_{5}\)&\(+e_{4}\)&\(-e_{3}\)&\(-e_{2}\)&\(+e_{1}\)&\(+\omega\)&\(-e_{7}\)\\
\(e_{1}\)&\(+e_{1}\)&\(+\omega\)&\(+e_{3}\)&\(-e_{2}\)&\(+e_{5}\)&\(-e_{4}\)&\(-e_{7}\)&\(+e_{6}\)&\(+1\)&\(-i_{3}\)&\(+i_{2}\)&\(-i_{5}\)&\(+i_{4}\)&\(+i_{7}\)&\(-i_{6}\)&\(+i_{1}\)\\
\(e_{2}\)&\(+e_{2}\)&\(-e_{3}\)&\(+\omega\)&\(+e_{1}\)&\(+e_{6}\)&\(+e_{7}\)&\(-e_{4}\)&\(-e_{5}\)&\(+i_{3}\)&\(+1\)&\(-i_{1}\)&\(-i_{6}\)&\(-i_{7}\)&\(+i_{4}\)&\(+i_{5}\)&\(+i_{2}\)\\
\(e_{3}\)&\(+e_{3}\)&\(+e_{2}\)&\(-e_{1}\)&\(+\omega\)&\(+e_{7}\)&\(-e_{6}\)&\(+e_{5}\)&\(-e_{4}\)&\(-i_{2}\)&\(+i_{1}\)&\(+1\)&\(-i_{7}\)&\(+i_{6}\)&\(-i_{5}\)&\(+i_{4}\)&\(+i_{3}\)\\
\(e_{4}\)&\(+e_{4}\)&\(-e_{5}\)&\(-e_{6}\)&\(-e_{7}\)&\(+\omega\)&\(+e_{1}\)&\(+e_{2}\)&\(+e_{3}\)&\(+i_{5}\)&\(+i_{6}\)&\(+i_{7}\)&\(+1\)&\(-i_{1}\)&\(-i_{2}\)&\(-i_{3}\)&\(+i_{4}\)\\
\(e_{5}\)&\(+e_{5}\)&\(+e_{4}\)&\(-e_{7}\)&\(+e_{6}\)&\(-e_{1}\)&\(+\omega\)&\(-e_{3}\)&\(+e_{2}\)&\(-i_{4}\)&\(+i_{7}\)&\(-i_{6}\)&\(+i_{1}\)&\(+1\)&\(+i_{3}\)&\(-i_{2}\)&\(+i_{5}\)\\
\(e_{6}\)&\(+e_{6}\)&\(+e_{7}\)&\(+e_{4}\)&\(-e_{5}\)&\(-e_{2}\)&\(+e_{3}\)&\(+\omega\)&\(-e_{1}\)&\(-i_{7}\)&\(-i_{4}\)&\(+i_{5}\)&\(+i_{2}\)&\(-i_{3}\)&\(+1\)&\(+i_{1}\)&\(+i_{6}\)\\
\(e_{7}\)&\(+e_{7}\)&\(-e_{6}\)&\(+e_{5}\)&\(+e_{4}\)&\(-e_{3}\)&\(-e_{2}\)&\(+e_{1}\)&\(+\omega\)&\(+i_{6}\)&\(-i_{5}\)&\(-i_{4}\)&\(+i_{3}\)&\(+i_{2}\)&\(-i_{1}\)&\(+1\)&\(+i_{7}\)\\
\(\omega\)&\(+\omega\)&\(-e_{1}\)&\(-e_{2}\)&\(-e_{3}\)&\(-e_{4}\)&\(-e_{5}\)&\(-e_{6}\)&\(-e_{7}\)&\(+i_{1}\)&\(+i_{2}\)&\(+i_{3}\)&\(+i_{4}\)&\(+i_{5}\)&\(+i_{6}\)&\(+i_{7}\)&\(-1\)\\
\end{tabular}
\end{center}

\pagebreak
There is an involutive automorphism of the sedenions defined by
\begin{equation*}
\sigma(x+\alpha_{1}\thinspace i_{1}+\cdots+\alpha_{7}\thinspace i_{7}+
\beta_{1}\thinspace e_{1}+\cdots+\beta_{7}\thinspace e_{7}+t\thinspace\omega)
= x-\alpha_{1}\thinspace i_{1}-\cdots-\alpha_{7}\thinspace i_{7}
-\beta_{1}\thinspace e_{1}-\cdots-\beta_{7}\thinspace e_{7}
+t\thinspace\omega
\end{equation*}
which \emph{does not} make the sedenions into a Cayley algebra (because the image of
\([x \mapsto x + \sigma(x)]\) is not contained in \(\bbR\)), but such that if we define
\(\rsfsN\) by \(\rsfsN(x) = x\thinspace\sigma(x)\), then
\(\rsfsN(x\thinspace y) = \rsfsN(x)\thinspace\rsfsN(y)\).
This, actually, is what seems to have driven the creation of the sedenions.

There are varying expectations on how a ``norm'' defined on an algebra should behave. For
norms taking their
values in a valuated field, the least one would expect would be that the norm of the product of
two elements to be \emph{smaller or equal to} the product of the norms of the elements (the
usual rule for square matrices of reals, seen as linear operators). A
more stringent requirement is the ``composition'' behavior, by which the norm of the
product of two elements is \emph{equal to} the product of the norms of the elements (the
usual rule for reals, complex, quaternions and octonions).

Without any assumption on associativity for the multiplication of that algebra, either behavior
can fail to be met.
For instance, consider in \(\bbX\),  \(x = i+j''\) and \(y = e'+k'''\), then
\(x \times y =  2j'-2j'''\) and \(\CayleyNorm{\bbX}(x \times y) = 8 > 4 =
\CayleyNorm{\bbX}(x) \times \CayleyNorm{\bbX}(y)\), where we have considered the
Cayley norm as taking its values in \(\bbR\), as explained earlier. Also, if we consider
the following example given by Carmody, \(z = i+e'e''=i+e'''\) and \(t = j+k'e''=j+k'''\),
then \(z \times t = 0\) and \(\CayleyNorm{\bbX}(z \times t) = 0 < 4 =
\CayleyNorm{\bbX}(z) \times \CayleyNorm{\bbX}(t)\). There is therefore really no
reasonable relation to be expected in \(\bbX\) between the norm of a product and the
product of the norms!

The sedenions, as we have said, display the ``composition'' behavior; they are not
a composition algebra (\cite{ZhevlakovSlin'koShestakovShirshov(1982), 
SpringerVeldkamp(2000)}) however, because the norm is not
non-degenerate.
According to \cite{Carmody(1988)}, the algebra is alternative, and verifies
(\ref{eq:M5}), (\ref{eq:M4}), (\ref{eq:M3}) and (\ref{eq:M2}). It is also riddled with
zero divisors (as it contains the counter-complex, counter-quaternions and counter-octonions;
it also contains the octonions).

While the construction of the sedenions could be carried out on other rings, it does appear
like a stunning \emph{ad'hoc} construction.


\pagebreak
\bibliography{Cayley_Integers}

\bibliographystyle{hplain}

 \end{document}